\documentclass[preprint,3p]{elsarticle}

\usepackage{graphicx}
\usepackage{amsmath}
\usepackage{array}
\usepackage{amssymb}
\usepackage{multirow}
\usepackage{subfig}
\usepackage{color}

\usepackage{amssymb}

\newcommand{\Cred}{black}
\newcommand{\Cblue}{black}

\journal{Elsevier}

\begin{document}

\begin{frontmatter}



\title{Stable evaluation of Green's functions in cylindrically stratified regions with uniaxial anisotropic layers}


\author[rvt1]{H. Moon\corref{cor1}}
\ead{haksu.moon@gmail.com}

\author[rvt2]{B. Donderici}
\ead{burkay.donderici@halliburton.com}

\author[rvt3]{F. L. Teixeira}
\ead{teixeira@ece.osu.edu}

\cortext[cor1]{Corresponding author}

\address[rvt1]{ElectroScience Laboratory, The Ohio State University, Columbus, OH 43212, USA (present address: Intel Corporation, Hillsboro, OR 97124, USA)}

\address[rvt2]{Sensor Physics \& Technology, Halliburton Energy Services, Houston, TX 77032, USA}

\address[rvt3]{ElectroScience Laboratory, The Ohio State University, Columbus, OH 43212, USA}

\begin{abstract}
We present a robust algorithm for the computation of electromagnetic fields radiated by point sources (Hertzian dipoles) in cylindrically stratified media where each layer may exhibit material properties (permittivity, permeability, and conductivity) with uniaxial anisotropy. Analytical expressions are obtained based on the spectral representation of the tensor Green's function based on cylindrical Bessel and Hankel eigenfunctions, and extended for layered uniaxial media. Due to the poor scaling of these eigenfunctions for extreme arguments and/or orders, direct numerical evaluation of such expressions can produce numerical instability, i.e., underflow, overflow, and/or round-off errors under finite precision arithmetic. To circumvent these problems, we develop a numerically stable formulation through suitable rescaling of various expressions involved in the computational chain, to yield a robust algorithm for all parameter ranges. Numerical results are presented to illustrate the robustness of the formulation including cases of practical interest to geophysical exploration.
\end{abstract}


\end{frontmatter}


%
\newcommand{\Jn}{J_n}
\newcommand{\Jnd}{J'_n}
\newcommand{\Hn}{H^{(1)}_n}
\newcommand{\Hnd}{H'^{(1)}_n}
\newcommand{\hJn}{\hat{J}_n}
\newcommand{\hJnd}{\hat{J}'_n}
\newcommand{\hHn}{\hat{H}^{(1)}_n}
\newcommand{\hHnd}{\hat{H}'^{(1)}_n}
%
\newcommand{\bJ}{\overline{\mathbf{J}}}
\newcommand{\bJn}{\overline{\mathbf{J}}\,_n}
\newcommand{\bJzn}{\overline{\mathbf{J}}\,_{zn}}
\newcommand{\bJpn}{\overline{\mathbf{J}}\,_{\phi n}}
\newcommand{\hbJ}{\hat{\overline{\mathbf{J}}}}
\newcommand{\hbJn}{\hat{\overline{\mathbf{J}}}\,_n}
\newcommand{\hbJzn}{\hat{\overline{\mathbf{J}}}\,_{zn}}
\newcommand{\hbJpn}{\hat{\overline{\mathbf{J}}}\,_{\phi n}}
\newcommand{\bH}{\overline{\mathbf{H}}}
\newcommand{\bHn}{\overline{\mathbf{H}}\,^{(1)}_n}
\newcommand{\bHzn}{\overline{\mathbf{H}}\,^{(1)}_{zn}}
\newcommand{\bHpn}{\overline{\mathbf{H}}\,^{(1)}_{\phi n}}
\newcommand{\hbH}{\hat{\overline{\mathbf{H}}}}
\newcommand{\hbHn}{\hat{\overline{\mathbf{H}}}\,^{(1)}_n}
\newcommand{\hbHzn}{\hat{\overline{\mathbf{H}}}\,^{(1)}_{zn}}
\newcommand{\hbHpn}{\hat{\overline{\mathbf{H}}}\,^{(1)}_{\phi n}}
\newcommand{\bD}{\overline{\mathbf{D}}}
\newcommand{\hbD}{\hat{\overline{\mathbf{D}}}}
\newcommand{\bR}{\overline{\mathbf{R}}\,}
\newcommand{\tbR}{\widetilde{\overline{\mathbf{R}}}\,}
\newcommand{\hbR}{\hat{\overline{\mathbf{R}}}\,}
\newcommand{\htbR}{\hat{\widetilde{\overline{\mathbf{R}}}}\,}
\newcommand{\bRn}{\overline{\mathbf{R}}_n}
\newcommand{\bT}{\overline{\mathbf{T}}\,}
\newcommand{\tbT}{\widetilde{\overline{\mathbf{T}}}\,}
\newcommand{\hbT}{\hat{\overline{\mathbf{T}}}\,}
\newcommand{\htbT}{\hat{\widetilde{\overline{\mathbf{T}}}}\,}
\newcommand{\bS}{\overline{\mathbf{S}}\,}
\newcommand{\hbS}{\hat{\overline{\mathbf{S}}}\,}
\newcommand{\bI}{\overline{\mathbf{I}}}
\newcommand{\bM}{\overline{\mathbf{M}}\,}
\newcommand{\tbM}{\widetilde{\overline{\mathbf{M}}}\,}
\newcommand{\hbM}{\hat{\overline{\mathbf{M}}}\,}
\newcommand{\htbM}{\hat{\widetilde{\overline{\mathbf{M}}}\,}}
\newcommand{\bMn}{\overline{\mathbf{M}}_n}
\newcommand{\bN}{\overline{\mathbf{N}}\,}
\newcommand{\tbN}{\widetilde{\overline{\mathbf{N}}}\,}
\newcommand{\hbN}{\hat{\overline{\mathbf{N}}}\,}
\newcommand{\bF}{\overline{\mathbf{F}}}
\newcommand{\bW}{\overline{\mathbf{W}}}
\newcommand{\bX}{\overline{\mathbf{X}}}
\newcommand{\obX}{\ddot{\overline{\mathbf{X}}}}
\newcommand{\bY}{\overline{\mathbf{Y}}}
\newcommand{\obY}{\ddot{\overline{\mathbf{Y}}}}
\newcommand{\bZ}{\overline{\mathbf{Z}}}
\newcommand{\obZ}{\ddot{\overline{\mathbf{Z}}}}
\newcommand{\bBn}{\overline{\mathbf{B}}_n}
\newcommand{\bCn}{\overline{\mathbf{C}}_n}
\newcommand{\bLn}{\overline{\mathbf{L}}_n}
\newcommand{\ba}{\overline{\boldsymbol{\alpha}}}
\newcommand{\As}[1]{\overline{\boldsymbol{\alpha}}_{#1}}      
\newcommand{\Ass}[2]{\overline{\boldsymbol{\alpha}}_{#1}^{#2}} 
\newcommand{\bb}{\overline{\boldsymbol{\beta}}}
\newcommand{\Bs}[1]{\overline{\boldsymbol{\beta}}_{#1}}        
\newcommand{\Bss}[2]{\overline{\boldsymbol{\beta}}_{#1}^{#2}}  

\newcommand{\btt}[1]{\textcolor{red}{d_{#1}}}

\newcommand{\rr}{\mathbf{r}}
\newcommand{\rp}{\mathbf{r'}}
\newcommand{\suma}{\sum_{n=-\infty}^{\infty}}
\newcommand{\sumb}{\sum_{n=1}^{\infty}}
\newcommand{\intmp}{\int_{-\infty}^{\infty}}
\newcommand{\Fn}{\overline{\mathbf{F}}_n(\rho,\rho')}
\newcommand{\Dj}{\overleftarrow{\mathbf{D}}'_{j}}
\newcommand{\Dja}{\overleftarrow{\mathbf{D}}'_{j1}}
\newcommand{\Djb}{\overleftarrow{\mathbf{D}}'_{j2}}
\newcommand{\Djc}{\overleftarrow{\mathbf{D}}'_{j3}}

\newcommand{\pa}{\partial}
\newcommand{\iu}{\mathrm{i}}
\newcommand{\na}{\boldsymbol{\nabla}}
\newcommand{\tna}{\boldsymbol{\widetilde{\nabla}}}

\newcommand{\trr}{\widetilde{\mathbf{r}}}
\newcommand{\tr}{\widetilde{r}}
\newcommand{\trp}{\widetilde{\mathbf{r}}'}
\newcommand{\tk}{\widetilde{k}}
\newcommand{\dk}{\ddot{k}}
\newcommand{\tE}{\widetilde{\mathbf{E}}}
\newcommand{\tH}{\widetilde{\mathbf{H}}}
\newcommand{\tJ}{\widetilde{\mathbf{J}}}
\newcommand{\teps}{\widetilde{\epsilon}}
\newcommand{\tmu}{\widetilde{\mu}}
\newcommand{\ue}{\overline{\overline{\epsilon}}}
\newcommand{\us}{\overline{\overline{\sigma}}}
\newcommand{\um}{\overline{\overline{\mu}}}
\newcommand{\ttGa}{\overline{\overline{\Gamma}}}
\newcommand{\ttS}{\overline{\overline{S}}}
\newcommand{\ttLa}{\overline{\overline{\Lambda}}}

\section{Introduction}
\label{ch1.intro}
Analysis of electromagnetic fields in cylindrically stratified media is of great importance in many applications, such as borehole geophysics~\cite{Wait:Geo, Telford:Applied, Ellis:Well}. 
This is a classical problem with separable geometry where the components of the tensor Green's function can be expressed 
in generic form as~\cite[Ch. 3]{Chew:Waves},\cite{Moon14:Stable}
\begin{flalign}
\suma e^{\iu n(\phi-\phi')}\intmp dk_z e^{\iu k_{z}(z-z')}
    \mathbf{\Phi}_n(\rho,\rho'),
\end{flalign}
where the integrand factor $\mathbf{\Phi}_n(\rho,\rho')$ contains various products of cylindrical Bessel and Hankel functions.
When applicable, such solutions are often preferred to brute-force numerical methods such as finite elements and finite difference 
 \cite{Wang01:3D, Weiss02:Electromagnetic, Weiss03:Electromagnetic, Lee07:Cylindrical, Lee12:Numerical,Pardo06:Two,Pardo06:Simulation,Novo07:Finite, Novo10:Three,Nam10:Simulation}
since the former can provide very accurate results with computational costs that are orders of magnitude smaller than 
the latter. This is especially important for inverse algorithms relying on repeated forward solutions and which seek to determine sought-after physical parameter values (say, layer resistivities) from the knowledge of the field values (measured) at certain subterranean locations.
  
However, numerical computations directly based on the canonical expressions of this problem can lead to underflow and overflow issues in finite precision arithmetic. This is caused by the poor scaling of
cylindrical Bessel and Hankel functions for extreme arguments and/or orders, which occur for low frequencies of operation and/or extreme values for layer resistivities.
In addition, convergence problems in the numerical evaluation of the spectral integral on the longitudinal wavenumber $k_z$ may occur depending on the 
separation distance between the source $(\rho',\phi',z')$ and observation point $(\rho,\phi,z)$ as well as on the operation frequency.
To circumvent these problems, a stable formulation based on a suitable analytical conditioning of the various factors in the computational chain and a proper choice of deformed integration paths in the complex $k_z$ plane 
was recently put forth in~\cite{Moon14:Stable}. This formulation was shown to be robust to variations on physical parameters that span several orders of magnitude. A related formulation to compute static fields (electric potentials) due to current electrodes in isotropic layers was described in~\cite{Moon15:Computation}.

In this work, we extend the formulation presented in~\cite{Moon14:Stable} to account for scenarios where the layers comprising the cylindrical stratified media may exhibit anisotropic properties. 
In borehole geophysics, anisotropy is quite common
\cite{Kunz58:Some,Teitler70:Refraction,Kong72:Electromagnetic,Moran79:Effects,Morgan87:Electromagnetic,Nekut94:Anisotropy,Bittar96:Effects,Howard00:Petrophysics,Yin01:Electromagnetic,Zhang04:Determination,Wang06:Weak,Hue07:Numerical,Wang08:Numerical,Zhong08:Computation,Yuan10:Simulation,Hagiwara11:Apparent,Hagiwara12:Determination,Liu12:Analysis,Luling13:Paradox,Sainath14:Robust,Sainath14:Tensor} and may result from geological factors affecting the various Earth layers such as salt water penetrating
porous fractured formations and thereby increasing the conductivity in
the direction parallel to the fracture and/or the presence of clay and sand laminates with directionally dependent
resistivities. Here, for generality, we assume each layer to be doubly uniaxial, i.e., both the complex permittivity tensor $\ue$ (which includes the conductivity tensor) and the permeability tensor $\um$ are independently uniaxial, which facilitates the analysis of equivalent problems using electromagnetic duality~\cite[Ch. 1]{Chew:Waves}.

\section{Fields in cylindrically-layered uniaxial media}
\label{sec.2.formul}
Most of the basic notation and terminology is adopted from \cite[Ch. 3]{Chew:Waves}. The section can be regarded as a generalization of the formulation presented for isotropic layers in \cite{Moon14:Stable} to uniaxial anisotropic layers.

\subsection{General solution in homogeneous, uniaxial media}
\label{sec.2.1}
Maxwell's curl equations in uniaxial, homogeneous, and source-free media (with time-harmonic dependence $e^{-\iu\omega t}$ assumed) read as
\begin{flalign}
\na\times\mathbf{E}&=\iu\omega\um\mathbf{H}, \label{ch1.1.E.Maxwell.1}\\
\na\times\mathbf{H}&=-\iu\omega\ue\mathbf{E}, \label{ch1.1.E.Maxwell.2}
\end{flalign}
where $\um$ and $\ue$  are the permeability tensor and complex permittivity tensor, respectively. In the unixial case, $\um$ is written as
\begin{flalign}
\renewcommand{\arraystretch}{1.2}
\um=
    \begin{bmatrix}
    \mu_h &     0 &     0 \\
        0 & \mu_h &     0 \\
        0 &     0 & \mu_v \\
    \end{bmatrix}, \label{ch1.1.E.uniaxial.mu}
\end{flalign}
where $\mu_h$ and $\mu_v$ are the horizontal and vertical permeabilities, resp. 
The complex permittivity tensor $\ue$ includes the electric conductivity and it is written as
\begin{flalign}
\renewcommand{\arraystretch}{1.2}
\ue
=   \begin{bmatrix}
    \epsilon_{h} &              0 &              0 \\
                 0 & \epsilon_{h} &              0 \\
                 0 &              0 & \epsilon_{v} \\
    \end{bmatrix}
=   \begin{bmatrix}
    \epsilon_{p,h} +\iu\sigma_h/\omega &                              0 &                             0 \\
                                 0 & \epsilon_{p,h} +\iu\sigma_h/\omega &                             0 \\
                                 0 &                              0 & \epsilon_{p,v} +\iu\sigma_v/\omega \\
    \end{bmatrix}, \label{ch1.1.E.uniaxial.epsilon}
\end{flalign}
where $\epsilon_{p,h}$ and $\epsilon_{p,v}$ are horizontal and vertical permittivities, and $\sigma_h$ and $\sigma_v$ are horizontal and vertical conductivities, resp. 
In such source-free media, the divergence equations can be written as
\begin{flalign}
\na\cdot\left( \ue \cdot \mathbf{E} \right) &=0, \label{ch1.1.E.Maxwell.3}\\
\na\cdot\left( \um \cdot \mathbf{H} \right) &=0. \label{ch1.1.E.Maxwell.4}
\end{flalign}
Note that in general $\na\cdot\mathbf{E}$ and $\na\cdot\mathbf{H}$ in uniaxial and source-free media are nonzero. Indeed, the left hand side of \eqref{ch1.1.E.Maxwell.3} in cylindrical coordinates is written as
\begin{flalign}
\na\cdot\ue\mathbf{E}
=\epsilon_h
	\left\{
		\frac{1}{\rho}\frac{\pa\left(\rho E_\rho\right)}{\pa\rho} 
		+ \frac{1}{\rho}\frac{\pa E_\phi}{\pa\phi}
		+ \frac{\pa E_z}{\pa z}
	-\left(1-\frac{\epsilon_v}{\epsilon_h}\right)\frac{\pa E_z}{\pa z}
	\right\}
=\epsilon_h
	\left\{
		\na\cdot\mathbf{E} - \left(1-\frac{\epsilon_v}{\epsilon_h}\right)\frac{\pa E_z}{\pa z}
	\right\}. \label{ch1.1.E.div.eE}
\end{flalign}
From \eqref{ch1.1.E.Maxwell.3} and \eqref{ch1.1.E.div.eE}, we can obtain
\begin{flalign}
\na\cdot\mathbf{E}
=\left(1-\frac{\epsilon_v}{\epsilon_h}\right)\frac{\pa E_z}{\pa z}. \label{ch1.1.E.div.E}
\end{flalign}
Similarly, we can obtain
\begin{flalign}
\na\cdot\mathbf{H}
=\left(1-\frac{\mu_v}{\mu_h}\right)\frac{\pa H_z}{\pa z}. \label{ch1.1.E.div.H}
\end{flalign}
To obtain the vector wave equation for $\mathbf{E}$, taking the curl of \eqref{ch1.1.E.Maxwell.1} and using \eqref{ch1.1.E.div.E} yields
\begin{flalign}
\na^2\mathbf{E} - \left(1-\frac{\epsilon_v}{\epsilon_h}\right)\na\frac{\pa E_z}{\pa z}
&=-\iu\omega
\left[\left(\na_s + \hat{z}\frac{\pa}{\pa z}\right) \times 
	      \left(\mu_h\mathbf{H}_s + \mu_v\mathbf{H}_z \right)
	\right], \label{ch1.1.E.curl.curl.E}
\end{flalign}
where $\na=\na_s + \hat{z}\frac{\pa}{\pa z}$ is used and subscript $s$ indicates the transverse components to the $z$-component. Similarly, the vector wave equation for $\mathbf{H}$ can be obtained by taking the curl of \eqref{ch1.1.E.Maxwell.2} and using \eqref{ch1.1.E.div.H} such that
\begin{flalign}
\na^2\mathbf{H} - \left(1-\frac{\mu_v}{\mu_h}\right)\na\frac{\pa H_z}{\pa z}
&=\iu\omega
\left[\left(\na_s + \hat{z}\frac{\pa}{\pa z}\right) \times 
	      \left(\epsilon_h\mathbf{E}_s + \epsilon_v\mathbf{E}_z \right)
	\right]. \label{ch1.1.E.curl.curl.H}
\end{flalign}
When the $z$-components are extracted from \eqref{ch1.1.E.curl.curl.E} and \eqref{ch1.1.E.curl.curl.H}, the equations for the $z$-components are written as
\begin{subequations}
\begin{flalign}
&\na^2 E_z - \left(1-\frac{\epsilon_v}{\epsilon_h}\right)\frac{\pa^2 E_z}{\pa z^2}
	+ \omega^2\mu_h\epsilon_v E_z = 0, \label{ch1.1.E.curl.curl.Ez}\\
&\na^2 H_z - \left(1-\frac{\mu_v}{\mu_h}\right)\frac{\pa^2 H_z}{\pa z^2}
	+ \omega^2\mu_v\epsilon_h H_z = 0. \label{ch1.1.E.curl.curl.Hz}
\end{flalign}
\end{subequations}
As usual, $E_z$ and $H_z$ can be solved for using the separation of variables technique. We define the propagation constant as 
$k = \omega\sqrt{\mu_h\epsilon_h}$ with
dispersion relation 
$\omega^2\mu_h\epsilon_h - k_z^2=k_\rho^2$ for the longitudinal (or vertical) $k_z$ and transverse (or radial) $k_\rho$ wavenumbers.
Two different anisotropic ratios can be defined in such media: the anisotropy ratio for the complex permittivity as
$\kappa_\epsilon=\sqrt{\frac{\epsilon_h}{\epsilon_v}}$,
and the anisotropy ratio for permeability as
$\kappa_\mu=\sqrt{\frac{\mu_h}{\mu_v}}$.
It is also convenient to define two scaled radial wavenumbers as
$\tk_\rho = \frac{k_\rho}{\kappa_\epsilon}$ and
$\dk_\rho = \frac{k_\rho}{\kappa_\mu}$.
With the above definitions, the general solution to the vector wave equation in such media becomes
\begin{subequations}
\begin{flalign}
E_z&=\left[A_n\Jn\left(\tk_\rho \rho\right)
	+B_n\Hn\left(\tk_\rho \rho\right)\right]e^{\iu n\phi}e^{\iu k_z z}, \label{ch1.1.E.gen.sol.Ez}\\
H_z&=\left[C_n\Jn\left(\dk_\rho \rho\right)
	+D_n\Hn\left(\dk_\rho \rho\right)\right]e^{\iu n\phi}e^{\iu k_z z}, \label{ch1.1.E.gen.sol.Hz}
\end{flalign}
\end{subequations}
with $A_n$, $B_n$, $C_n$, and $D_n$ determined by boundary conditions.
The dispersion relations for $\tk_\rho$ and $\dk_\rho$ are
\begin{subequations}
\begin{flalign}
\frac{\epsilon_v}{\epsilon_h}\left(\omega^2\mu_h\epsilon_h - k_z^2\right)
&= \tk_\rho^2, \label{ch1.1.E.dispersion.relation.eps}\\
\frac{\mu_v}{\mu_h}\left(\omega^2\mu_h\epsilon_h - k_z^2\right)
&= \dk_\rho^2. \label{ch1.1.E.dispersion.relation2.mu}
\end{flalign}
\end{subequations}

The transverse ($\rho$ and $\phi$) field components can be expressed in terms of the above longitudinal components directly by using Maxwell's equations \cite{Chew:Waves}, to yield
\begin{subequations}
\begin{flalign}
\mathbf{E}_s
	&=\frac{1}{\omega^2\mu_h\epsilon_h-k_z^2}
		\left[\na_s\frac{\pa E_z}{\pa z} - \iu\omega\mu_h\hat{z}\times\na_s H_z\right]
	=\frac{1}{k_\rho^2}
		\left[\iu k_z\na_s E_z - \iu\omega\mu_h\hat{z}\times\na_s H_z\right], \label{ch1.1.E.Etrans}\\
\mathbf{H}_s
	&=\frac{1}{\omega^2\mu_h\epsilon_h-k_z^2}
		\left[\na_s\frac{\pa H_z}{\pa z} + \iu\omega\epsilon_h\hat{z}\times\na_s E_z\right]
	=\frac{1}{k_\rho^2}
		\left[\iu k_z\na_s H_z + \iu\omega\epsilon_h\hat{z}\times\na_s E_z\right]. \label{ch1.1.E.Htrans}
\end{flalign}
\end{subequations}
or, in a convenient matrix form,
\renewcommand{\arraystretch}{1.4}
\begin{subequations}
\begin{flalign}
    \begin{bmatrix}
    E_\rho \\ H_\rho
    \end{bmatrix}
&=\frac{1}{k^2_\rho}
    \begin{bmatrix}
    \iu k_z\frac{\pa}{\pa\rho} & -\frac{n\omega\mu_h}{\rho} \\
    \frac{n\omega\epsilon_h}{\rho} & \iu k_z\frac{\pa}{\pa\rho}
    \end{bmatrix}
    \begin{bmatrix}
    E_z \\ H_z
    \end{bmatrix}
=\frac{1}{k^2_\rho}\bBn
	\begin{bmatrix}
    E_z \\ H_z
    \end{bmatrix}, \label{ch1.1.E.EHrho}\\
    \begin{bmatrix}
    E_\phi \\ H_\phi
    \end{bmatrix}
&=\frac{1}{k^2_\rho}
    \begin{bmatrix}
    -\frac{nk_z}{\rho} & -\iu\omega\mu_h\frac{\pa}{\pa\rho} \\
    \iu\omega\epsilon_h\frac{\pa}{\pa\rho} & -\frac{nk_z}{\rho} \\
    \end{bmatrix}
    \begin{bmatrix}
    E_z \\ H_z
    \end{bmatrix}
=\frac{1}{k^2_\rho}\bCn
	\begin{bmatrix}
    E_z \\ H_z
    \end{bmatrix}. \label{ch1.1.E.EHphi}
\end{flalign}
\end{subequations}
It should be noted that $\bBn$ and $\bCn$ depend on horizontal medium properties $\epsilon_h$ and $\mu_h$.

\subsection{Local reflection and transmission coefficients}
\label{sec.2.2}
When a number of cylindrical layers with different properties are present, the appropriate boundary conditions need to be enforced into the solutions. This is typically done via reflections and transmissions coefficients. In this section, local reflection and transmission coefficients are first derived for the two-layer case. Later, this will be extended to the case with arbitrary number of layers. 

\begin{figure}[t]
	\centering
	\subfloat[\label{F.2layers.out}]{%
      \includegraphics[width=2.2in]{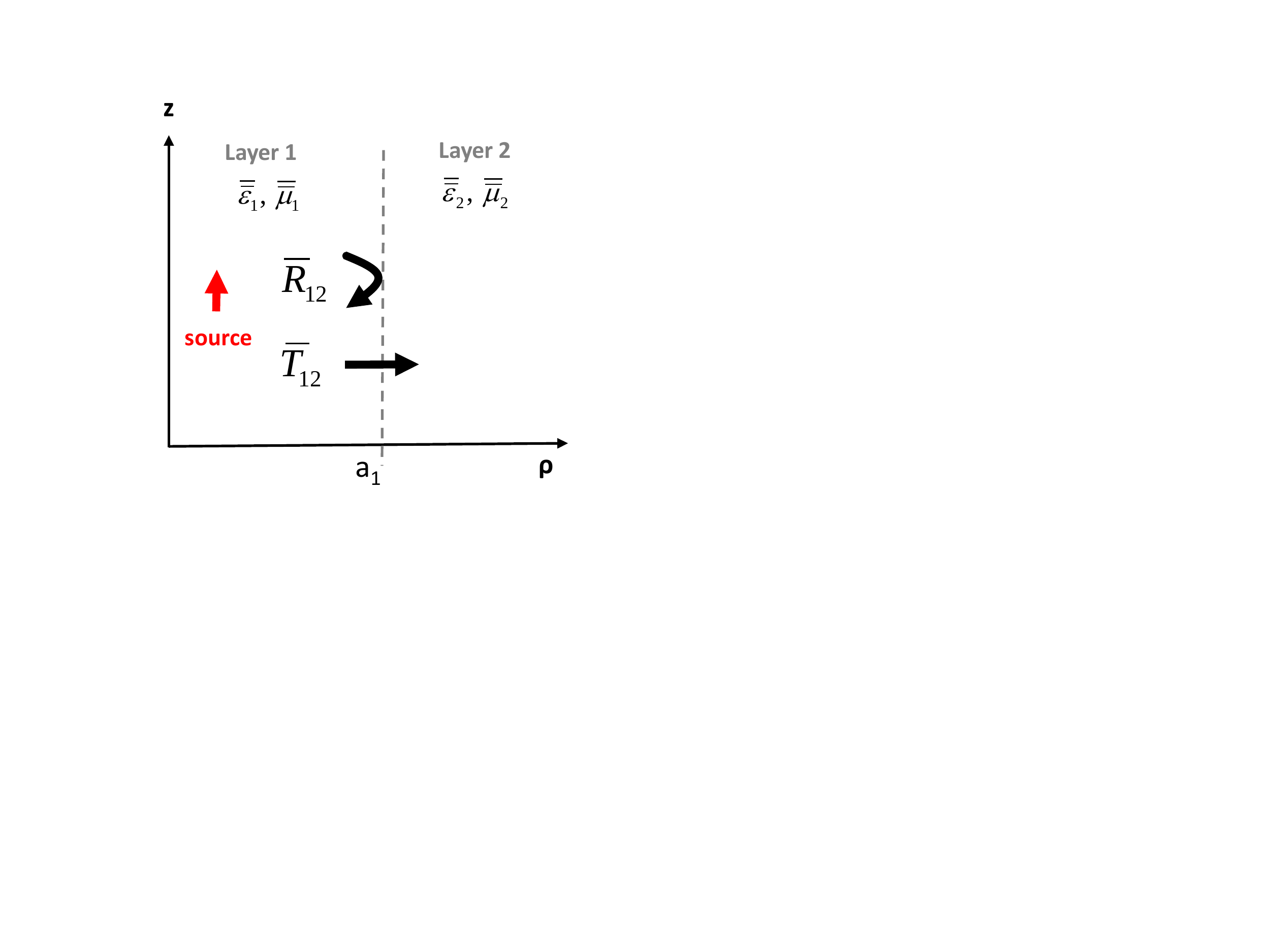}
    }
    \hspace{2.5 cm}    
    \subfloat[\label{F.2layers.stand}]{%
      \includegraphics[width=2.2in]{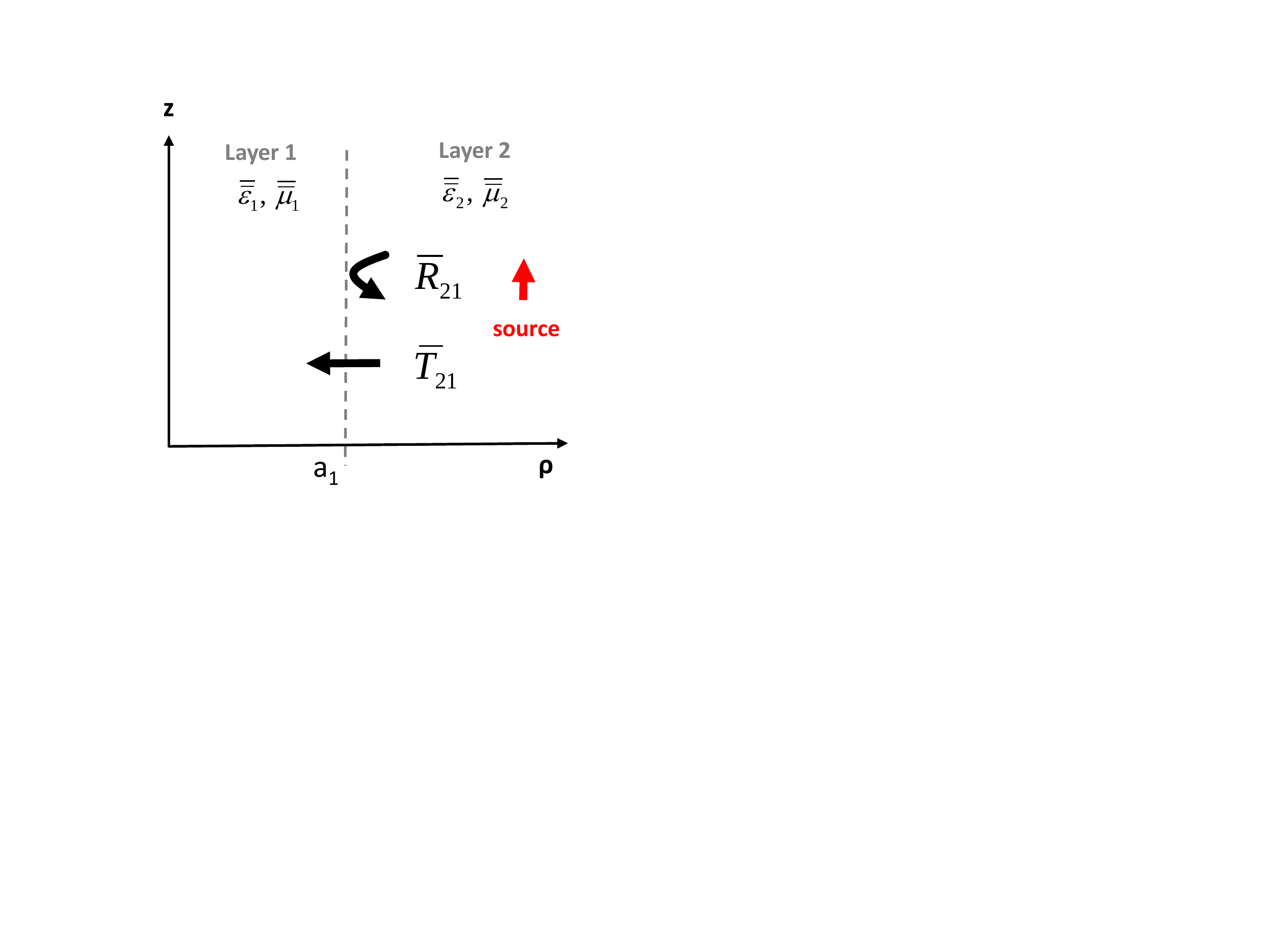}
	}    
    \caption{Two different cases of two uniaxial cylindrical layers with relevant reflection and transmission coefficients in the $\rho z$-plane: (a) Outgoing-wave case and (b) Standing-wave case.}
    \label{F.2layers}
\end{figure}

The local coefficients can be classified into two types: outgoing-wave type and standing-wave type depending on the relative location of the source versus the observation point, as illustrated in Fig. \ref{F.2layers}. As the general solution of $E_z$ and $H_z$ for uniaxial media are slightly different from those for isotropic media, so are the local reflection and transmission coefficients.  Nevertheless, the expressions for the {\it generalized} reflection and transmission coefficients (to account for more than two layers) {\it in terms of local ones} remain the same as those in~\cite{Moon14:Stable}. 

\subsubsection{Outgoing-wave case}
\label{sec.2.2.1}
Based on \eqref{ch1.1.E.gen.sol.Ez} and \eqref{ch1.1.E.gen.sol.Hz}, outgoing waves in a uniaxial medium  can be expressed as
\begin{flalign}
\renewcommand{\arraystretch}{1.2}
    \begin{bmatrix}
    E_z \\ H_z
    \end{bmatrix}
=   \begin{bmatrix}
    \Hn(\tk_\rho \rho) & 0 \\
     0 & \Hn(\dk_\rho \rho) \\
    \end{bmatrix}
    \begin{bmatrix}
    e_z \\ h_z
    \end{bmatrix}    
= \bHzn(k_\rho \rho)\cdot\mathbf{a}, \label{ch1.2.E.outgoing.EHz.gen}
\end{flalign}
where the column vector $\mathbf{a}$ includes $e^{\iu n\phi + \iu k_z z}$ dependence. 
Since the source is embedded in layer 1 for the outgoing-wave case depicted in Fig. \ref{F.2layers.out}, the $z$-components of the total fields in layer 1 and layer 2 is expressed as
\begin{subequations}
\begin{flalign}
    \begin{bmatrix}
    E_{z1} \\ H_{z1}
    \end{bmatrix}
&= \bHzn(k_{1\rho}\rho)\cdot\mathbf{a}_1
	+ \bJzn(k_{1\rho}\rho)\cdot\bR_{12}\cdot\mathbf{a}_1, \label{ch1.2.E.out.EHz.1}\\
    \begin{bmatrix}
    E_{z2} \\ H_{z2}
    \end{bmatrix}
&= \bHzn(k_{2\rho}\rho)\cdot\bT_{12}\cdot\mathbf{a}_1. \label{ch1.2.E.out.EHz.2}
\end{flalign}
\end{subequations}
Likewise, using \eqref{ch1.1.E.EHphi}, we have
\begin{subequations}
\begin{flalign}
    \begin{bmatrix}
    H_{\phi 1} \\ E_{\phi 1}
    \end{bmatrix}
&= \bHpn(k_{1\rho}\rho)\cdot\mathbf{a}_1
	+ \bJpn(k_{1\rho}\rho)\cdot\bR_{12}\cdot\mathbf{a}_1, \label{ch1.2.E.out.EHphi.1} \\
    \begin{bmatrix}
    H_{\phi 2} \\ E_{\phi 2}
    \end{bmatrix}
&= \bHpn(k_{2\rho}\rho)\cdot\bT_{12}\cdot\mathbf{a}_1. \label{ch1.2.E.out.EHphi.2}
\end{flalign}
\end{subequations}
From \eqref{ch1.2.E.out.EHz.1} through \eqref{ch1.2.E.out.EHphi.2}, two types of matrices only depending on $\rho$ are defined as
\begin{subequations}
\begin{flalign}
\overline{\mathbf{B}}_{zn}(k_{i\rho}\rho)
&=
   \begin{bmatrix}
    B_n(\tk_{i\rho}\rho) & 0 \\
     0 & B_n(\dk_{i\rho}\rho) \\
    \end{bmatrix}
,\label{ch1.2.E.Bzn}\\
\overline{\mathbf{B}}_{\phi n}(k_{i\rho}\rho)
&=\frac{1}{k_{i\rho}^2\rho}
    \begin{bmatrix}
    \iu\omega\epsilon_{hi}\tk_{i\rho}\rho B'_n(\tk_{i\rho}\rho) & -n k_z B_n(\dk_{i\rho}\rho) \\
    -n k_z B_n(\tk_{i\rho}\rho) & -\iu\omega\mu_{hi} \dk_{i\rho}\rho B'_n(\dk_{i\rho}\rho) \\
    \end{bmatrix}
,\label{ch1.2.E.Bpn}
\end{flalign}
\end{subequations}
where $B_n$ is either $\Hn$ or $\Jn$, $k_{i\rho}=\omega^2\mu_{hi}\epsilon_{hi} - k_z^2$, and $\epsilon_{hi}$ and $\mu_{hi}$ are the horizontal complex permittivity and permeability in layer $i$, respectively. Applying the pertinent boundary conditions, viz., continuity of $z$- and $\phi$-components at $\rho=a_1$, to \eqref{ch1.2.E.out.EHz.1}--\eqref{ch1.2.E.out.EHphi.2} yields
\begin{subequations}
\begin{flalign}
\left[\bHzn(k_{1\rho}a_1) + \bJzn(k_{1\rho}a_1)\cdot\bR_{12}\right]\cdot\mathbf{a}_1
&= \bHzn(k_{2\rho}a_1)\cdot\bT_{12}\cdot\mathbf{a}_1, \label{ch1.2.E.BC1}\\
\left[\bHpn(k_{1\rho}a_1)	+ \bJpn(k_{1\rho}a_1)\cdot\bR_{12}\right]\cdot\mathbf{a}_1
&= \bHpn(k_{2\rho}a_1)\cdot\bT_{12}\cdot\mathbf{a}_1. \label{ch1.2.E.BC2}
\end{flalign}
\end{subequations}
For notational simplicity, a shorthand notation is defined such that
\begin{flalign}
\overline{\mathbf{B}}_{\alpha n}(k_{i\rho}a_j)=\overline{\mathbf{B}}_{\alpha ij}. \label{ch1.2.E.shorthands}
\end{flalign}
In the right hand side of \eqref{ch1.2.E.shorthands}, the first, second, and third subscripts indicate the relevant components of the fields, radial wavenumbers, and radial distances, respectively. For notational simplicity, the Hankel function superscript and subscript (kind and modal number) are suppressed in the following. Consequently, \eqref{ch1.2.E.BC1} and \eqref{ch1.2.E.BC2} are re-expressed as
\begin{subequations}
\begin{flalign}
\left[\bH_{z11}+\bJ_{z11}\cdot\bR_{12}\right] &= \bH_{z21}\cdot\bT_{12}, \label{ch1.2.E.BC1.R12}\\
\left[\bH_{\phi 11}+\bJ_{\phi 11}\cdot\bR_{12}\right] &= \bH_{\phi 21}\cdot\bT_{12}. \label{ch1.2.E.BC2.R12}
\end{flalign}
\end{subequations}
From \eqref{ch1.2.E.BC1.R12} and \eqref{ch1.2.E.BC2.R12}, we obtain
\begin{subequations}
\begin{flalign}
\bR_{12}
&=\left[\bJ_{z11} - \bH_{z21}\cdot\bH_{\phi 21}^{-1}\cdot\bJ_{\phi 11}\right]^{-1}\cdot
  \left[\bH_{z21}\cdot\bH_{\phi 21}^{-1}\cdot\bH_{\phi 11}-\bH_{z11}\right], \label{ch1.2.E.R12}\\
\bT_{12}
&=\left[\bH_{z21} - \bJ_{z11}\cdot\bJ_{\phi 11}^{-1}\cdot\bH_{\phi 21}\right]^{-1}\cdot
  \left[\bH_{z11} - \bJ_{z11}\cdot\bJ_{\phi 11}^{-1}\cdot\bH_{\phi 11}\right].	\label{ch1.2.E.T12}
\end{flalign}
\end{subequations}

\subsubsection{Standing-wave case}
\label{sec.2.2.2}
For the standing-wave case depicted in Fig. \ref{F.2layers.stand}, we now have
\begin{subequations}
\begin{flalign}
    \begin{bmatrix}
    E_{z1} \\ H_{z1}
    \end{bmatrix}
&= \bJzn(k_{1\rho}\rho)\cdot\bT_{21}\cdot\mathbf{a}_2,
\label{ch1.2.E.stand.EHz.1}\\
    \begin{bmatrix}
    E_{z2} \\ H_{z2}
    \end{bmatrix}
&= \bHzn(k_{2\rho}\rho)\cdot\bR_{21}\cdot\mathbf{a}_2
	+ \bJzn(k_{2\rho}\rho)\cdot\mathbf{a}_2,
\label{ch1.2.E.stand.EHz.2}
\end{flalign}
\end{subequations}
and, using \eqref{ch1.1.E.EHphi},
\begin{subequations}
\begin{flalign}
    \begin{bmatrix}
    H_{\phi 1} \\ E_{\phi 1}
    \end{bmatrix}
&= \bJpn(k_{1\rho}\rho)\cdot\bT_{21}\cdot\mathbf{a}_2,
\label{ch1.2.E.stand.EHphi.1} \\
    \begin{bmatrix}
    H_{\phi 2} \\ E_{\phi 2}
    \end{bmatrix}
&= \bHpn(k_{2\rho}\rho)\cdot\bR_{21}\cdot\mathbf{a}_2
	+ \bJpn(k_{2\rho}\rho)\cdot\mathbf{a}_2.
\label{ch1.2.E.stand.EHphi.2}
\end{flalign}
\end{subequations}
Applying the boundary conditions at $\rho=a_1$ to \eqref{ch1.2.E.stand.EHz.1}--\eqref{ch1.2.E.stand.EHphi.2} yields
\begin{subequations}
\begin{flalign}
\bJzn(k_{1\rho}a_1)\cdot\bT_{21}\cdot\mathbf{a}_2
&=\left[\bHzn(k_{2\rho}a_1)\cdot\bR_{21} + \bJzn(k_{2\rho}a_1)\right]\cdot\mathbf{a}_2, \label{ch1.2.E.BC3}\\
\bJpn(k_{1\rho}a_1)\cdot\bT_{21}\cdot\mathbf{a}_2
&=\left[\bHpn(k_{2\rho}a_1)\cdot\bR_{21} + \bJpn(k_{2\rho}a_1)\right]\cdot\mathbf{a}_2, \label{ch1.2.E.BC4}
\end{flalign}
\end{subequations}
which can be re-expressed as
\begin{subequations}
\begin{flalign}
\bJ_{z11}\cdot\bT_{21} &= \left[\bH_{z21}\cdot\bR_{21}+\bJ_{z21}\right] , \label{ch1.2.E.BC1.R21}\\
\bJ_{\phi 11}\cdot\bT_{21} &= \left[\bH_{\phi 21}\cdot\bR_{21}+\bJ_{\phi 21}\right]. \label{ch1.2.E.BC2.R21}
\end{flalign}
\end{subequations}
Therefore, $\bR_{21}$ and $\bT_{21}$ are written as
\begin{subequations}
\begin{flalign}
\bR_{21}
&=\left[\bH_{z21} - \bJ_{z11}\cdot\bJ_{\phi 11}^{-1}\cdot\bH_{\phi 21}\right]^{-1}\cdot
  \left[\bJ_{z11}\cdot\bJ_{\phi 11}^{-1}\cdot\bJ_{\phi 21} - \bJ_{z21}\right], \label{ch1.2.E.R21}\\
\bT_{21}
&=\left[\bJ_{z11} - \bH_{z21}\cdot\bH_{\phi 21}^{-1}\cdot\bJ_{\phi 11}\right]^{-1}\cdot
  \left[\bJ_{z21} - \bH_{z21}\cdot\bH_{\phi 21}^{-1}\cdot\bJ_{\phi 21}\right]. \label{ch1.2.E.T21}
\end{flalign}
\end{subequations}
Furthermore, the local reflection and transmission coefficients \eqref{ch1.2.E.R12}, \eqref{ch1.2.E.T12}, \eqref{ch1.2.E.R21}, and \eqref{ch1.2.E.T21} can be succinctly rewritten as
\begin{subequations}
\begin{flalign}
\bR_{12}
&=\bD_A^{-1}\cdot
  \left[\bH_{z21}\cdot\bH_{\phi 21}^{-1}\cdot\bH_{\phi 11}-\bH_{z11}\right], \label{ch1.2.E.R12.short}\\
\bR_{21}
&= \bD_B^{-1}\cdot
  \left[\bJ_{z11}\cdot\bJ_{\phi 11}^{-1}\cdot\bJ_{\phi 21} - \bJ_{z21}\right], \label{ch1.2.E.R21.short}\\
\bT_{12}
&= \bD_B^{-1}\cdot
  \left[\bH_{z11} - \bJ_{z11}\cdot\bJ_{\phi 11}^{-1}\cdot\bH_{\phi 11}\right], \label{ch1.2.E.T12.short}\\
\bT_{21}
&= \bD_A^{-1}\cdot
  \left[\bJ_{z21} - \bH_{z21}\cdot\bH_{\phi 21}^{-1}\cdot\bJ_{\phi 21}\right], \label{ch1.2.E.T21.short}
\end{flalign}
\end{subequations}
where
\begin{subequations}
\begin{flalign}
\bD_A&=\left[\bJ_{z11} - \bH_{z21}\cdot\bH_{\phi 21}^{-1}\cdot\bJ_{\phi 11}\right], \label{ch1.2.E.DA}\\
\bD_B&=\left[\bH_{z21} - \bJ_{z11}\cdot\bJ_{\phi 11}^{-1}\cdot\bH_{\phi 21}\right]. \label{ch1.2.E.DB}
\end{flalign}
\end{subequations}

\subsection{Spectral representation of the Green's function}
\label{sec.2.3}
In this section, we derive convenient analytical expressions for the Green's function expressed in cylindrical coordinates that will facilitate further analysis in such uniaxial media. Specifically, we obtain expressions for the $z$-components of the fields produced by a point source (arbitrarily-oriented Hertzian electric dipole). In this case, the source writes as
\begin{flalign}
\mathbf{J}(\rr)=Il\hat{\alpha}'\delta(\rr-\rp), \label{ch1.3.E.current.source}
\end{flalign}
where $Il$ is the dipole moment. By taking the curl of Faraday's law, we obtain
\begin{flalign}
\na\times\na\times\mathbf{E}
=\iu\omega\na\times\um\mathbf{H}. \label{ch1.3.E.curl.curl.w.source} 
\end{flalign}
With the presence of the source, \eqref{ch1.3.E.curl.curl.w.source} becomes
\begin{flalign}
\na^2\mathbf{E}-\left(1-\frac{\epsilon_v}{\epsilon_h}\right)\na\frac{\pa E_z}{\pa z}
=-\iu\omega\na\times\um\mathbf{H}
+\na\left(\frac{\na\cdot\mathbf{J}}{\iu\omega\epsilon_h}\right), \label{ch1.3.E.2nd.order.E.b}
\end{flalign}
where the divergence of $\mathbf{E}$ is deduced from \eqref{ch1.1.E.div.eE} and the continuity equation $\na\cdot\mathbf{J} - \iu\omega\rho_v=0$ is applied.
Extracting the $z$-component of the above, we can easily show that
\begin{flalign}
\na^2 E_z + \omega^2\mu_h\epsilon_v E_z - \left(1-\frac{\epsilon_v}{\epsilon_h}\right)\frac{\pa^2 E_z}{\pa z^2}
&=-\iu\omega\mu_h \hat{z}\cdot\mathbf{J}
	+ \frac{\pa}{\pa z}\left(\frac{\na\cdot\mathbf{J}}{\iu\omega\epsilon_h}\right). \label{ch1.3.E.2nd.order.Ez.a}
\end{flalign}
Using \eqref{ch1.3.E.current.source} and $e^{\iu k_z z}$ dependence shown in \eqref{ch1.1.E.gen.sol.Ez}, we obtain
\begin{flalign}
\na^2 E_z+\tk^2 E_z = -\frac{\iu Il}{\omega\epsilon_h}
\left[k^2(\hat{z}\cdot\hat{\alpha}')+\frac{\pa}{\pa z}\na\cdot\hat{\alpha}'\right]\delta(\rr-\rp),
	\label{ch1.3.E.2nd.order.Ez.c}
\end{flalign}
where
\begin{flalign}
\tk^2 &= \omega^2\mu_h\epsilon_v + \left(1-\frac{\epsilon_v}{\epsilon_h}\right)k_z^2
       = \frac{\epsilon_v}{\epsilon_h}\left(\omega^2\mu_h\epsilon_h-k_z^2\right)+k_z^2
	   = \tk_\rho^2 + k_z^2, \label{ch1.3.E.dis.rel.tk}\\
k^2   &= \omega^2\mu_h\epsilon_h. \label{ch1.3.E.dis.rel.k}
\end{flalign}
Therefore, $E_z$ is obtained via the scalar Green's function,
\begin{flalign}
E_z=\frac{\iu Il}{\omega\epsilon_h}
	\left[k^2(\hat{z}\cdot\hat{\alpha}')+\frac{\pa}{\pa z'}\na'\cdot\hat{\alpha}'
	\right]\widetilde{g}(\rr-\rp), \label{ch1.3.E.Ez}
\end{flalign}
where
\begin{flalign}
\widetilde{g}(\rr-\rp)=
    \frac{e^{\iu \tk|\rr-\rp|}}{4\pi|\rr-\rp|}. \label{ch1.3.E.scalar.stretched.Green.1}
\end{flalign}

The derivation for the $z$-component of the magnetic field can be done quite similarly by taking the curl of Ampere's law. The equivalent to
\eqref{ch1.3.E.2nd.order.Ez.a} becomes
\begin{flalign}
\na^2 H_z + \omega^2\mu_v\epsilon_h H_z- \left(1-\frac{\mu_v}{\mu_h}\right)\frac{\pa^2 H_z}{\pa z^2}
&= - \hat{z}\cdot\na\times\mathbf{J}. \label{ch1.3.E.2nd.order.Hz.a}
\end{flalign}
Since $H_z$ also has $e^{\iu k_z z}$ dependence as shown in \eqref{ch1.1.E.gen.sol.Hz}, \eqref{ch1.3.E.2nd.order.Hz.a} reduces to
\begin{flalign}
\na^2 H_z + \dk^2 H_z
&= - Il\hat{z}\cdot\na\times\hat{\alpha}'\delta(\rr-\rp), \label{ch1.3.E.2nd.order.Hz.b}
\end{flalign}
where
\begin{flalign}
\dk^2 = \omega^2\mu_v\epsilon_h+\left(1-\frac{\mu_v}{\mu_h}\right)k_z^2
	  = \frac{\mu_v}{\mu_h}\left(\omega^2\mu_h\epsilon_h - k_z^2\right) + k_z^2
	  = \dk_\rho^2 + k_z^2. \label{ch1.3.E.dis.rel.dk}
\end{flalign}
Now, $H_z$ is obtained via
\begin{flalign}
H_z = -Il\hat{z}\cdot\na'\times\hat{\alpha}' \ddot{g}(\rr-\rp), \label{ch1.3.E.Hz}
\end{flalign}
where
\begin{flalign}
\ddot{g}(\rr-\rp)=
    \frac{e^{\iu \dk|\rr-\rp|}}{4\pi|\rr-\rp|}. \label{ch1.3.E.scalar.stretched.Green.2}
\end{flalign}

The $z$-components of electromagnetic fields expressed in \eqref{ch1.3.E.Ez} and \eqref{ch1.3.E.Hz} can be expanded as the linear combination of all spectral components. In other words, using the spectral representation of the scalar Green's function,
\begin{flalign}
\frac{e^{\iu k|\rr-\rp|}}{|\rr-\rp|}
    =\suma\frac{\iu e^{\iu n(\phi-\phi')}}{2}\intmp dk_z
    e^{\iu k_z(z-z')}\Jn(k_\rho\rho_<)\Hn(k_\rho\rho_>), \label{ch1.3.E.expansion.scalar.Green}
\end{flalign}
$E_z$ and $H_z$ in homogeneous, uniaxial media can be written as
\begin{flalign}
    \begin{bmatrix}
    E_z \\ H_z
    \end{bmatrix}
=\frac{\iu Il}{4\pi\omega\epsilon_h}\suma e^{\iu n(\phi-\phi')}\intmp dk_z e^{\iu k_z(z-z')}
	\begin{bmatrix}
    \Jn(\tk_\rho \rho_<)\Hn(\tk_\rho \rho_>) & 0 \\ 0 & \Jn(\dk_\rho \rho_<)\Hn(\dk_\rho \rho_>)
    \end{bmatrix}
	\cdot\overleftarrow{\mathbf{D}}', \label{ch1.3.E.EzHz.single}
\end{flalign}
where $\overleftarrow{\mathbf{D}}'$ is an operator acting on the primed variables on the left and written as
\begin{flalign}
\overleftarrow{\mathbf{D}}'
&=\frac{\iu}{2}
    \begin{bmatrix}
    (\hat{z}k^2+\frac{\pa}{\pa z'}\na')\cdot\hat{\alpha}' \\
    \iu\omega\epsilon_h\hat{\alpha}'\cdot\hat{z}\times\na'
    \end{bmatrix}. \label{ch1.3.E.D.single}
\end{flalign}

In order to account for multilayers, let us consider cylindrically stratified media with the source in layer $j$. Then, the reflection terms from the boundaries at $\rho=a_j$ and $\rho=a_{j-1}$ are added to \eqref{ch1.3.E.EzHz.single} such that
\begin{flalign}
    \begin{bmatrix}
    E_{zj} \\ H_{zj}
    \end{bmatrix}
=\frac{\iu Il}{4\pi\omega\epsilon_{hj}}\suma e^{\iu n(\phi-\phi')}\intmp dk_z e^{\iu k_z(z-z')}
	\left\{\bJ_{zj\rho_<}\cdot\bH_{zj\rho_>}
+ \bH_{zj\rho}\cdot\overline{\mathbf{a}}_{jn}(\rho')
		      + \bJ_{zj\rho}\cdot\overline{\mathbf{b}}_{jn}(\rho')
	\right\}\cdot\Dj, \label{ch1.3.E.EzHz.multi}
\end{flalign}
where
\begin{subequations}
\begin{flalign}
\bJ_{zj\rho_<}\cdot\bH_{zj\rho_>}
&=	\begin{bmatrix}
    \Jn(\tk_{j\rho} \rho_<)\Hn(\tk_{j\rho} \rho_>) & 0 \\ 0 & \Jn(\dk_{j\rho} \rho_<)\Hn(\dk_{j\rho} \rho_>)
    \end{bmatrix}, \label{ch1.3.E.JH.zjrho}\\
\bH_{zj\rho}
&=	\begin{bmatrix}
    \Hn(\tk_{j\rho} \rho) & 0 \\ 0 & \Hn(\dk_{j\rho} \rho)
    \end{bmatrix}, \label{ch1.3.E.H.zjrho}\\
\bJ_{zj\rho}
&=  \begin{bmatrix}
    \Jn(\tk_{j\rho} \rho) & 0 \\ 0 & \Jn(\dk_{j\rho} \rho)
    \end{bmatrix}, \label{ch1.3.E.J.zjrho}\\
\Dj&=\frac{\iu}{2}
    \begin{bmatrix}
    (\hat{z}k_j^2-\iu k_z\na')\cdot\hat{\alpha}' \\
    \iu\omega\epsilon_{hj}\hat{\alpha}'\cdot\hat{z}\times\na'
    \end{bmatrix}. \label{ch1.3.E.Dj.multi}
\end{flalign}
\end{subequations}
Recall that notations of \eqref{ch1.3.E.JH.zjrho}, \eqref{ch1.3.E.H.zjrho}, and \eqref{ch1.3.E.J.zjrho} are based on \eqref{ch1.2.E.Bzn} and \eqref{ch1.2.E.shorthands}. In \eqref{ch1.3.E.Dj.multi}, $k_j=\omega\mu_{hj}\epsilon_{hj}$, and $\mu_{hj}$ and $\epsilon_{hj}$ represent horizontal permeability and complex permittivity in layer $j$, respectively. When $\hat{\alpha}'$ is represented in cylindrical coordinates as $\hat{\alpha}'=\hat{\rho}'\alpha_{\rho'}+\hat{\phi}'\alpha_{\phi'}+\hat{z}'\alpha_{z'}$,
\begin{flalign}
\Dj=
    \frac{\iu}{2}
    \left(\Dja+\Djb+\frac{\pa}{\pa\rho'}\Djc\right)
=\frac{\iu}{2}
\left(
    \begin{bmatrix}
    (k^2_{j\rho})\alpha_{z'} \\ 0
    \end{bmatrix}
    +
    \begin{bmatrix}
    -\frac{nk_z}{\rho'}\alpha_{\phi'} \\
    -\frac{n\omega\epsilon_{hj}}{\rho'}\alpha_{\rho'}
    \end{bmatrix}
    +\frac{\pa}{\pa\rho'}
    \begin{bmatrix}
    -\iu k_z\alpha_{\rho'} \\
    \iu\omega\epsilon_{hj}\alpha_{\phi'}
    \end{bmatrix}
\right). \label{ch1.3.E.Dj.cylin}
\end{flalign}
Using the two constraint conditions at $\rho=a_j$ and $\rho=a_{j-1}$ detailed in \cite{Chew:Waves}, two unknowns $\overline{\mathbf{a}}_{jn}(\rho')$ and $\overline{\mathbf{b}}_{jn}(\rho')$ appeared in \eqref{ch1.3.E.EzHz.multi} can be determined such that
\begin{subequations}
\begin{flalign}
\overline{\mathbf{a}}_{jn}
&=\left[\bI-\tbR_{j,j-1}\cdot\tbR_{j,j+1}\right]\cdot\tbR_{j,j-1}\cdot
	\left[\bH_{zj\rho'}+\tbR_{j,j+1}\cdot\bJ_{zj\rho'}\right], \label{ch1.3.E.an}\\
\overline{\mathbf{b}}_{jn}
&=\left[\bI-\tbR_{j,j+1}\cdot\tbR_{j,j-1}\right]\cdot\tbR_{j,j+1}\cdot
	\left[\bJ_{zj\rho'}+\tbR_{j,j-1}\cdot\bH_{zj\rho'}\right]. \label{ch1.3.E.bn}
\end{flalign}
\end{subequations}
Note that the second brackets in the right hand sides of \eqref{ch1.3.E.an} and \eqref{ch1.3.E.bn} are slightly different from those for isotropic media. When $\rho>\rho'$, $\rho_<=\rho'$ and $\rho_>=\rho$. The curly bracket in \eqref{ch1.3.E.EzHz.multi} is expressed as
\begin{flalign}
\bJ_{zj\rho_<}\cdot\bH_{zj\rho_>} &+ \bH_{zj\rho}\cdot\overline{\mathbf{a}}_{jn}(\rho')
		      + \bJ_{zj\rho}\cdot\overline{\mathbf{b}}_{jn}(\rho') \notag\\
&=\left[\bH_{zj\rho}+\bJ_{zj\rho}\cdot\tbR_{j,j+1}\right]\cdot\tbM_{j+}\cdot
	\left[\bJ_{zj\rho'}+\tbR_{j,j-1}\cdot\bH_{zj\rho'}\right]. \label{ch1.3.E.bracket.case1}
\end{flalign}
On the other hand, when $\rho<\rho'$, $\rho_<=\rho$ and $\rho_>=\rho'$. The curly bracket \eqref{ch1.3.E.EzHz.multi} is now expressed as
\begin{flalign}
\bJ_{zj\rho_<}\cdot\bH_{zj\rho_>} &+ \bH_{zj\rho}\cdot\overline{\mathbf{a}}_{jn}(\rho')
		      + \bJ_{zj\rho}\cdot\overline{\mathbf{b}}_{jn}(\rho') \notag\\
&=\left[\bJ_{zj\rho}+\bH_{zj\rho}\cdot\tbR_{j,j-1}\right]\cdot\tbM_{j-}\cdot
	\left[\bH_{zj\rho'}+\tbR_{j,j+1}\cdot\bJ_{zj\rho'}\right]. \label{ch1.3.E.bracket.case2}
\end{flalign}
Again, \eqref{ch1.3.E.bracket.case1} and \eqref{ch1.3.E.bracket.case2} are slightly different from those for isotropic media. When the field layer is not the same as the source layer, the approach is the same as that in \cite{Chew:Waves}. In summary,
\begin{flalign}
    \begin{bmatrix}
    E_z \\ H_z
    \end{bmatrix}
=\frac{\iu Il}{4\pi\omega\epsilon_{hj}}\suma e^{\iu n(\phi-\phi')}\intmp dk_z e^{\iu k_{z}(z-z')}
    \Fn\cdot\Dj, \label{ch1.3.E.EzHz.general}
\end{flalign}
where: \\
\begin{subequations}
\text{for Case 1: $\rho$ and $\rho'$ are in the same region and $\rho\geq\rho'$.}
\begin{flalign}
&\qquad\Fn=
	\left[\bH_{zj\rho}+\bJ_{zj\rho}\cdot\tbR_{j,j+1}\right]\cdot
	\tbM_{j+}\cdot
	\left[\bJ_{zj\rho'}+\tbR_{j,j-1}\cdot\bH_{zj\rho'}\right],&
    \label{ch1.3.E.EzHz.general.case1}
\end{flalign}
\text{for Case 2: $\rho$ and $\rho'$ are in the same region and $\rho<\rho'$.}
\begin{flalign}
&\qquad\Fn=
	\left[\bJ_{zj\rho}+\bH_{zj\rho}\cdot\tbR_{j,j-1}\right]\cdot
	\tbM_{j-}\cdot
	\left[\bH_{zj\rho'}+\tbR_{j,j+1}\cdot\bJ_{zj\rho'}\right],&	
    \label{ch1.3.E.EzHz.general.case2}
\end{flalign}
\text{for Case 3: $\rho$ and $\rho'$ are in different regions and $\rho>\rho'$.}
\begin{flalign}
&\qquad\Fn=
	\left[\bH_{zi\rho}+\bJ_{zi\rho}\cdot\tbR_{i,i+1}\right]\cdot
	\bN_{i+}\cdot\tbT_{ji}\cdot\tbM_{j+}\cdot
	\left[\bJ_{zj\rho'}+\tbR_{j,j-1}\cdot\bH_{zj\rho'}\right],&
    \label{ch1.3.E.EzHz.general.case3}
\end{flalign}
\text{for Case 4: $\rho$ and $\rho'$ are in different regions and $\rho<\rho'$.}
\begin{flalign}
&\qquad\Fn=
	\left[\bJ_{zi\rho}+\bH_{zi\rho}\cdot\tbR_{i,i-1}\right]\cdot
	\bN_{i-}\cdot\tbT_{ji}\cdot\tbM_{j-}\cdot	
	\left[\bH_{zj\rho'}+\tbR_{j,j+1}\cdot\bJ_{zj\rho'}\right].&
    \label{ch1.3.E.EzHz.general.case4}
\end{flalign}
\end{subequations}

\section{Range-conditioned formulation}
\label{sec.3.range.formul}
As noted before, the poor scaling behavior of Bessel and Hankel functions for extreme arguments and/or orders causes instabilities in the numerical computation of the field expressions above under some parameter ranges.  This section discusses how to stabilize the computation and provides relevant mathematical derivations.

\subsection{Range-conditioned cylindrical functions and matrices}
\label{sec.3.1}
Range-conditioned cylindrical functions derived in~\cite{Moon14:Stable} are modified for uniaxial media because two different scaled radial wavenumbers $\tk_{i\rho}$ and $\dk_{i\rho}$ appear as part of the function arguments. Here,
the subscript $i$ is the layer index, and, to recall, $\kappa_{i\epsilon}=\sqrt{\epsilon_{hi}/\epsilon_{vi}}$ and $\kappa_{i\mu}=\sqrt{\mu_{hi}/\mu_{vi}}$ are the anisotropy ratios of complex permittivity and permeability in layer $i$, respectively. Here, we will only focus on those aspects that differ from~\cite{Moon14:Stable}. The reader is referred to~\cite{Moon14:Stable} for the fundamentals of this stabilization approach.
Table \ref{ch2.1.T.rccf} shows the definitions of the range-conditioned cylindrical functions, indicated by a hat, for uniaxial complex permittivity media, where:
\begin{table}[t]
\begin{center}
\renewcommand{\arraystretch}{1.6}
\setlength{\tabcolsep}{8pt}
\caption{Definition of range-conditioned cylindrical functions for uniaxial complex permittivity media.}
\begin{tabular}{ccc}
    \hline
    Small arguments & Moderate arguments & Large arguments \\
    \hline
    $\Jn(\tk_{i\rho}a_j)=\widetilde{G}_i a_j^n\hJn(\tk_{i\rho}a_j)$ &
    $\Jn(\tk_{i\rho}a_j)=\widetilde{P}_{ij}\hJn(\tk_{i\rho}a_j)$ &
    $\Jn(\tk_{i\rho}a_j)=e^{|\tk''_{i\rho}|a_j}\hJn(\tk_{i\rho}a_j)$ \\       
    $\Jnd(\tk_{i\rho}a_j)=\widetilde{G}_i a_j^n\hJnd(\tk_{i\rho}a_j)$ &
    $\Jnd(\tk_{i\rho}a_j)=\widetilde{P}_{ij}\hJnd(\tk_{i\rho}a_j)$ &
    $\Jnd(\tk_{i\rho}a_j)=e^{|\tk''_{i\rho}|a_j}\hJnd(\tk_{i\rho}a_j)$ \\    
    $\Hn(\tk_{i\rho}a_j)=\widetilde{G}^{-1}_i a_j^{-n}\hHn(\tk_{i\rho}a_j)$ &
    $\Hn(\tk_{i\rho}a_j)=\widetilde{P}_{ij}^{-1}\hHn(\tk_{i\rho}a_j)$ &
    $\Hn(\tk_{i\rho}a_j)=e^{-\tk''_{i\rho}a_j}\hHn(\tk_{i\rho}a_j)$ \\    
    $\Hnd(\tk_{i\rho}a_j)=\widetilde{G}^{-1}_i a_j^{-n}\hHnd(\tk_{i\rho}a_j)$ &
    $\Hnd(\tk_{i\rho}a_j)=\widetilde{P}_{ij}^{-1}\hHnd(\tk_{i\rho}a_j)$ &
    $\Hnd(\tk_{i\rho}a_j)=e^{-\tk''_{i\rho}a_j}\hHnd(\tk_{i\rho}a_j)$ \\
    \hline
\end{tabular}
\label{ch2.1.T.rccf}
\end{center}
\end{table}
\begin{flalign}
\widetilde{G}_i=\frac{1}{n!}\left(\frac{\tk_{i\rho}}{2}\right)^n, \label{ch2.1.E.Gi}
\end{flalign}
\begin{flalign} 
\widetilde{P}_{ij}=
	\begin{cases}
		1,& \text{if } |\Jn(\tk_{i\rho} a_j)|^{-1}<T_m, \\
		|\Jn(\tk_{i\rho} a_j)|, & \text{if } |\Jn(\tk_{i\rho} a_j)|^{-1} \ge T_m,
	\end{cases} \label{ch2.1.E.Pii}
\end{flalign}
\begin{flalign}
\tk''_{i\rho}
=\Im m\left[\frac{k_{i\rho}}{\kappa_{i\epsilon}}\right]
=\Im m\left[\frac{k'_{i\rho}+\iu k''_{i\rho}}{\kappa_{i\epsilon}'+\iu \kappa_{i\epsilon}'' }\right]
=\Im m\left[\frac{\Big(k'_{i\rho}+\iu k''_{i\rho}\Big)\Big(\kappa_{i\epsilon}'-\iu \kappa_{i\epsilon}'' \Big)}
	{\Big(\kappa_{i\epsilon}'\Big)^2+\Big(\kappa_{i\epsilon}''\Big)^2 }\right]
=\frac{\kappa_{i\epsilon}' k''_{i\rho} - \kappa_{i\epsilon}'' k'_{i\rho}}
	{\left| \kappa_{i\epsilon} \right|^2}, \label{ch2.1.E.kirho}
\end{flalign}
where $T_m$ is the magnitude threshold for moderate arguments \cite{Moon14:Stable}. Note that subscripts $i$ and $j$ are arbitrary.   Range-conditioned functions for uniaxial permeability media can be similarly constructed. The multiplicative factors associated with the new functions shown in Table \ref{ch2.1.T.rccf} can be classified into two types: $\alpha$-type and $\beta$-type, whereby the relationship between the original cylindrical functions and the range-conditioned ones can be succinctly expressed as
\begin{subequations}
\begin{flalign}
\Jn(\tk_{i\rho}a_j) &= \widetilde\beta_{ij}\hJn(\tk_{i\rho}a_j), \\
\Jnd(\tk_{i\rho}a_j) &= \widetilde\beta_{ij}\hJnd(\tk_{i\rho}a_j), \\
\Hn(\tk_{i\rho}a_j) &= \widetilde\alpha_{ij}\hHn(\tk_{i\rho}a_j), \\
\Hnd(\tk_{i\rho}a_j) &= \widetilde\alpha_{ij}\hHnd(\tk_{i\rho}a_j).
\end{flalign}
\end{subequations}
\begin{table}[t]
\begin{center}
\renewcommand{\arraystretch}{1.6}
\setlength{\tabcolsep}{10pt}
\caption{Definition of $\widetilde\alpha_{ij}$ and $\widetilde\beta_{ij}$.}
\begin{tabular}{ccc}
    \hline
    Argument type & $\widetilde\alpha_{ij}$ & $\widetilde\beta_{ij}$ \\
    \hline
    Small  & $\widetilde G^{-1}_i a^{-n}_j$ & $\widetilde G_i a^n_j$ \\
    Moderate & $\widetilde P_{ij}^{-1}$ & $\widetilde P_{ij}$\\
    Large  & $e^{-\widetilde k''_{i\rho}a_j}$ & $e^{\widetilde k''_{i\rho}a_j}$\\
    \hline
\end{tabular}
\label{ch2.1.T.alpha.beta}
\end{center}
\end{table}
The definitions of $\widetilde\alpha_{ij}$ and $\widetilde\beta_{ij}$ are provided in Table \ref{ch2.1.T.alpha.beta}.
Similarly to the isotropic case, $\widetilde\alpha_{ij}$ and $\widetilde\beta_{ij}$ exhibit two important properties to ensure a stable computation~\cite{Moon14:Stable}:

1. {\it Reciprocity.} 
\begin{flalign}
\widetilde\alpha_{ii}=\widetilde\beta_{ii}^{-1}. \label{ch2.1.E.alpha.beta.1}
\end{flalign}

2. {\it Boundness.}
\begin{flalign}
|\widetilde\beta_{im}\,\widetilde\alpha_{in}|\leq 1, \quad\text{for } a_m < a_n. \label{ch2.1.E.alpha.beta.2}
\end{flalign}

For the anisotropic case, it is convenient to also derive range-conditioned cylindrical `matrices' because 2$\times$2 matrices rather than scalar factors appear in the computation of the reflection and transmission coefficients in layered media. We use hats to denote those matrices as well. From \eqref{ch1.2.E.Bzn} and \eqref{ch1.2.E.Bpn},  we obtain
\begin{subequations}
\begin{flalign}
\bJzn(k_{i\rho} a_j)
&=	\begin{bmatrix}
    \hJn(\tk_{i\rho} a_j) & 0 \\
     0 & \hJn(\dk_{i\rho}a_j) \\
    \end{bmatrix}
\cdot
\begin{bmatrix}
    \widetilde{\beta}_{ij} & 0 \\
     0 & \ddot{\beta}_{ij} \\
    \end{bmatrix}
=\bb_{ij}\cdot\hbJzn=\hbJzn\cdot\bb_{ij} ,\label{ch2.1.E.Jzn}\\
\bHzn(k_{i\rho} a_j)
&=  \begin{bmatrix}
    \hHn(\tk_{i\rho} a_j) & 0 \\
     0 & \hHn(\dk_{i\rho}a_j) \\
    \end{bmatrix}
\cdot
	\begin{bmatrix}
    \widetilde{\alpha}_{ij} & 0 \\
     0 & \ddot{\alpha}_{ij} \\
    \end{bmatrix}
=\ba_{ij}\cdot\hbHzn=\hbHzn\cdot\ba_{ij} ,\label{ch2.1.E.Hzn}\\
\bJpn(k_{i\rho}a_j)
&=\frac{1}{k_{i\rho}^2 a_j}
	\begin{bmatrix}
    \iu\omega\epsilon_{hi} \tk_{i\rho}a_j \hJnd(\tk_{i\rho}a_j) & -n k_z \hJn(\dk_{i\rho}a_j) \\
    -n k_z \hJn(\tk_{i\rho}a_j) & -\iu\omega\mu_{hi} \dk_{i\rho}a_j \hJnd(\dk_{i\rho}a_j) \\
    \end{bmatrix}
\cdot
	\begin{bmatrix}
    \widetilde{\beta}_{ij} & 0 \\
     0 & \ddot{\beta}_{ij} \\
    \end{bmatrix}
=\hbJpn\cdot\bb_{ij} ,\label{ch2.1.E.Jpn}\\
\bHpn(k_{i\rho}a_j)
&=\frac{1}{k_{i\rho}^2 a_j}
	\begin{bmatrix}
    \iu\omega\epsilon_{hi} \tk_{i\rho}a_j \hHnd(\tk_{i\rho}a_j) & -n k_z \hHn(\dk_{i\rho}a_j) \\
    -n k_z \hHn(\tk_{i\rho}a_j) & -\iu\omega\mu_{hi} \dk_{i\rho}a_j \hHnd(\dk_{i\rho}a_j) \\
    \end{bmatrix}
\cdot
	\begin{bmatrix}
    \widetilde{\alpha}_{ij} & 0 \\
     0 & \ddot{\alpha}_{ij} \\
    \end{bmatrix}
=\hbHpn\cdot\ba_{ij}. \label{ch2.1.E.Hpn}
\end{flalign}
\end{subequations}
Note that two matrices in \eqref{ch2.1.E.Jzn} and \eqref{ch2.1.E.Hzn} are diagonal, so they commute.

\subsection{Range-conditioned reflection and transmission coefficients}
\label{sec.3.2}
\begin{figure}[t]
  \centering
  \includegraphics[height=2.0in]{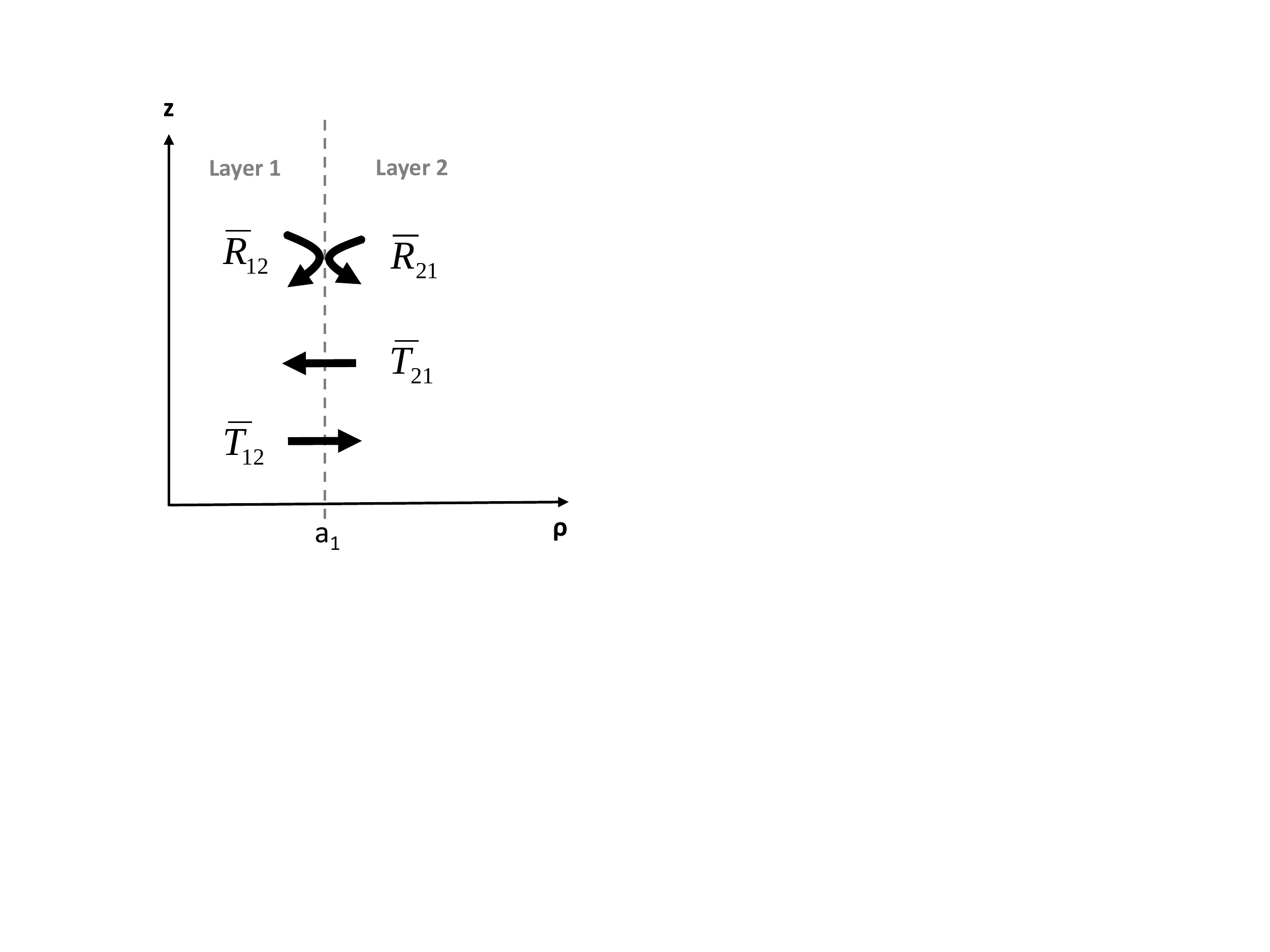}
  \caption{Reflection and transmission coefficients for two cylindrical layers.}
  \label{F.local.region12}
\end{figure}
Using the redefined matrices in  \eqref{ch2.1.E.Jzn} -- \eqref{ch2.1.E.Hpn}, local reflection and transmission coefficients for the two-layer medium depicted in Fig. \ref{F.local.region12} can be also redefined. First of all, two intermediate matrices, $\bD_A$ and $\bD_B$, are redefined as
\begin{subequations}
\begin{flalign}
\bD_A
=\left[\hbJ_{z11}\cdot\Bs{11}
	-\hbH_{z21}\cdot\As{21}\cdot\Ass{21}{-1}\cdot\hbH_{\phi 21}^{-1}\cdot\hbJ_{\phi 11}\cdot\Bs{11}
  \right]
=\left[\hbJ_{z11} - \hbH_{z21}\cdot\hbH_{\phi 21}^{-1}\cdot\hbJ_{\phi 11}
  \right]\cdot\Bs{11}
=\hbD_A\cdot\Bs{11}, \label{ch2.2.E.DA}\\
\bD_B
=\left[\hbH_{z21}\cdot\As{21}
	-\hbJ_{z11}\cdot\Bs{11}\cdot\Bss{11}{-1}\cdot\hbJ_{\phi 11}^{-1}\cdot\hbH_{\phi 21}\cdot\As{21}
  \right]
=\left[\hbH_{z21} - \hbJ_{z11}\cdot\hbJ_{\phi 11}^{-1}\cdot\hbH_{\phi 21}
  \right]\cdot\As{21}
=\hbD_B\cdot\As{21}. \label{ch2.2.E.DB}
\end{flalign}
\end{subequations}
Therefore, the local reflection and transmission coefficient are redefined as
\begin{subequations}
\begin{flalign}
\bR_{12}
&=\bD_A^{-1}\cdot
	\left[\bH_{z21}\cdot\bH_{\phi 21}^{-1}\cdot\bH_{\phi 11}-\bH_{z11}\right]
=\Bss{11}{-1}\cdot\hbD_A^{-1}\cdot
	\left[\hbH_{z21}\cdot\As{21}\cdot\Ass{21}{-1}\cdot\hbH_{\phi 21}^{-1}\cdot\hbH_{\phi 11}\cdot\As{11}
	-\hbH_{z11}\cdot\As{11}\right] \notag\\
&=\As{11}\cdot\hbD_A^{-1}\cdot
	\left[\hbH_{z21}\cdot\hbH_{\phi 21}^{-1}\cdot\hbH_{\phi 11}-\hbH_{z11}\right]\cdot\As{11}
=\As{11}\cdot\hbR_{12}\cdot\As{11}, \label{ch2.2.E.R12}
\end{flalign}
\begin{flalign}
\bR_{21}
&=\bD_B^{-1}\cdot
	\left[\bJ_{z11}\cdot\bJ_{\phi 11}^{-1}\cdot\bJ_{\phi 21} - \bJ_{z21}\right]
=\Ass{21}{-1}\cdot\hbD_B^{-1}\cdot
	\left[\hbJ_{z11}\cdot\Bs{11}\cdot\Bss{11}{-1}\cdot\hbJ_{\phi 11}^{-1}\cdot\hbJ_{\phi 21}\cdot\Bs{21}
	-\hbJ_{z21}\cdot\Bs{21}\right] \notag\\
&=\Bs{21}\cdot\hbD_B^{-1}\cdot
	\left[\hbJ_{z11}\cdot\hbJ_{\phi 11}^{-1}\cdot\hbJ_{\phi 21} - \hbJ_{z21}\right]\cdot\Bs{21}
=\Bs{21}\cdot\hbR_{21}\cdot\Bs{21}, \label{ch2.2.E.R21}
\end{flalign}
\begin{flalign}
\bT_{12}
&=\bD_B^{-1}\cdot
	\left[\bH_{z11} - \bJ_{z11}\cdot\bJ_{\phi 11}^{-1}\cdot\bH_{\phi 11}\right]
=\Ass{21}{-1}\cdot\hbD_B^{-1}\cdot
	\left[\hbH_{z11}\cdot\As{11}
		-\hbJ_{z11}\cdot\As{11}\cdot\Bss{11}{-1}\cdot\hbJ_{\phi 11}^{-1}\cdot\hbH_{\phi 11}\cdot\As{11}
	\right] \notag\\
&=\Bs{21}\cdot\hbD_B^{-1}\cdot
	\left[\hbH_{z11} - \hbJ_{z11}\cdot\hbJ_{\phi 11}^{-1}\cdot\hbH_{\phi 11}\right]\cdot\As{11}
=\Bs{21}\cdot\hbT_{12}\cdot\As{11}, \label{ch2.2.E.T12}
\end{flalign}
\begin{flalign}
\bT_{21}
&=\bD_A^{-1}\cdot
	\left[\bJ_{z21} - \bH_{z21}\cdot\bH_{\phi 21}^{-1}\cdot\bJ_{\phi 21}\right]
=\Bss{11}{-1}\cdot\hbD_A^{-1}\cdot
	\left[\hbJ_{z21}\cdot\Bs{21}
		-\hbH_{z21}\cdot\As{21}\cdot\Ass{21}{-1}\cdot\hbH_{\phi 21}^{-1}\cdot\hbJ_{\phi 21}\cdot\Bs{21}
	\right] \notag\\
&=\As{11}\cdot\hbD_A^{-1}\cdot
	\left[\hbJ_{z21} - \hbH_{z21}\cdot\hbH_{\phi 21}^{-1}\cdot\hbJ_{\phi 21}\right]\cdot\Bs{21}
=\As{11}\cdot\hbT_{21}\cdot\Bs{21}. \label{ch2.2.E.T21}
\end{flalign}
\end{subequations}

\begin{figure}[t]
	\centering
	\subfloat[\label{ch2.2.F.RF.Out}]{%
      \includegraphics[width=2.5in]{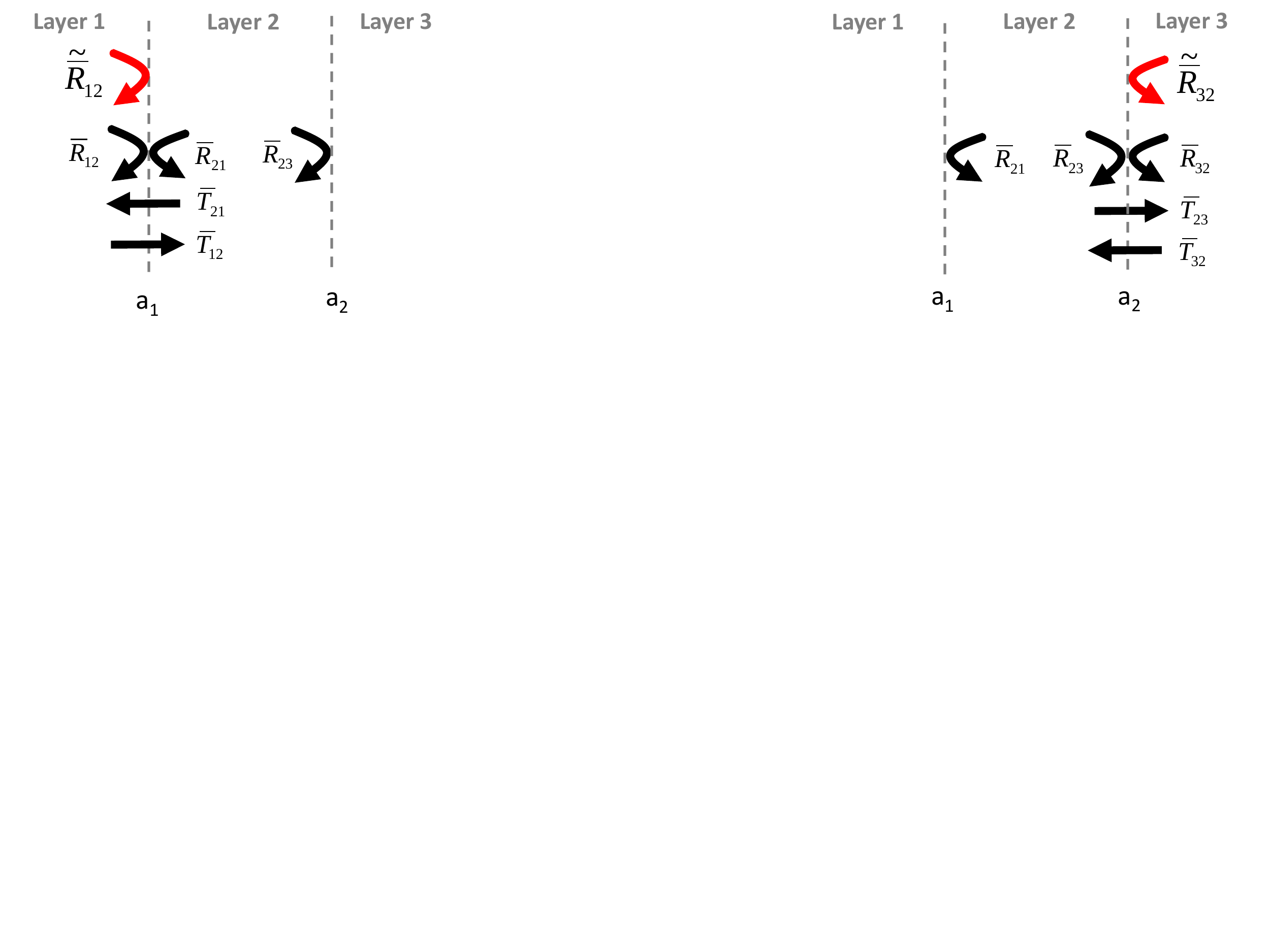}
    }
    \hspace{2.0cm}
    \subfloat[\label{ch2.2.F.RF.Stand}]{%
      \includegraphics[width=2.5in]{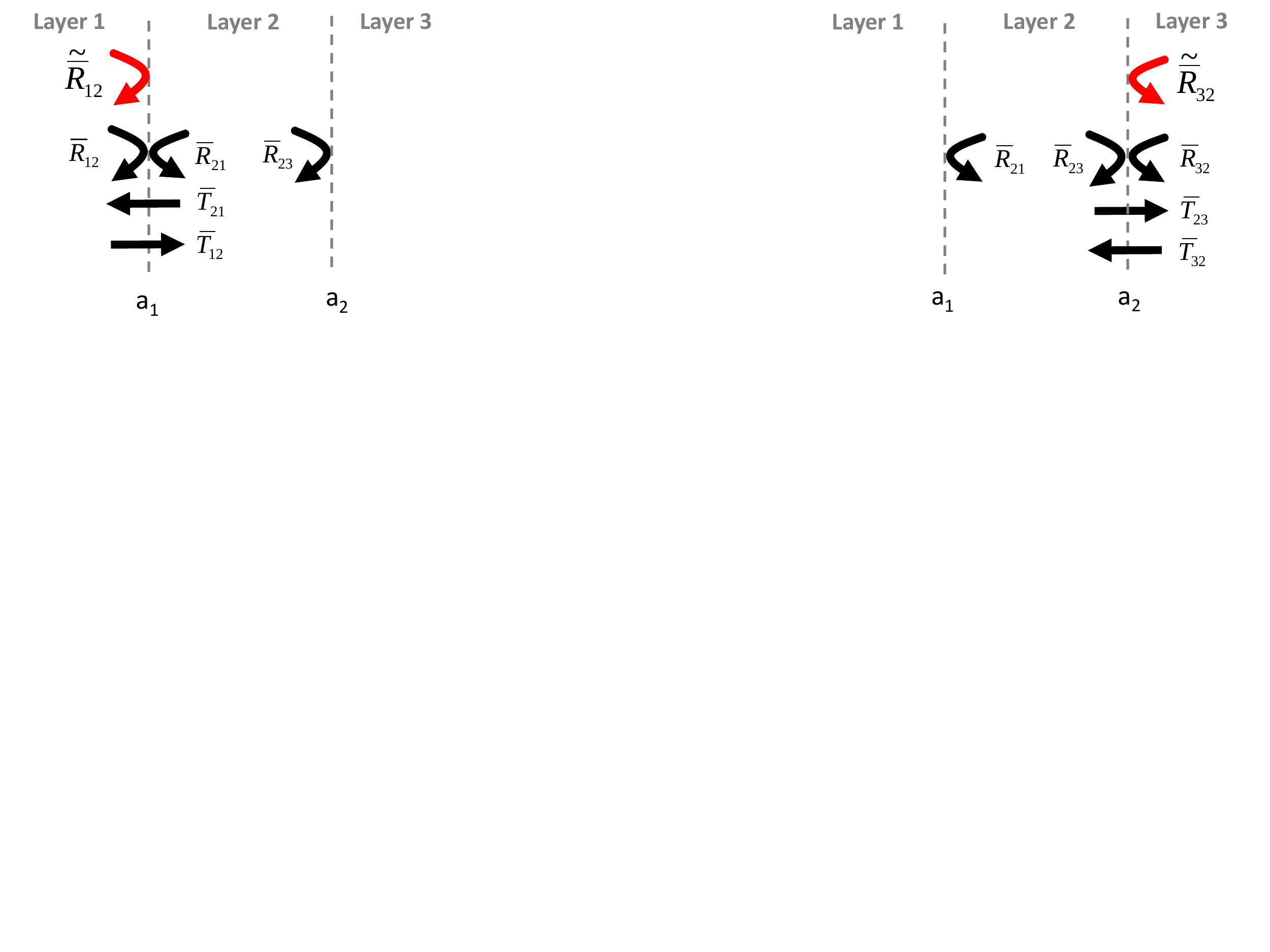}
    }
    \caption{Generalized reflection coefficients for three cylindrical layers: (a) $\tbR_{12}$ for the outgoing-wave case and (b) $\tbR_{32}$ for the standing-wave case.}
    \label{ch2.2.F.RF.3layers}
\end{figure}

We can proceed to redefine generalized reflection coefficients for multilayers, which are functions of local reflection and transmission coefficients. The generalized reflection coefficient for the outgoing-wave case for three cylindrical layers depicted in Fig. \ref{ch2.2.F.RF.Out} is modified to
\begin{flalign}
\tbR_{12}
&=\bR_{12}+\bT_{21}\cdot\bR_{23}\cdot
	\left[\bI-\bR_{21}\cdot\bR_{23}\right]^{-1}\cdot\bT_{12} \notag\\
&=\As{11}\cdot
	\left\{
		\hbR_{12}+\hbT_{21}\cdot\Bs{21}\cdot\As{22}\cdot\hbR_{23}\cdot\Bs{21}\cdot\As{22}\cdot
		\left[\bI-\hbR_{21}\cdot\Bs{21}\cdot\As{22}\cdot\hbR_{23}\cdot\Bs{21}\cdot\As{22}
		\right]^{-1}\cdot\hbT_{12}
	\right\}\cdot\As{11} \notag\\
&=\As{11}\cdot\htbR_{12}\cdot\As{11}. \label{ch2.2.E.gen.R12}
\end{flalign}
Note that the magnitude of the multiplicative factor $\Bs{21}\cdot\As{22}$ in \eqref{ch2.2.E.gen.R12} is never greater than one due to the boundness property, which guarantees moderate magnitude for $\htbR_{12}$ in any case. Since the associated multiplicative factors $\As{11}$ shown in \eqref{ch2.2.E.gen.R12} are the same as those for $\bR_{12}$ (see \eqref{ch2.2.E.R12}), they do not change when more than three layers are present. Therefore, the redefined generalized reflection coefficient between two arbitrarily-indexed adjacent layers for the outgoing-wave case can be expressed in general as
\begin{flalign}
\tbR_{i,i+1}=\As{ii}\cdot\htbR_{i,i+1}\cdot\As{ii}, \label{ch2.2.E.gen.Rij}
\end{flalign}
where
\begin{flalign}
\htbR_{i,i+1}
&=\hbR_{i,i+1}+\hbT_{i+1,i}\cdot\Bs{i+1,i}\cdot\As{i+1,i+1}\cdot
	\hbR_{i+1,i+2}\cdot\Bs{i+1,i}\cdot\As{i+1,i+1} \notag\\
&\qquad\qquad\cdot
	\left[\bI-\hbR_{i+1,i}\cdot\Bs{i+1,i}\cdot\As{i+1,i+1}\cdot\hbR_{i+1,i+2}\cdot\Bs{i+1,i}\cdot\As{i+1,i+1}
	\right]^{-1}\cdot\hbT_{i,i+1}. \label{ch2.2.E.hRij}
\end{flalign}
The generalized reflection coefficient for the standing-wave case for three cylindrical layers depicted in Fig. \ref{ch2.2.F.RF.Stand} is modified to
\begin{flalign}
\tbR_{32}
&=\bR_{32}+\bT_{23}\cdot\bR_{21}\cdot\left[\bI-\bR_{23}\cdot\bR_{21}\right]^{-1}\cdot\bT_{32} \notag\\
&=\Bs{32}\cdot
	\left\{
		\hbR_{32}+\hbT_{23}\cdot\Bs{21}\cdot\As{22}\cdot\hbR_{21}\cdot\Bs{21}\cdot\As{22}\cdot
		\left[\bI-\hbR_{23}\cdot\Bs{21}\cdot\As{22}\cdot\hbR_{21}\cdot\Bs{21}\cdot\As{22}
		\right]^{-1}\cdot\hbT_{32}
	\right\}\cdot\Bs{32} \notag\\
&=\Bs{32}\cdot\htbR_{32}\cdot\Bs{32}. \label{ch2.2.E.gen.R32}
\end{flalign}
Again, the multiplicative factor $\Bs{21}\cdot\As{22}$ in \eqref{ch2.2.E.gen.R32} is never greater than one in magnitude due to the boundness property, which also guarantees moderate magnitude for $\htbR_{32}$ in any case. Again, when more than three layers are present, the multiplicative factors in \eqref{ch2.2.E.gen.R32} are not affected. Therefore, the redefined generalized reflection coefficient between two arbitrarily-indexed adjacent layers for the standing-wave case is expressed as
\begin{flalign}
\tbR_{i+1,i}&=\Bs{i+1,i}\cdot\htbR_{i+1,i}\cdot\Bs{i+1,i}, \label{ch2.2.E.gen.Rji}
\end{flalign}
where
\begin{flalign}
\htbR_{i+1,i}
&=\hbR_{i+1,i}+\hbT_{i,i+1}\cdot\Bs{i,i-1}\cdot\As{ii}\cdot\hbR_{i,i-1}\cdot
	\Bs{i,i-1}\cdot\As{ii} \notag\\
&\qquad\qquad\cdot
	\left[\bI-\hbR_{i,i+1}\cdot\Bs{i,i-1}\cdot\As{ii}\cdot\hbR_{i,i-1}\cdot\Bs{i,i-1}\cdot\As{ii} 
	\right]^{-1}\cdot\hbT_{i+1,i}. \label{ch2.2.E.hRji}
\end{flalign}

\begin{figure}[t]
	\centering
	\subfloat[\label{ch2.2.F.S.Out}]{%
      \includegraphics[width=2.5in]{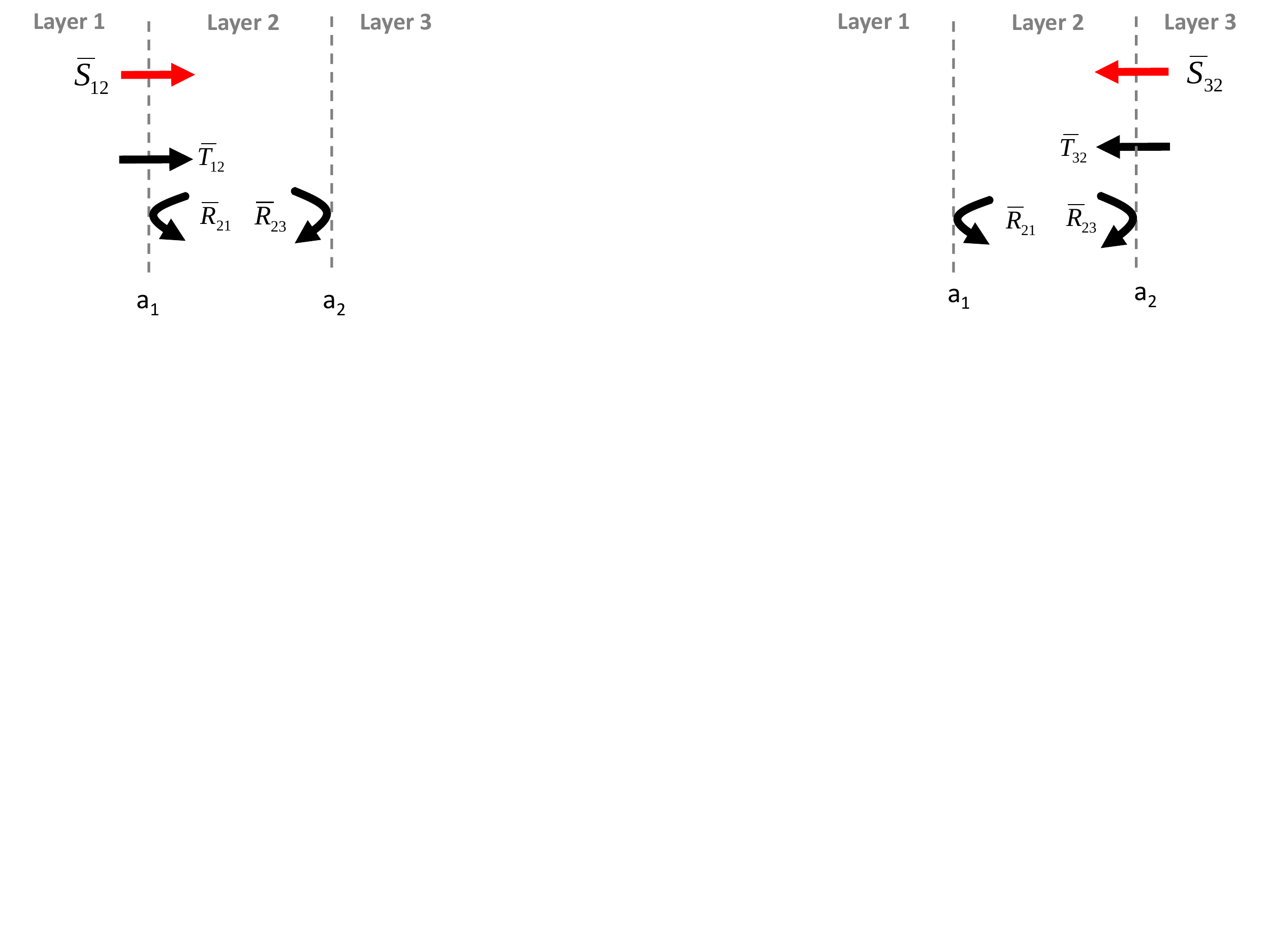}
    }
    \hspace{2.0cm}
    \subfloat[\label{ch2.2.F.S.Stand}]{%
      \includegraphics[width=2.5in]{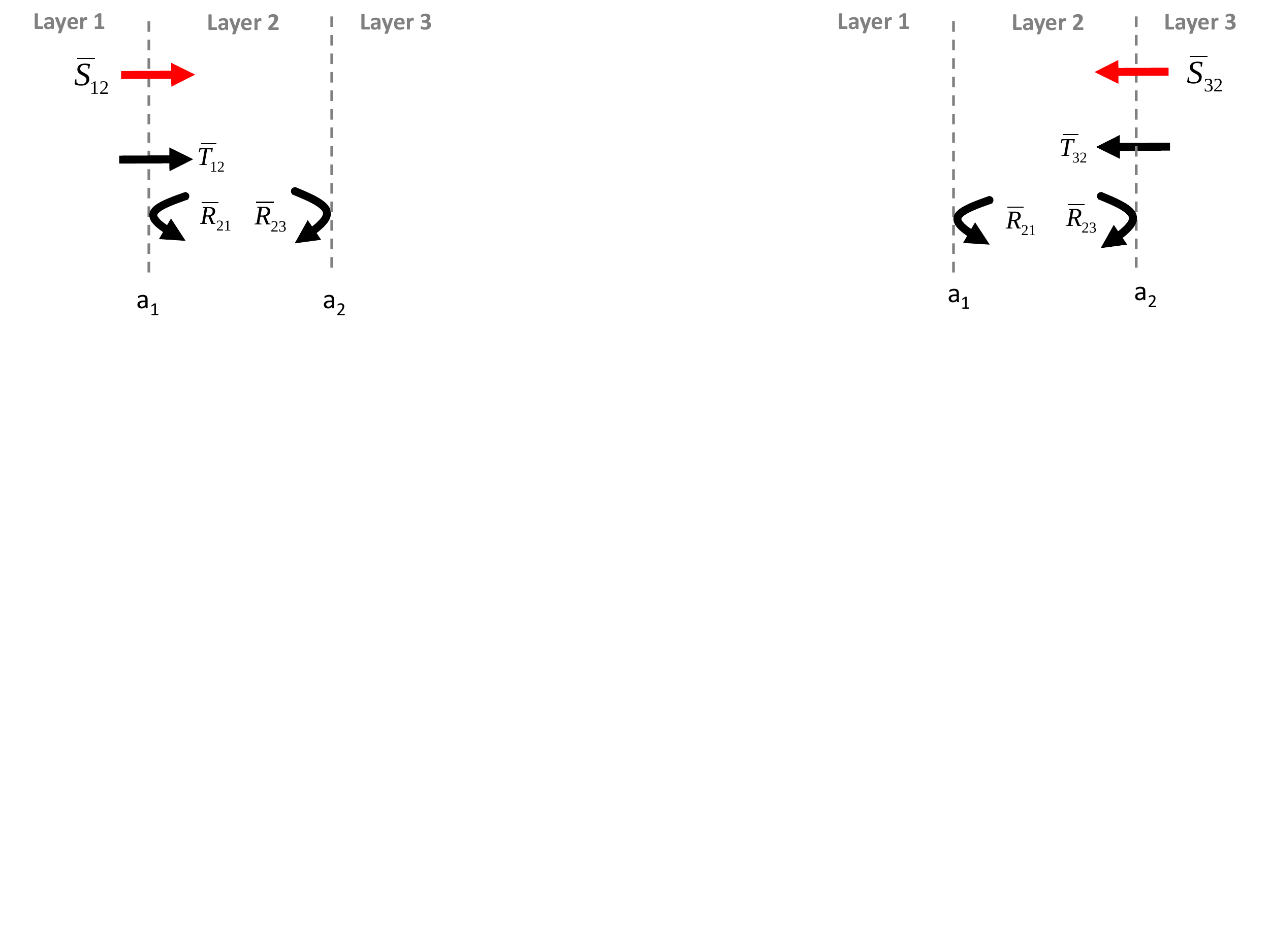}
    }
    \caption{$S$-coefficients for three cylindrical layers: (a) $\bS_{12}$ for the outgoing-wave case and (b) $\bS_{32}$ for the standing-wave case.}
    \label{ch2.2.F.S.3layers}
\end{figure}
Before we proceed to obtain generalized transmission coefficients, the redefinition of the so-called $S$-coefficients~\cite[Ch. 3]{Chew:Waves} is necessary, as they represent local transmission factors in the presence of multilayers. The $S$-coefficient for the outgoing-wave case for three cylindrical layers depicted in Fig. \ref{ch2.2.F.S.Out} is modified to
\begin{flalign}
\bS_{12}
&=\left[\bI-\bR_{21}\cdot\bR_{23}\right]^{-1}\cdot\bT_{12} \notag\\
&=\Bs{21}\cdot
	\left[\bI-\hbR_{21}\cdot\Bs{21}\cdot\As{22}\cdot\hbR_{23}\cdot\Bs{21}\cdot\As{22}
	\right]^{-1}\cdot\hbT_{12}\cdot\As{11} \notag\\
&=\Bs{21}\cdot\hbS_{12}\cdot\As{11}. \label{ch2.2.E.S12}
\end{flalign}
Therefore, the redefined arbitrarily-indexed $S$-coefficient for the outgoing-wave case is written as
\begin{flalign}
\bS_{i,i+1}=\Bs{i+1,i}\cdot\hbS_{i,i+1}\cdot\As{ii}, \label{ch2.2.E.gen.Sij}
\end{flalign}
where
\begin{flalign}
\hbS_{i,i+1}
&=\left[\bI-\hbR_{i+1,i}\cdot\Bs{i+1,i}\cdot
		\As{i+1,i+1}\cdot\hbR_{i+1,i+2}\cdot\Bs{i+1,i}\cdot\As{i+1,i+1}
	\right]^{-1}\cdot\hbT_{i,i+1}. \label{ch2.2.E.hSij}
\end{flalign}

The $S$-coefficient for the standing-wave case for three cylindrical layers depicted in Fig. \ref{ch2.2.F.S.Stand} is modified to
\begin{flalign}
\bS_{32}
&=\left[\bI-\bR_{23}\cdot\bR_{21}\right]^{-1}\cdot\bT_{32} \notag\\
&=\As{22}\cdot
	\left[\bI-\hbR_{23}\cdot\Bs{21}\cdot\As{22}\cdot\hbR_{21}\cdot\Bs{21}\cdot\As{22}
	\right]^{-1}\cdot\hbT_{32}\cdot\Bs{32} \notag\\
&=\As{22}\cdot\hbS_{32}\cdot\Bs{32}. \label{ch2.2.E.S32}
\end{flalign}
As a result, the redefined arbitrarily-indexed $S$-coefficient for the standing-wave case is written as
\begin{flalign}
\bS_{i+1,i}=\As{ii}\cdot\hbS_{i+1,i}\cdot\Bs{i+1,i}, \label{ch2.2.E.gen.Sji}
\end{flalign}
where
\begin{flalign}
\hbS_{i+1,i}
&=\left[\bI-\hbR_{i,i+1}\cdot\Bs{i,i-1}\cdot\As{ii}\cdot\hbR_{i,i-1}\cdot\Bs{i,i-1}\cdot\As{ii}
	\right]^{-1}\cdot\hbT_{i+1,i}. \label{ch2.2.E.hSji}
\end{flalign}
Let us now consider the generalized transmission coefficient for the outgoing-wave case ($i>j$) in cylindrically stratified media, which is expressed as
\begin{flalign}
\tbT_{ji}&=\bT_{i-1,i}\cdot\bS_{i-2,i-1}\cdots\bS_{j,j+1}. \label{ch2.2.E.Tji.Out.orig}
\end{flalign}
\eqref{ch2.2.E.Tji.Out.orig} can be modified in a way that
\begin{flalign}
\tbT_{ji}
&=\bT_{i-1,i}\cdot\bS_{i-2,i-1}\cdots\bS_{j,j+1} \notag\\
&=\Bs{i,i-1}\cdot\hbT_{i-1,i}\cdot
    \left(\prod_{k=j}^{i-2}\Bs{k+1,k}\cdot\As{k+1,k+1}\cdot\hbS_{k,k+1}\right)\cdot\As{jj}
\label{ch2.2.E.Tji.Out.part}\\
&=\Bs{i,i-1}\cdot\htbT_{ji}\cdot\As{jj}. \label{ch2.2.E.Tji.Out}
\end{flalign}
The magnitude of the multiplicative factors $\Bs{k+1,k}\cdot\As{k+1,k+1}$ in \eqref{ch2.2.E.Tji.Out.part} is never greater than one, which stabilizes the computation of $\htbT_{ji}$. The product in \eqref{ch2.2.E.Tji.Out.part} is the product of a number of 2$\times$2 matrices, so the order of the product should be specified. The 2$\times$2 matrix for $k=j$ and 2$\times$2 matrix for $k=i-2$ should be placed in the rightmost and leftmost in the matrix product, respectively. Furthermore, when $i=j+1$, the matrix product reduces to an identity matrix. It should be also noted that the associated multiplicative factors shown in \eqref{ch2.2.E.Tji.Out} are the generalized version of those shown in \eqref{ch2.2.E.T12}.

Next, the generalized transmission coefficient for the standing-wave case ($i<j$) in cylindrically stratified media is expressed as
\begin{flalign}
\tbT_{ji}&=\bT_{i+1,i}\cdot\bS_{i+2,i+1}\cdots\bS_{j,j-1}. \label{ch2.2.E.Tji.Stand.orig}
\end{flalign}
Similarly, \eqref{ch2.2.E.Tji.Stand.orig} is modified to
\begin{flalign}
\tbT_{ji}
&=\bT_{i+1,i}\cdot\bS_{i+2,i+1}\cdots\bS_{j,j-1} \notag\\
&=\As{ii}\cdot\hbT_{i+1,i}\cdot
    \left(\prod_{k=i+1}^{j-1}\Bs{k,k-1}\cdot\As{kk}\cdot\hbS_{k+1,k}\right)\cdot\Bs{j,j-1}
\label{ch2.2.E.Tji.Stand.part}\\
&=\As{ii}\cdot\htbT_{ji}\cdot\Bs{j,j-1}. \label{ch2.2.E.Tji.Stand}
\end{flalign}
Again, the magnitudes of the multiplicative factors $\Bs{k,k-1}\cdot\As{kk}$ in \eqref{ch2.2.E.Tji.Stand.part} are never greater than one. For the matrix product in \eqref{ch2.2.E.Tji.Stand.part}, the 2$\times$2 matrix for $k=i+1$ and 2$\times$2 matrix for $k=j-1$ should be placed in the leftmost and rightmost, which is opposite to the outgoing-wave case. Furthermore, when $j=i+1$, the matrix product reduces to an identity matrix. The associated multiplicative factors shown in \eqref{ch2.2.E.Tji.Stand} are the generalized version of those shown in \eqref{ch2.2.E.T21}.

Several auxiliary coefficients appeared in \eqref{ch1.3.E.EzHz.general.case1} -- \eqref{ch1.3.E.EzHz.general.case4} should be redefined properly as well. For the first integrand type, shown in \eqref{ch1.3.E.EzHz.general.case1},  $\tbM_{j+}$ is redefined as
\begin{flalign}
\tbM_{j+}
&=\left[\bI-\tbR_{j,j-1}\cdot\tbR_{j,j+1}\right]^{-1}
=\left[\bI-\Bs{j,j-1}\cdot\htbR_{j,j-1}\cdot\Bs{j,j-1}\cdot\As{jj}\cdot
	\htbR_{j,j+1}\cdot\As{jj}\right]^{-1} \notag\\
&=\Bs{j,[j-1,j]}\cdot
	\left[\bI-\Bs{j,j-1}\cdot\As{j,[j-1,j]}\cdot\htbR_{j,j-1}\cdot\Bs{j,j-1}\cdot
		\As{jj}\cdot\htbR_{j,j+1}\cdot\Bs{j,[j-1,j]}\cdot\As{jj}
	\right]^{-1}\cdot\As{j,[j-1,j]} \notag\\
&=\Bs{j,[j-1,j]}\cdot\htbM_{j+}\cdot\As{j,[j-1,j]}, \label{ch2.2.E.Mj.plus}
\end{flalign}
where the radial distance corresponding to subscript $[j-1,j]$ is $a_{[j-1,j]}=ca_{j-1}+(1-c)a_j$, $0\leq c\leq 1$. Two extreme choices of $a_{[j-1,j]}$ ($a_{[j-1,j]}=a_{j-1}$ and $a_{[j-1,j]}=a_j$) can be used for notational convenience but these are not useful in the redefinition of the integrand, as clarified below in Section \ref{sec.3.3}.

For the second integrand type, shown in \eqref{ch1.3.E.EzHz.general.case2}, $\tbM_{j-}$ is redefined as
\begin{flalign}
\tbM_{j-}
&=\left[\bI-\tbR_{j,j+1}\cdot\tbR_{j,j-1}\right]^{-1}
=\left[\bI-\As{jj}\cdot\htbR_{j,j+1}\cdot\As{jj}\cdot\Bs{j,j-1}\cdot
	\htbR_{j,j-1}\cdot\Bs{j,j-1}\right]^{-1} \notag\\
&=\As{j,[j-1,j]}\cdot
	\left[\bI-\Bs{j,[j-1,j]}\cdot\As{jj}\cdot\htbR_{j,j+1}\cdot\Bs{j,j-1}\cdot
		\As{jj}\cdot\htbR_{j,j-1}\cdot\Bs{j,j-1}\cdot\As{j,[j-1,j]}
	\right]^{-1}\cdot\Bs{j,[j-1,j]} \notag\\
&=\As{j,[j-1,j]}\cdot\htbM_{j-}\cdot\Bs{j,[j-1,j]}. \label{ch2.2.E.Mj.minus}
\end{flalign}
Again, the two extreme cases of $a_{[j-1,j]}$ are undesired for the proper redefinition of the integrand as shown in Section \ref{sec.3.3}.

For the third integrand type, shown in \eqref{ch1.3.E.EzHz.general.case3}, $\bN_{i+}$ is redefined as
\begin{flalign}
\bN_{i+}
&=\left[\bI-\bR_{i,i-1}\cdot\tbR_{i,i+1}\right]^{-1}
=\left[\bI-\Bs{i,i-1}\cdot\hbR_{i,i-1}\cdot\Bs{i,i-1}\cdot\As{ii}\cdot
	\htbR_{i,i+1}\cdot\As{ii}\right]^{-1} \notag\\
&=\Bs{i,i-1}\cdot
	\left[\bI-\hbR_{i,i-1}\cdot\Bs{i,i-1}\cdot\As{ii}\cdot\htbR_{i,i+1}\cdot\Bs{i,i-1}\cdot\As{ii}
	\right]^{-1}\cdot\As{i,i-1} \notag\\
&=\Bs{i,i-1}\cdot\hbN_{i+}\cdot\As{i,i-1}. \label{ch2.2.E.Ni.plus}
\end{flalign}

Finally, for the fourth integrand type, shown in \eqref{ch1.3.E.EzHz.general.case4}, $\bN_{i-}$ is redefined as
\begin{flalign}
\bN_{i-}
&=\left[\bI-\bR_{i,i+1}\cdot\tbR_{i,i-1}\right]^{-1}
=\left[\bI-\As{ii}\cdot\hbR_{i,i+1}\cdot\As{ii}\cdot\Bs{i,i-1}\cdot
	\htbR_{i,i-1}\cdot\Bs{i,i-1}\right]^{-1} \notag\\
&=\As{ii}\cdot
	\left[\bI-\hbR_{i,i+1}\cdot\Bs{i,i-1}\cdot\As{ii}\cdot\htbR_{i,i-1}\cdot\Bs{i,i-1}\cdot\As{ii}
	\right]^{-1}\cdot\Bs{ii} \notag\\
&=\As{ii}\cdot\hbN_{i-}\cdot\Bs{ii}. \label{ch2.2.E.Ni.minus}
\end{flalign}

\subsection{Range-conditioned integrand}
\label{sec.3.3}
For Case 1 in \eqref{ch1.3.E.EzHz.general.case1}, there are four arguments of interest: $k_{j\rho}a_{j-1}$, $k_{j\rho}\rho'$, $k_{j\rho}\rho$, and $k_{j\rho}a_j$. For convenience, we let $a_{j-1}=a_1$, $\rho'=a_2$, $\rho=a_3$, and $a_j=a_4$ so that $a_1<a_2<a_3<a_4$. The integrand is redefined as
\begin{flalign}
\Fn
&=\left[\bH_{zj\rho}+\bJ_{zj\rho}\cdot\tbR_{j,j+1}\right]\cdot
	\tbM_{j+}\cdot
	\left[\bJ_{zj\rho'}+\tbR_{j,j-1}\cdot\bH_{zj\rho'}\right] \notag\\
&=\left[\hbH_{zj\rho}\cdot\As{j3}+\hbJ_{zj\rho}\cdot\Bs{j3}\cdot\As{j4}\cdot\htbR_{j,j+1}\cdot\As{j4}
	\right]\cdot\Bs{j2}\cdot\htbM_{j+}\cdot\As{j2} \label{ch2.3.E.EzHz.general.case1.m}\\
&\qquad\qquad\cdot
	\left[\hbJ_{zj\rho'}\cdot\Bs{j2}+\Bs{j1}\cdot\htbR_{j,j-1}\cdot\Bs{j1}\cdot\hbH_{zj\rho'}\cdot\As{j2}
	\right] \notag\\
&=\left[\Bs{j2}\cdot\As{j3}\cdot\hbH_{zj\rho}
	+\Bs{j3}\cdot\As{j4}\cdot\hbJ_{zj\rho}\cdot\htbR_{j,j+1}\cdot\Bs{j2}\cdot\As{j4}
  \right]\cdot\htbM_{j+} \notag\\
&\qquad\qquad\cdot
	\left[\hbJ_{zj\rho'}
		+\Bs{j1}\cdot\As{j2}\cdot\htbR_{j,j-1}\cdot\hbH_{zj\rho'}\cdot\Bs{j1}\cdot\As{j2}
	\right]. \label{ch2.3.E.EzHz.general.case1}
\end{flalign}
Note that, as \eqref{ch2.3.E.EzHz.general.case1.m} shows, the corresponding radial distance to the subscript $[j-1,j]$ in \eqref{ch2.2.E.Mj.plus} is chosen to be $a_2$, neither $a_1$ nor $a_4$. To be more specific, the radial distance of the source $\rho'$ is selected. The choice enables the left and right squared bracket factors in \eqref{ch2.3.E.EzHz.general.case1} to be balanced and yields a stable computation.

For Case 2 in \eqref{ch1.3.E.EzHz.general.case2}, four arguments are of interest: $k_{j\rho}a_{j-1}$, $k_{j\rho}\rho$, $k_{j\rho}\rho'$, and $k_{j\rho}a_j$. Similarly, we let $a_{j-1}=a_1$, $\rho=a_2$, $\rho'=a_3$, and $a_j=a_4$ so that $a_1<a_2<a_3<a_4$. The integrand is redefined as
\begin{flalign}
\Fn
&=\left[\bJ_{zj\rho}+\bH_{zj\rho}\cdot\tbR_{j,j-1}\right]\cdot
	\tbM_{j-}\cdot
	\left[\bH_{zj\rho'}+\tbR_{j,j+1}\cdot\bJ_{zj\rho'}\right] \notag\\
&=\left[\hbJ_{zj\rho}\cdot\Bs{j2}+\hbH_{zj\rho}\cdot\As{j2}\cdot\Bs{j1}\cdot\htbR_{j,j-1}\cdot\Bs{j1}
	\right]\cdot\As{j3}\cdot\htbM_{j-}\cdot\Bs{j3} \label{ch2.3.E.EzHz.general.case2.m}\\
&\qquad\qquad\cdot
	\left[\hbH_{zj\rho'}\cdot\As{j3}+\As{j4}\cdot\htbR_{j,j+1}\cdot\As{j4}\cdot\hbJ_{zj\rho'}\cdot\Bs{j3}
	\right] \notag\\
&=\left[\Bs{j2}\cdot\As{j3}\cdot\hbJ_{zj\rho}
	+\Bs{j1}\cdot\As{j2}\cdot\hbH_{zj\rho}\cdot\htbR_{j,j-1}\cdot\Bs{j1}\cdot\As{j3}
  \right]\cdot\htbM_{j-} \notag\\
&\qquad\qquad\cdot
	\left[\hbH_{zj\rho'}
		+\Bs{j3}\cdot\As{j4}\cdot\htbR_{j,j+1}\cdot\hbJ_{zj\rho'}\cdot\Bs{j3}\cdot\As{j4}
	\right]. \label{ch2.3.E.EzHz.general.case2}
\end{flalign}
It should be noted that, as \eqref{ch2.3.E.EzHz.general.case2.m} shows, the corresponding radial distance to the subscript $[j-1,j]$ in \eqref{ch2.2.E.Mj.minus} is chosen to be $a_3$, the radial distance of the source (neither $a_1$ nor $a_4$). Again, this choice enables the left and right squared bracket factors in \eqref{ch2.3.E.EzHz.general.case2} to be balanced and yields a stable computation.

For Case 3 in \eqref{ch1.3.E.EzHz.general.case3}, there are 6 arguments of interest: $k_{i\rho}a_{i-1}$, $k_{i\rho}\rho$, $k_{i\rho}a_i$, $k_{j\rho}a_{j-1}$, $k_{j\rho}\rho'$, and $k_{j\rho}a_j$. We let $a_{i-1}=a_1$, $\rho=a_2$, $a_i=a_3$, $a_{j-1}=b_1$, $\rho'=b_2$, and $a_j=b_3$ so that $a_1<a_2<a_3$ and $b_1<b_2<b_3$. The integrand is redefined as
\begin{flalign}
\Fn
&=\left[\bH_{zi\rho}+\bJ_{zi\rho}\cdot\tbR_{i,i+1}\right]\cdot
	\bN_{i+}\cdot\tbT_{ji}\cdot\tbM_{j+}\cdot
	\left[\bJ_{zj\rho'}+\tbR_{j,j-1}\cdot\bH_{zj\rho'}\right] \notag\\
&=\left[\hbH_{zi\rho}\cdot\As{i2}+\hbJ_{zi\rho}\cdot\Bs{i2}\cdot\As{i3}\cdot\htbR_{i,i+1}\cdot\As{i3}
  \right] \notag\\
&\qquad\qquad\cdot
	\left(\Bs{i1}\cdot\hbN_{i+}\cdot\As{i1}\right)\cdot
    \left(\Bs{i1}\cdot\htbT_{ji}\cdot\As{j3}\right)\cdot
	\left(\Bs{j2}\cdot\htbM_{j+}\cdot\As{j2}\right) \label{ch2.3.E.EzHz.general.case3.m}\\
&\qquad\qquad\qquad\cdot
	\left[\hbJ_{zj\rho'}\cdot\Bs{j2}+\Bs{j1}\cdot\htbR_{j,j-1}\cdot
		\Bs{j1}\cdot\hbH_{zj\rho'}\cdot\As{j2}
	\right] \notag\\
&=\left[\Bs{i1}\cdot\As{i2}\cdot\hbH_{zi\rho}
		+\Bs{i2}\cdot\As{i3}\cdot\hbJ_{zi\rho}\cdot\htbR_{i,i+1}\cdot\Bs{i1}\cdot\As{i3}
  \right]\cdot\hbN_{i+}\cdot\htbT_{ji}\cdot\Bs{j2}\cdot\As{j3}\cdot\htbM_{j+} \notag\\
&\qquad\qquad\qquad\cdot
	\left[\hbJ_{zj\rho'}
		+\Bs{j1}\cdot\As{j2}\cdot\htbR_{j,j-1}\cdot\hbH_{zj\rho'}\cdot\Bs{j1}\cdot\As{j2}
	\right]. \label{ch2.3.E.EzHz.general.case3}
\end{flalign}
In should be stressed that the corresponding radial distance for $\htbM_{j+}$ in \eqref{ch2.3.E.EzHz.general.case3.m} is now $b_2$, which is the radial distance of the source.

For Case 4 in \eqref{ch1.3.E.EzHz.general.case4}, the arguments of interest are the same as those for Case 3. The integrand is redefined as
\begin{flalign}
\Fn
&=\left[\bJ_{zi\rho}+\bH_{zi\rho}\cdot\tbR_{i,i-1}\right]\cdot
	\bN_{i-}\cdot\tbT_{ji}\cdot\tbM_{j-}\cdot	
	\left[\bH_{zj\rho'}+\tbR_{j,j+1}\cdot\bJ_{zj\rho'}\right] \notag\\
&=\left[\hbJ_{zi\rho}\cdot\Bs{i2}+\hbH_{zi\rho}\cdot\As{i2}\cdot\Bs{i1}\cdot\htbR_{i,i-1}\cdot\Bs{i1}
  \right] \notag\\
&\qquad\qquad\cdot
	\left(\As{i3}\cdot\hbN_{i-}\cdot\Bs{i3}\right)\cdot
	\left(\As{i3}\cdot\htbT_{ji}\cdot\Bs{j1}\right)\cdot
	\left(\As{j2}\cdot\htbM_{j-}\cdot\Bs{j2}\right) \label{ch2.3.E.EzHz.general.case4.m}\\
&\qquad\qquad\qquad\cdot
	\left[\hbH_{zj\rho'}\cdot\As{j2}+\As{j3}\cdot\htbR_{j,j+1}\cdot
		\As{j3}\cdot\hbJ_{zj\rho'}\cdot\Bs{j2}
	\right] \notag\\	
&=\left[\Bs{i2}\cdot\As{i3}\cdot\hbJ_{zi\rho}
		+\Bs{i1}\cdot\As{i2}\cdot\hbH_{zi\rho}\cdot\htbR_{i,i-1}\cdot\Bs{i1}\cdot\As{i3}
  \right]\cdot\hbN_{i-}\cdot\htbT_{ji}\cdot\Bs{j1}\cdot\As{j2}\cdot\htbM_{j-} \notag\\
&\qquad\qquad\qquad\cdot
	\left[\hbH_{zj\rho'}
		+\Bs{j2}\cdot\As{j3}\cdot\htbR_{j,j+1}\cdot\hbJ_{zj\rho'}\cdot\Bs{j2}\cdot\As{j3}
	\right]. \label{ch2.3.E.EzHz.general.case4}
\end{flalign}
Again, the radial distance of the source $b_2$ is chosen for the corresponding radial distance for $\htbM_{j-}$ in \eqref{ch2.3.E.EzHz.general.case4.m}.

\subsection{Azimuth modal summation}
\label{sec.3.4}
The spectral representations of electromagnetic fields involve an infinite series as the azimuthal summation. The three components of the electromagnetic fields are expressed as \cite{Moon14:Stable}
\begin{subequations}
\begin{flalign}
    \begin{bmatrix}
    E_z \\ H_z
    \end{bmatrix}
&=\frac{\iu Il}{4\pi\omega\epsilon_{hj}}
    \intmp dk_ze^{\iu k_z(z-z')} \left[\suma e^{\iu n(\phi-\phi')}\Fn\cdot\Dj\right], \label{ch3.1.EH.z.compnts}\\
	\begin{bmatrix}
    E_\rho \\ H_\rho
    \end{bmatrix}
&=\frac{\iu Il}{4\pi\omega\epsilon_{hj}}\intmp dk_z e^{\iu k_z(z-z')}\frac{1}{k^2_{\rho}}
    \left[\suma e^{\iu n(\phi-\phi')}\bBn\cdot\bLn(\rho)\cdot\bMn\cdot
        \bRn(\rho')\cdot\Dj\right], \label{ch3.2.EH.rho.compnts}\\
    \begin{bmatrix}
    E_\phi \\ H_\phi
    \end{bmatrix}
&=\frac{\iu Il}{4\pi\omega\epsilon_{hj}}\intmp dk_z e^{\iu k_z(z-z')}\frac{1}{k^2_{\rho}}
    \left[\suma e^{\iu n(\phi-\phi')}\bCn\cdot\bLn(\rho)\cdot\bMn\cdot
        \bRn(\rho')\cdot\Dj\right]. \label{ch4.3.EH.phi.compnts}
\end{flalign}
\end{subequations}
 To expedite the computation, it is possible to fold the series above by exploiting symmetries of the cylindrical eigenfunctions so that only zero and positive orders remain. The expressions of folded sums are quite similar to those for isotropic media in \cite{Moon14:Stable}, so final results are not provided here. In addition, the spectral integrals \eqref{ch3.1.EH.z.compnts}, \eqref{ch3.2.EH.rho.compnts}, and \eqref{ch4.3.EH.phi.compnts} typically cannot be computed along the real axis in a robust fashion. These integrals should instead be numerically evaluated along suitably deformed integration paths in the complex $k_z$ plane. To this end, either the so-called Sommerfeld integration path (SIP) and the deformed SIP (DSIP) can be used with the optimal choice among the two, depending on the longitudinal distance between source and observation points~\cite{Moon14:Stable}.
\textcolor{\Cblue}{
Moreover, it should be noted that the concept of the direct field subtraction when source and observation points are in the same layer, discussed in \cite{Moon14:Stable}, still applies to uniaxial media. As analytical field expressions in such media are available only in certain cases (see \ref{app1}), direct field terms for isotropic media are used in the computation.
}

\section{Numerical validation results}
\label{sec.4.results}
This section provides some validation results for the formulations detailed above. In Section \ref{sec.4.1}, the results are compared against closed-form analytical solutions due to point dipole sources, available in homogeneous uniaxial media. In Section \ref{sec.4.2}, the results are compared against Finite Element Method (FEM) results for several selected cases of practical interest in geophysical exploration. 
These results also examine the effect of anisotropy ratios in the surrounding cylindrical layers on the resulting electromagnetic fields.
Throughout this section, field values are expressed in a phasor form under $e^{\iu\omega t}$ convention.

\subsection{Homogeneous uniaxial media}
\label{sec.4.1}
Numerical results from the present algorithm are compared to closed-form field expressions due to point dipole sources in homogeneous uniaxial media. For a derivation of such analytical solutions, refer to \ref{app1}.
The fields are evaluated throughout a square region of 10 cm $\times$ 10 cm in the $\rho z$-plane. The source is a $z$-directed Hertzian electric dipole with unit dipole moment and operating frequency of 36 kHz. The medium has $\epsilon_{p,h}=16\epsilon_0$ [F/m], $\mu_h=16\mu_0$ [H/m], $\sigma_h=16$ [S/m], where $\epsilon_0$ and $\mu_0$ denote free-space permittivity and permeability values. These horizontal values are fixed whereas different $\epsilon_{p,v}$, $\mu_v$, and $\sigma_v$ values are considered to yield different anisotropy ratios.

\begin{figure}[t]
	\centering
	\subfloat[\label{ch4.F.4k.10.1000}]{%
      \includegraphics[width=3.0in]{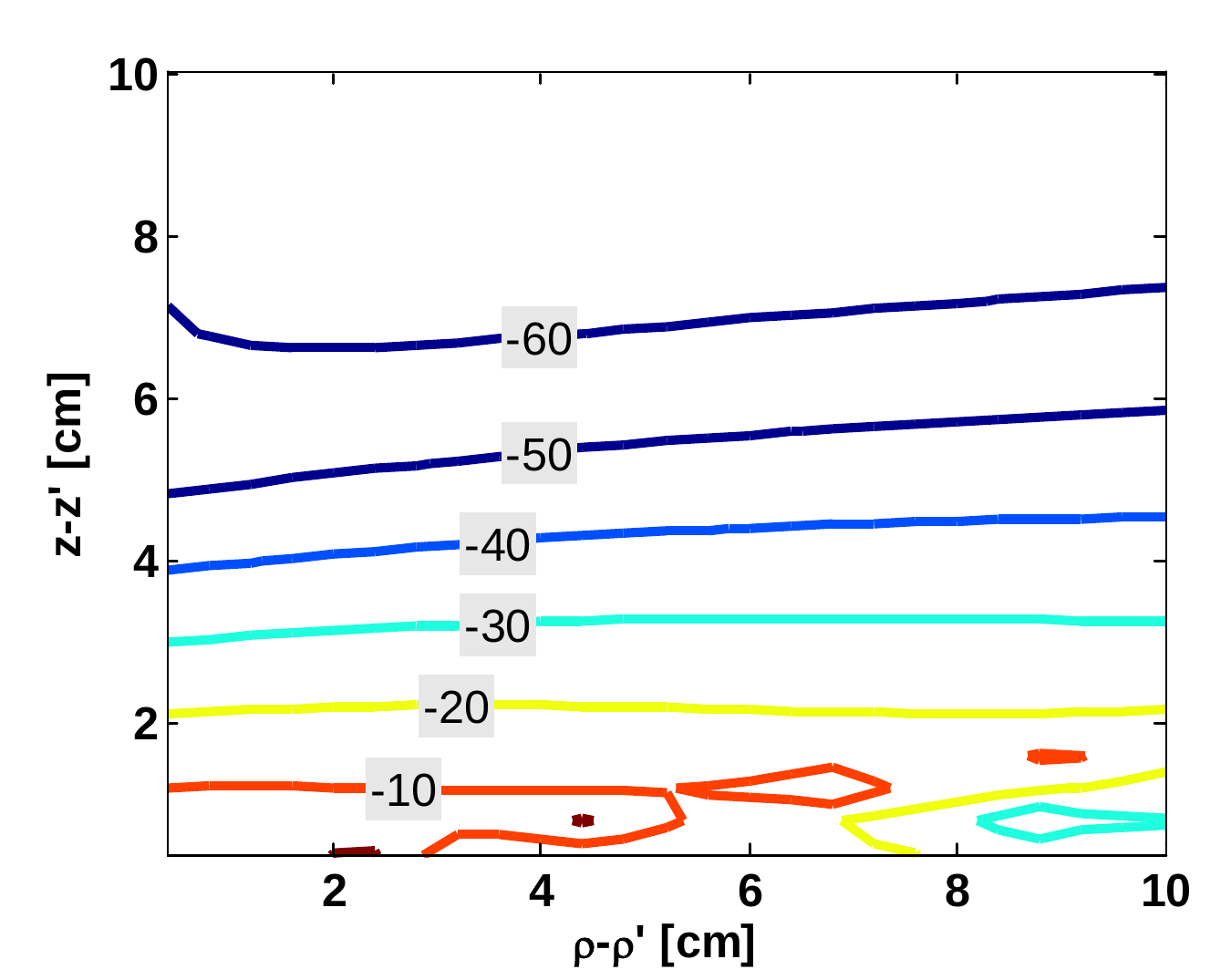}
    }
    \hfill
    \subfloat[\label{ch4.F.4k.20.1000}]{%
      \includegraphics[width=3.0in]{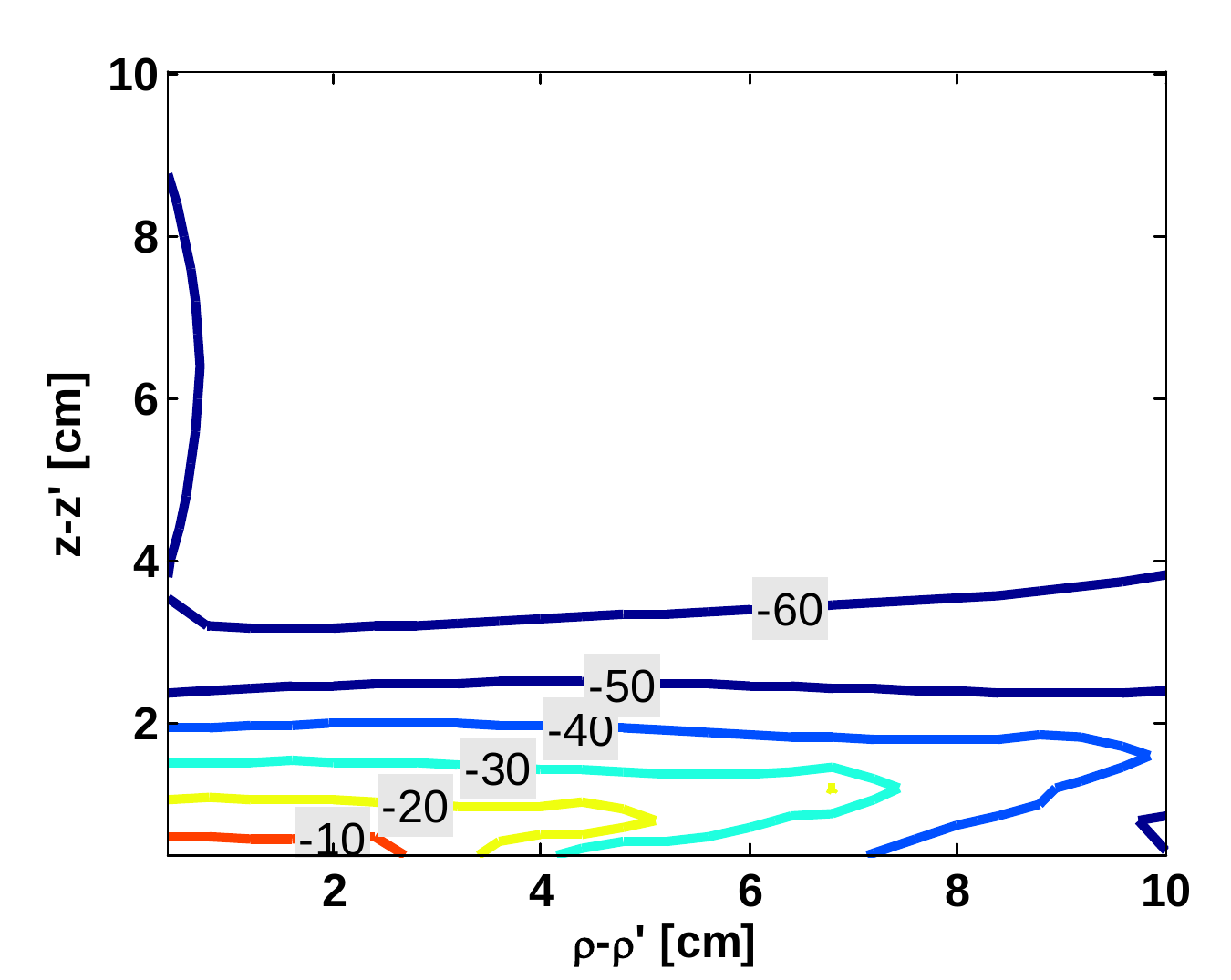}
    }\\
    \subfloat[\label{ch4.F.4k.10.2000}]{%
      \includegraphics[width=3.0in]{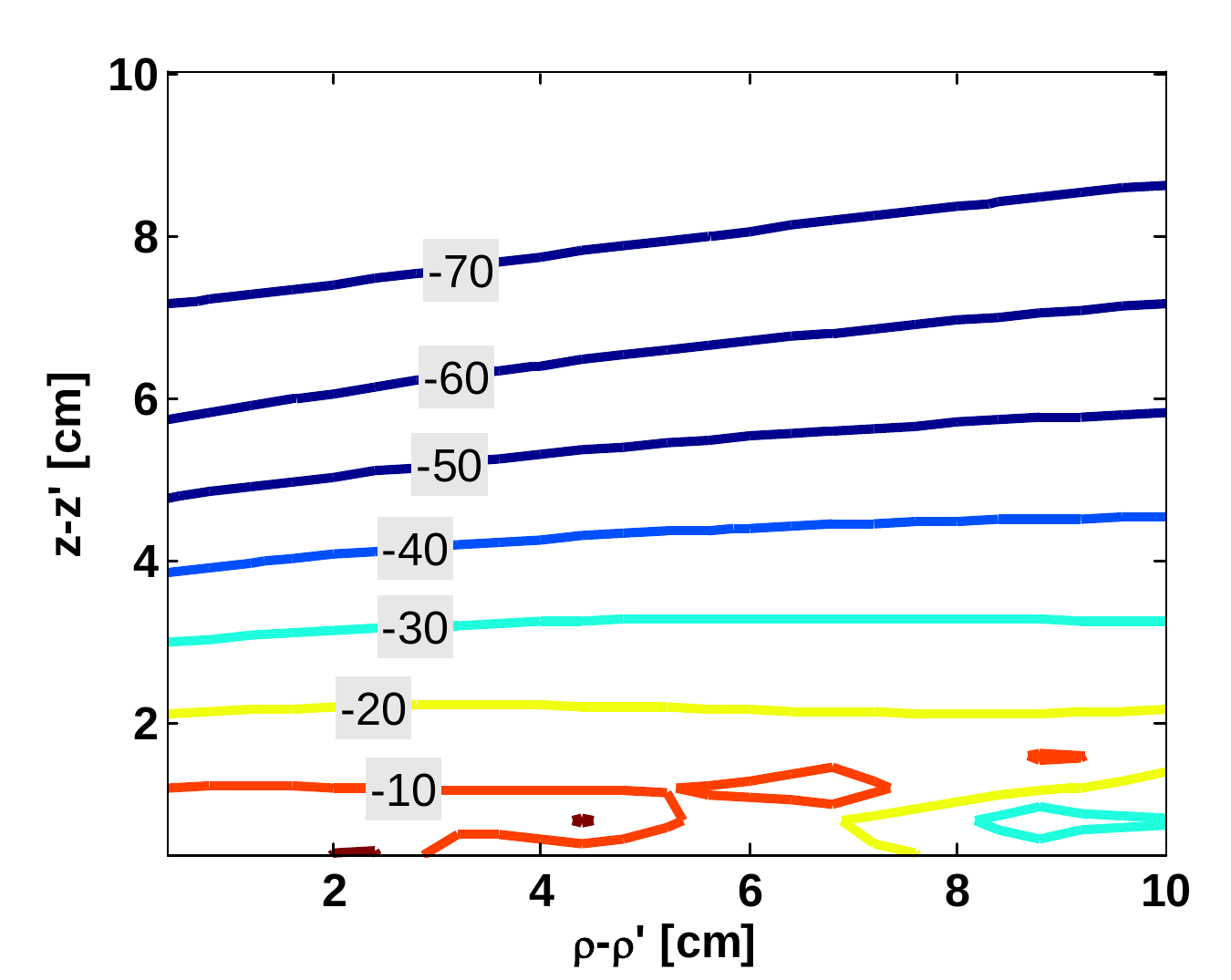}
    }
    \hfill
    \subfloat[\label{ch4.F.4k.20.2000}]{%
      \includegraphics[width=3.0in]{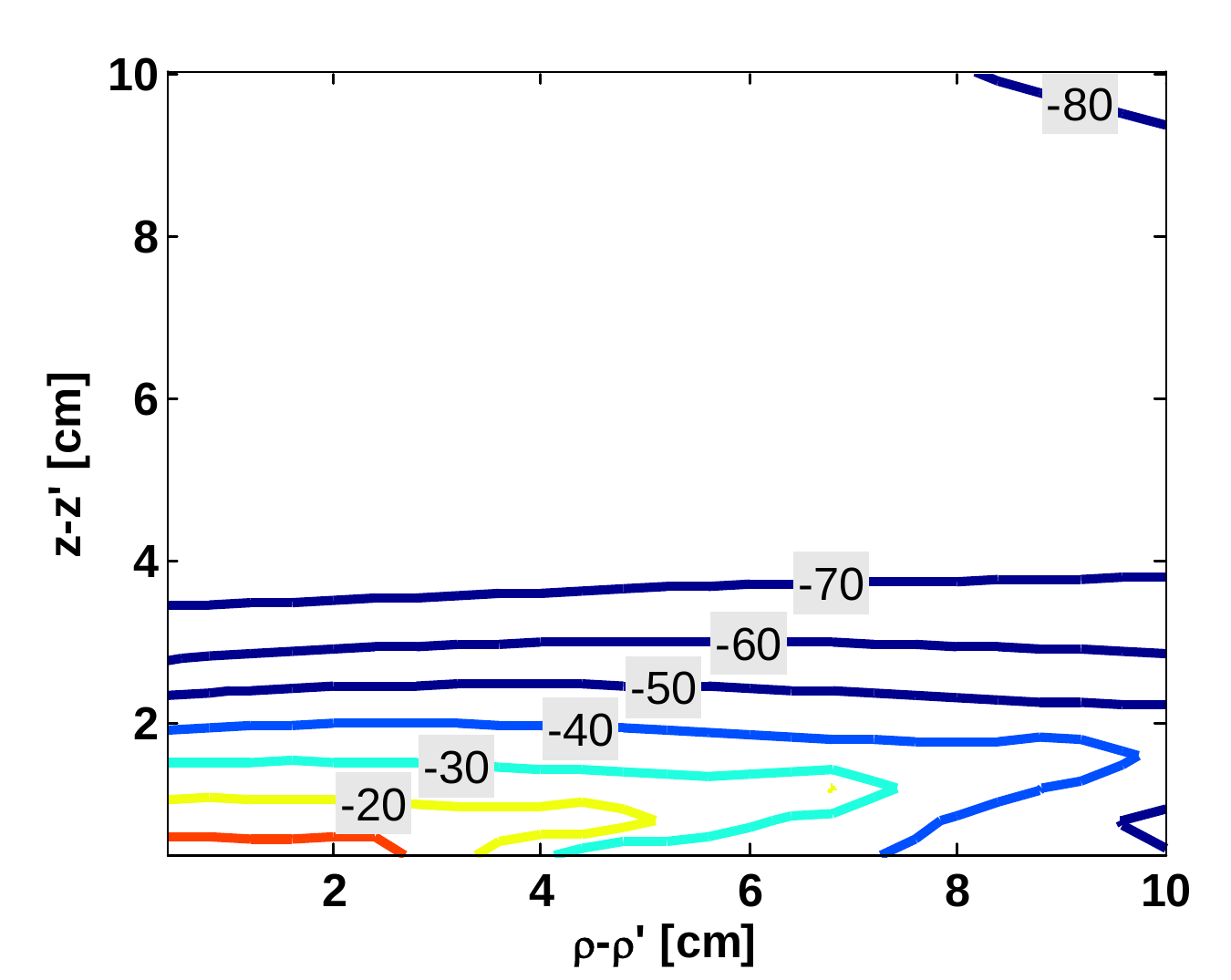}
	}    
    \caption{Relative error distribution with $\kappa=4$: (a) $n_{max}=10$, $n_{int}=1000$, (b) $n_{max}=20$, $n_{int}=1000$, (c) $n_{max}=10$, $n_{int}=2000$, and (d) $n_{max}=20$, $n_{int}=2000$.}
    \label{ch4.F.error.4k}
\end{figure}

Fig. \ref{ch4.F.4k.10.1000} -- \ref{ch4.F.4k.20.2000} show the relative error between the present algorithm and the closed-form analytical solution for different maximum orders $n_{max}$ employed in the azimuth summation and for various numbers of quadrature points $n_{int}$ employed in the numerical integration. It is assumed that $\epsilon_{p,v}=\epsilon_0$ [F/m], $\mu_v=\mu_0$ [H/m], and $\sigma_v=1$ [S/m], with $\kappa_\epsilon=\kappa_\mu=\kappa=4$. The relative error is defined as
\begin{flalign}
\text{relative error}_{dB} = 10\log_{10}\frac{|E_{z,a}-E_{z,n}|}{|E_{z,a}|},
\label{ch4.E.relative.error}
\end{flalign}
where $E_{z,a}$ and $E_{z,n}$ indicate analytical and numerical results, respectively. As expected, smaller relative errors 
are obtained for larger number of quadrature points or summation terms.
\begin{figure}[t]
	\centering
	\subfloat[\label{ch4.F.2k.10.1000}]{%
      \includegraphics[width=3.0in]{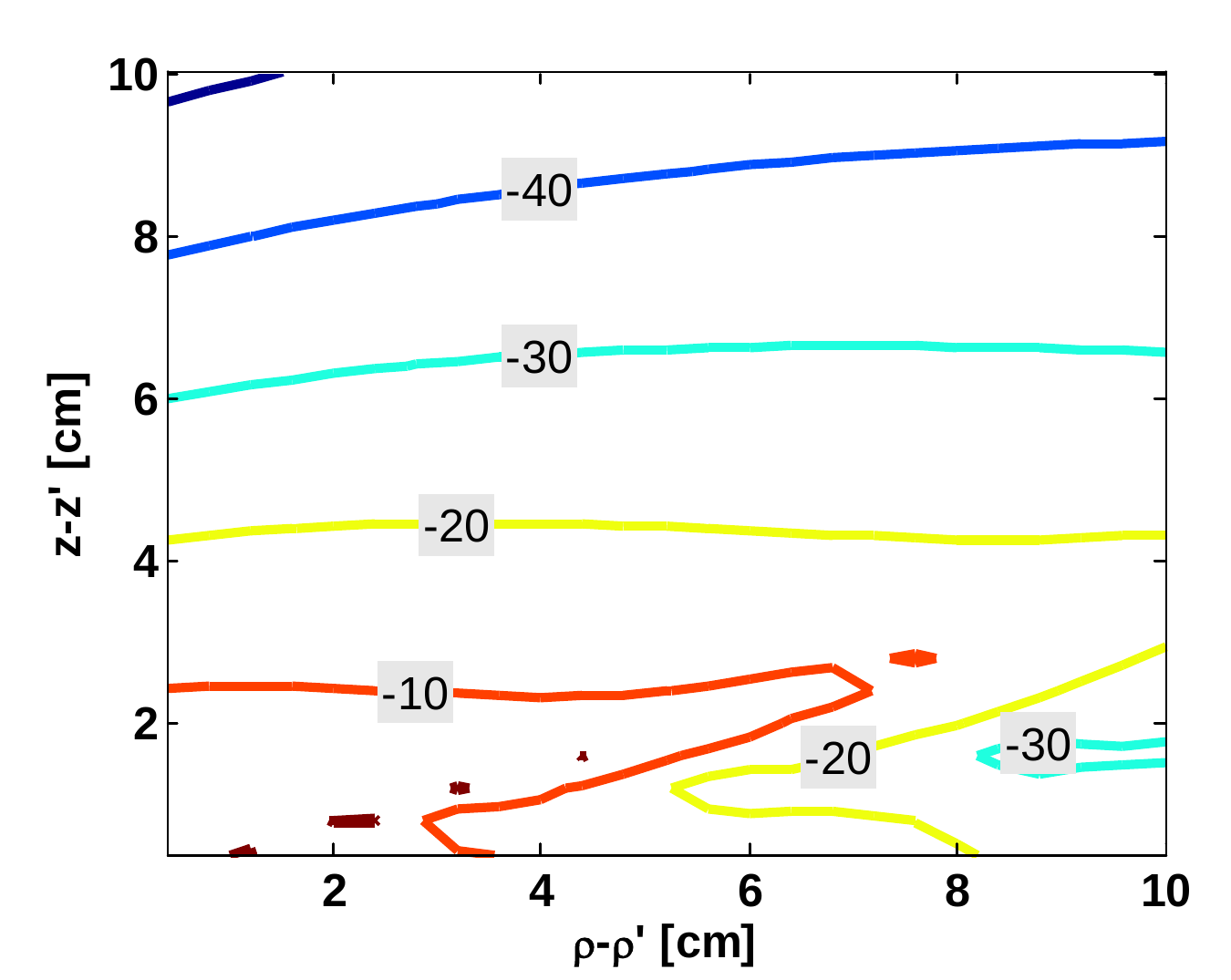}
    }
    \hfill
    \subfloat[\label{ch4.F.2k.20.1000}]{%
      \includegraphics[width=3.0in]{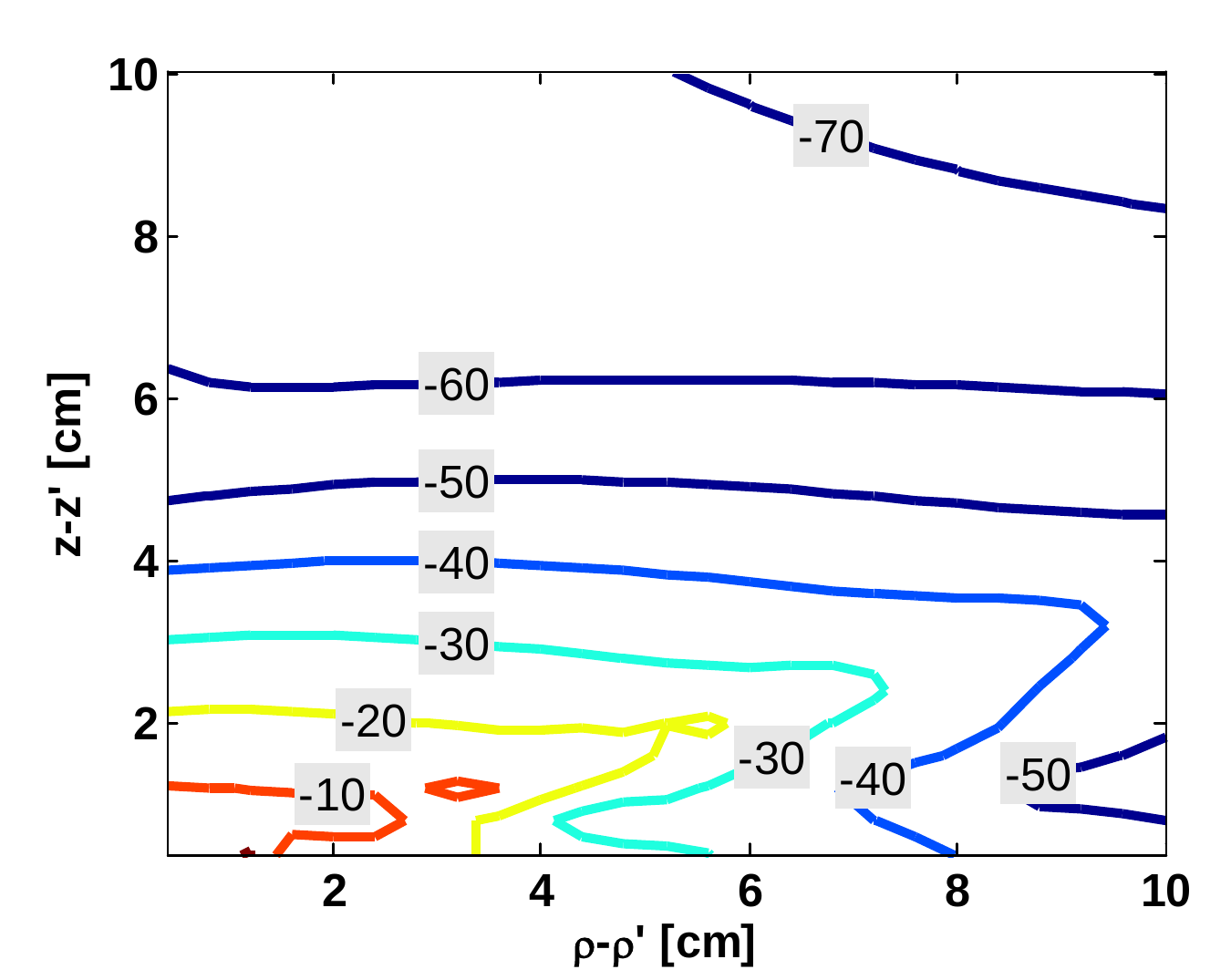}
    }\\
    \subfloat[\label{ch4.F.2k.10.2000}]{%
      \includegraphics[width=3.0in]{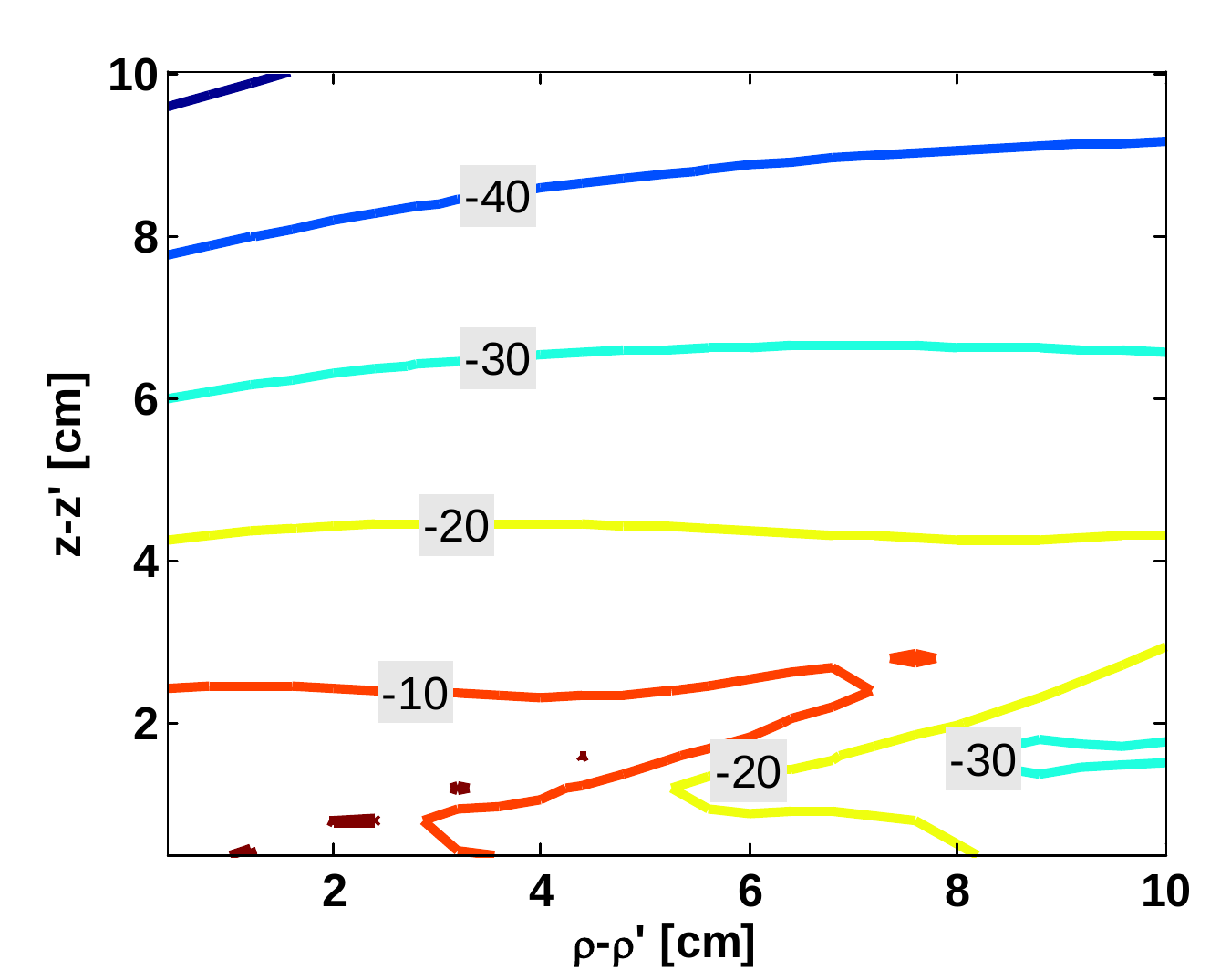}
    }
    \hfill
    \subfloat[\label{ch4.F.2k.20.2000}]{%
      \includegraphics[width=3.0in]{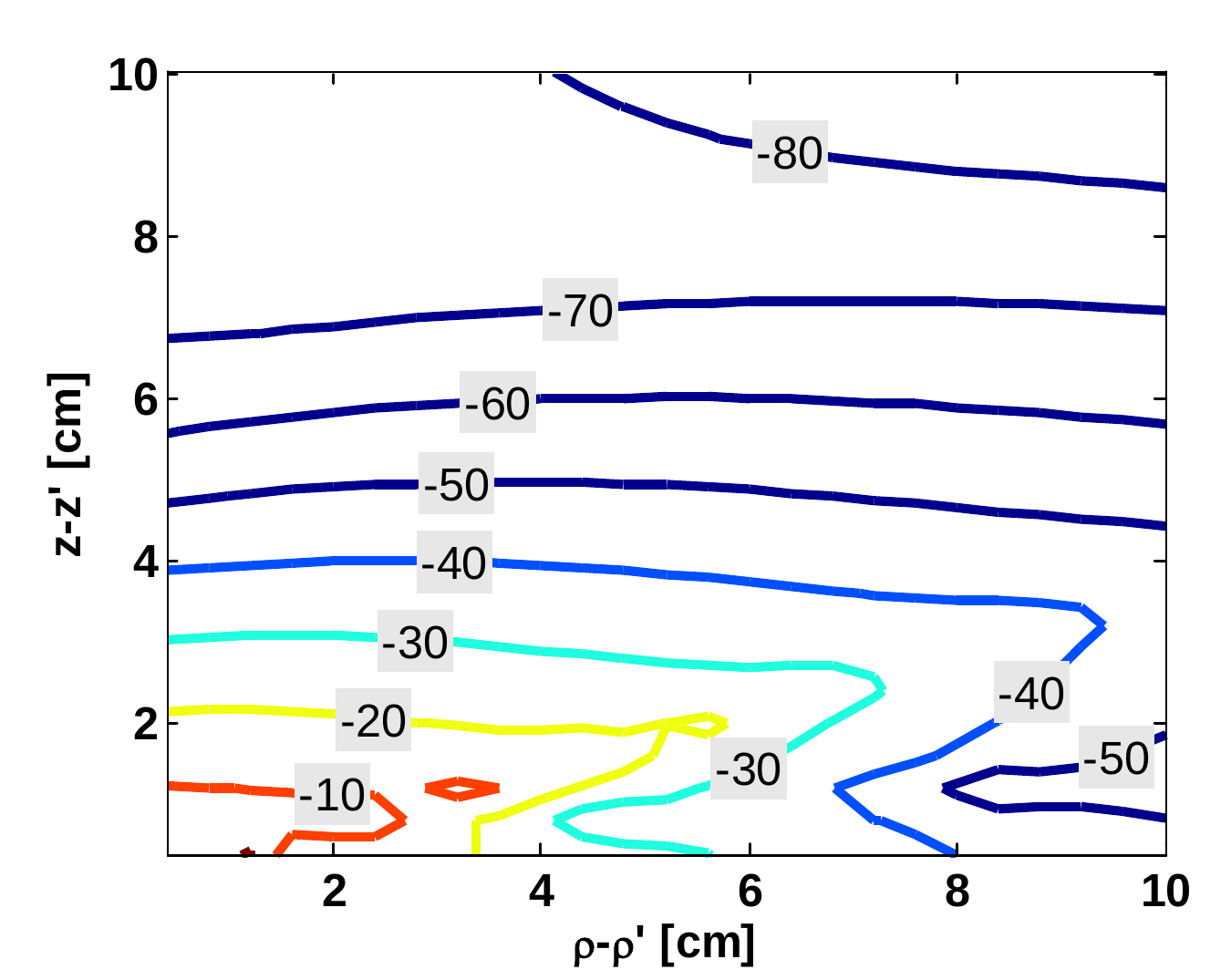}
	}    
    \caption{Relative error distribution with $\kappa=2$: (a) $n_{max}=10$, $n_{int}=1000$, (b) $n_{max}=20$, $n_{int}=1000$, (c) $n_{max}=10$, $n_{int}=2000$, and (d) $n_{max}=20$, $n_{int}=2000$.}
    \label{ch4.F.error.2k}
\end{figure}
Fig. \ref{ch4.F.2k.10.1000} -- \ref{ch4.F.2k.20.2000} show the relative error distribution for $\kappa_\epsilon=\kappa_\mu=\kappa=2$, under the assumption of $\epsilon_{p,v}=4\epsilon_0$ [F/m], $\mu_v=4\mu_0$ [H/m], and $\sigma_v=4$ [S/m]. As expected, higher $n_{max}$ and $n_{int}$ produce smaller relative errors.
\begin{figure}[t]
	\centering
	\subfloat[\label{ch4.F.sqrt2k.10.1000}]{%
      \includegraphics[width=3.0in]{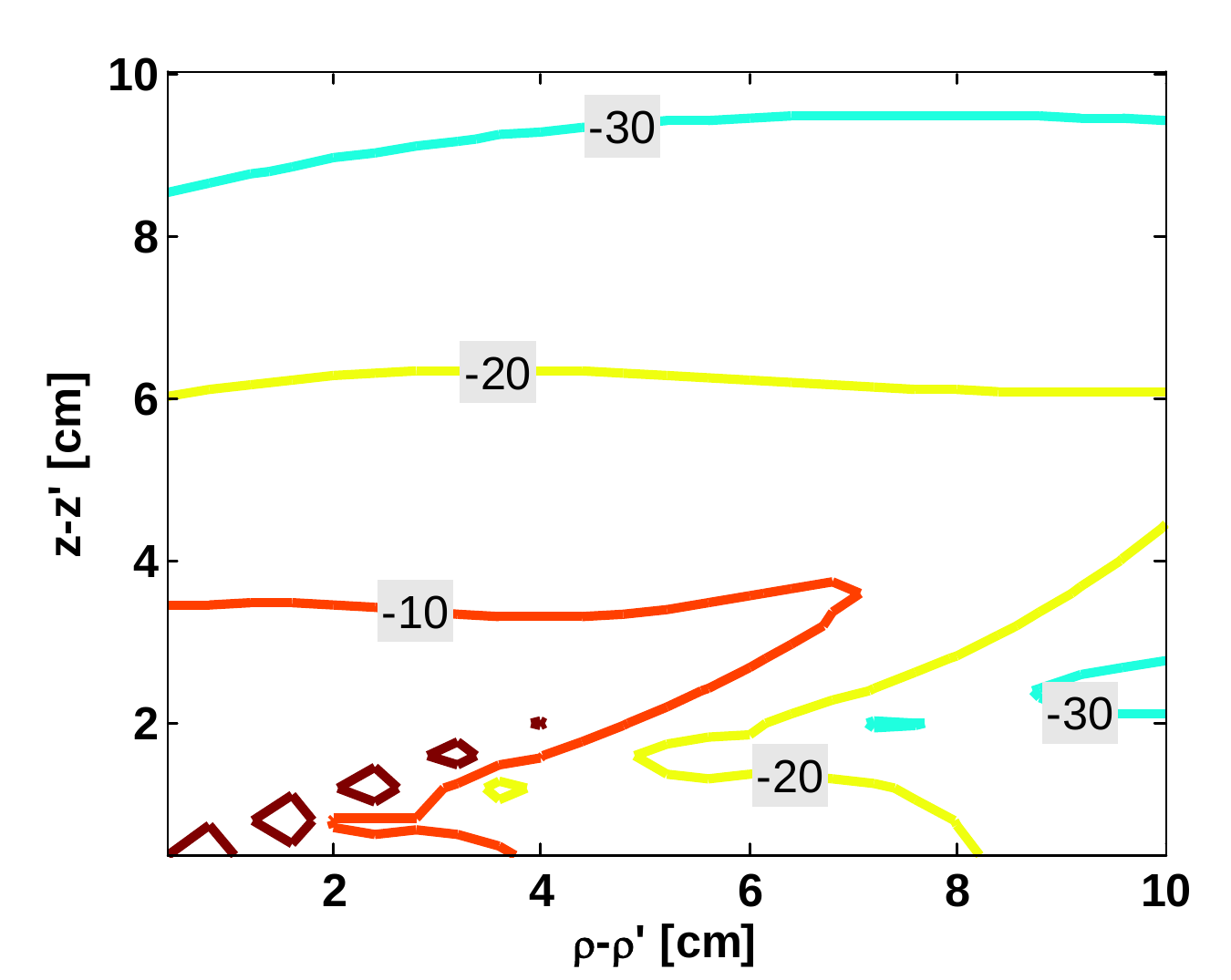}
    }
    \hfill
    \subfloat[\label{ch4.F.sqrt2k.20.1000}]{%
      \includegraphics[width=3.0in]{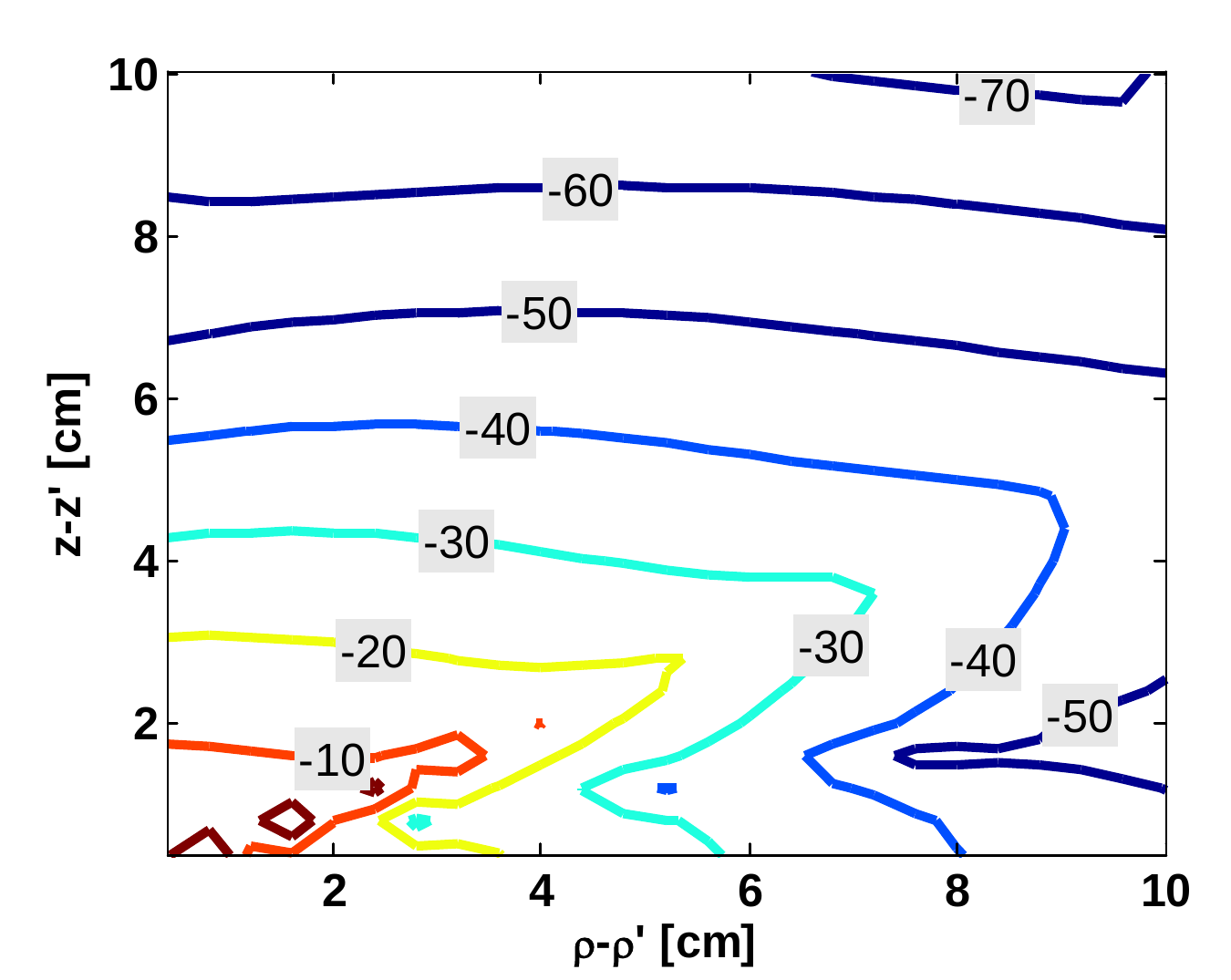}
    }\\
    \subfloat[\label{ch4.F.sqrt2k.10.2000}]{%
      \includegraphics[width=3.0in]{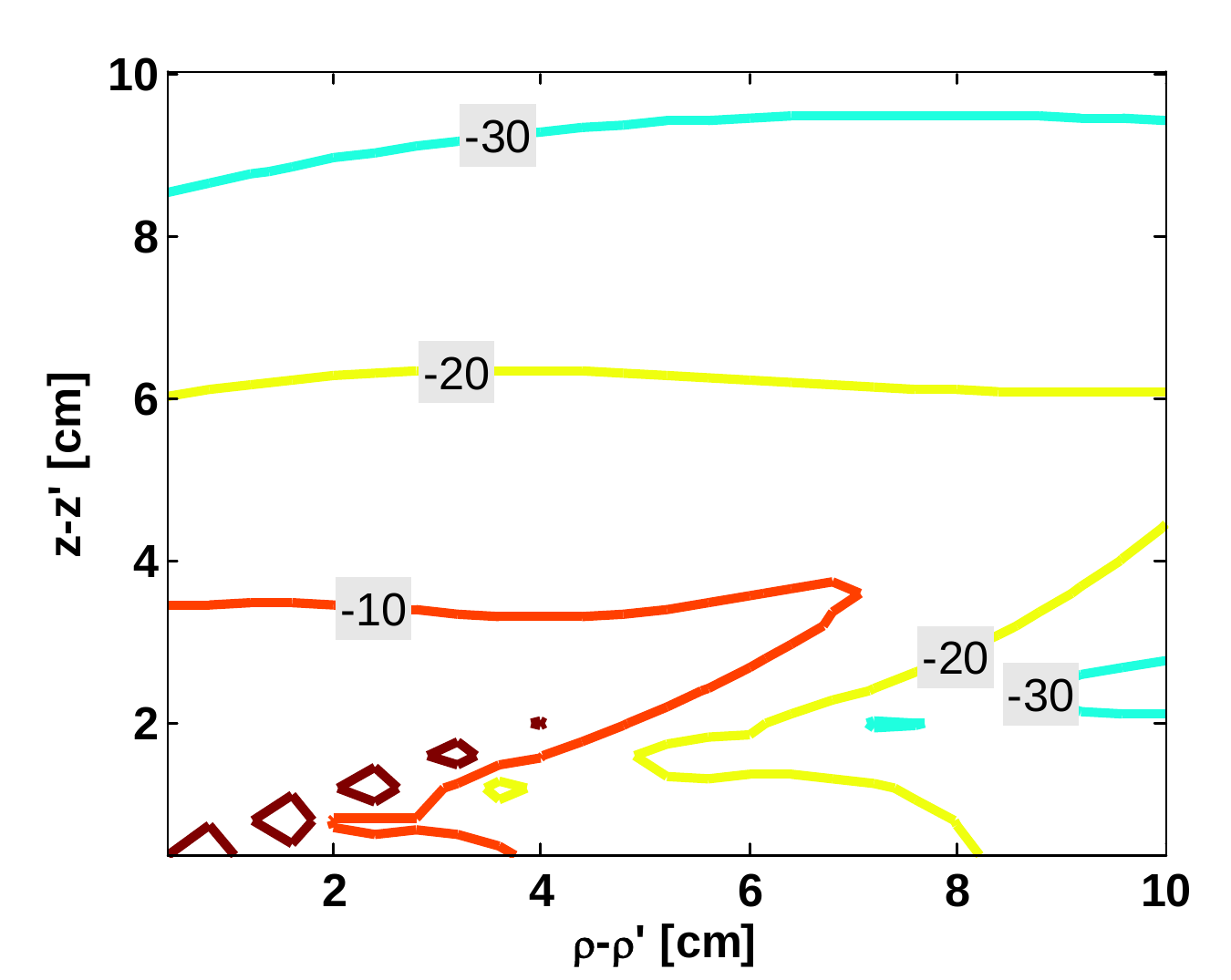}
    }
    \hfill
    \subfloat[\label{ch4.F.sqrt2k.20.2000}]{%
      \includegraphics[width=3.0in]{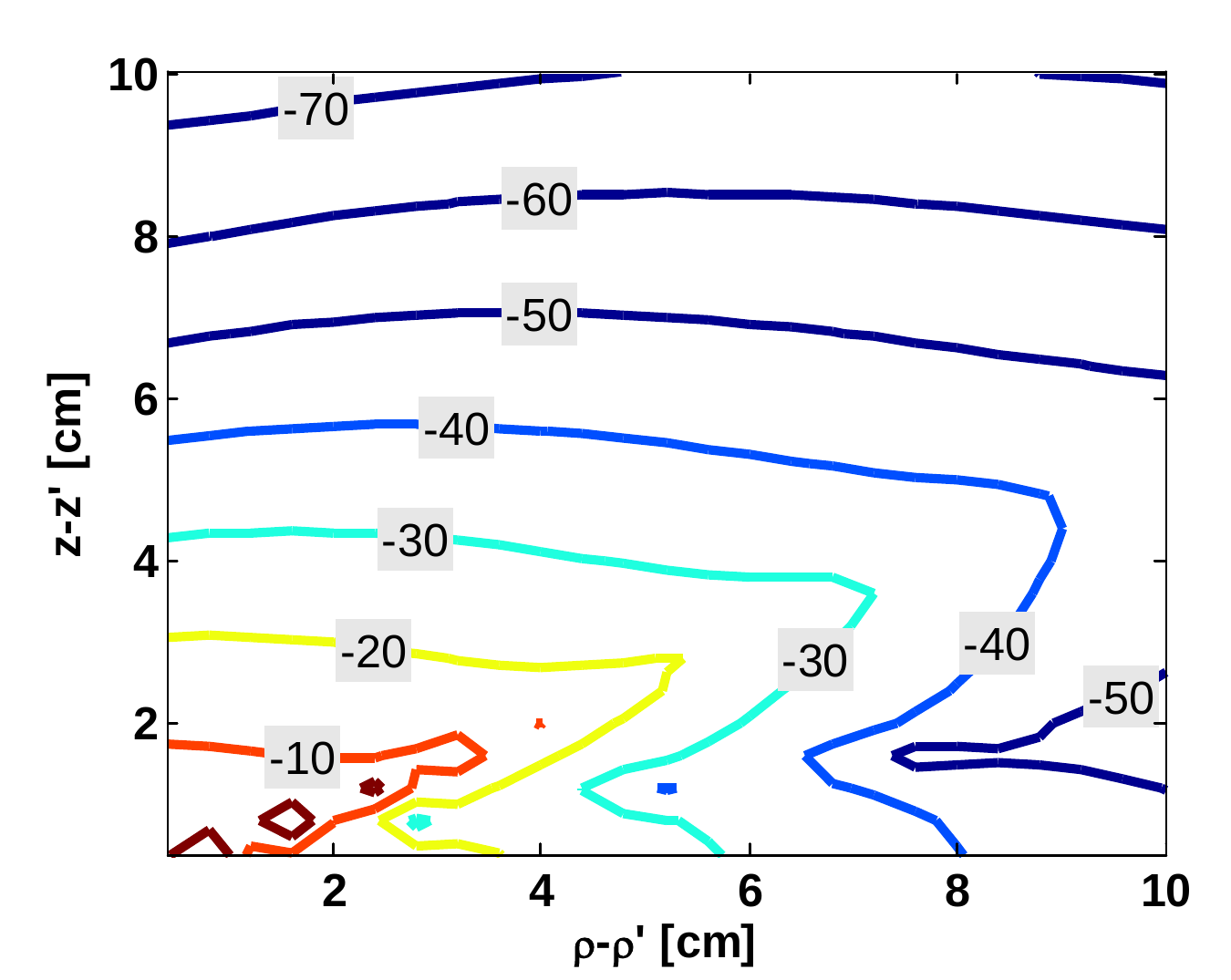}
	}    
    \caption{Relative error distribution with $\kappa=\sqrt{2}$: (a) $n_{max}=10$, $n_{int}=1000$, (b) $n_{max}=20$, $n_{int}=1000$, (c) $n_{max}=10$, $n_{int}=2000$, and (d) $n_{max}=20$, $n_{int}=2000$.}
    \label{ch4.F.error.sqrt2k}
\end{figure}
Fig. \ref{ch4.F.sqrt2k.10.1000} -- \ref{ch4.F.sqrt2k.20.2000} show the relative error distribution for $\kappa_\epsilon=\kappa_\mu=\kappa=\sqrt{2}$, under the assumption of $\epsilon_{p,v}=8\epsilon_0$ [F/m], $\mu_v=8\mu_0$ [H/m], and $\sigma_v=8$ [S/m]. Comparing the cases with $\kappa=4$, $\kappa=2$, $\kappa=\sqrt{2}$, we observe that the error distribution shows a faster rate of decay along the vertical spatial direction for larger $\kappa$.
\begin{figure}[t]
	\centering
	\subfloat[\label{ch4.F.1k.10.1000}]{%
      \includegraphics[width=3.0in]{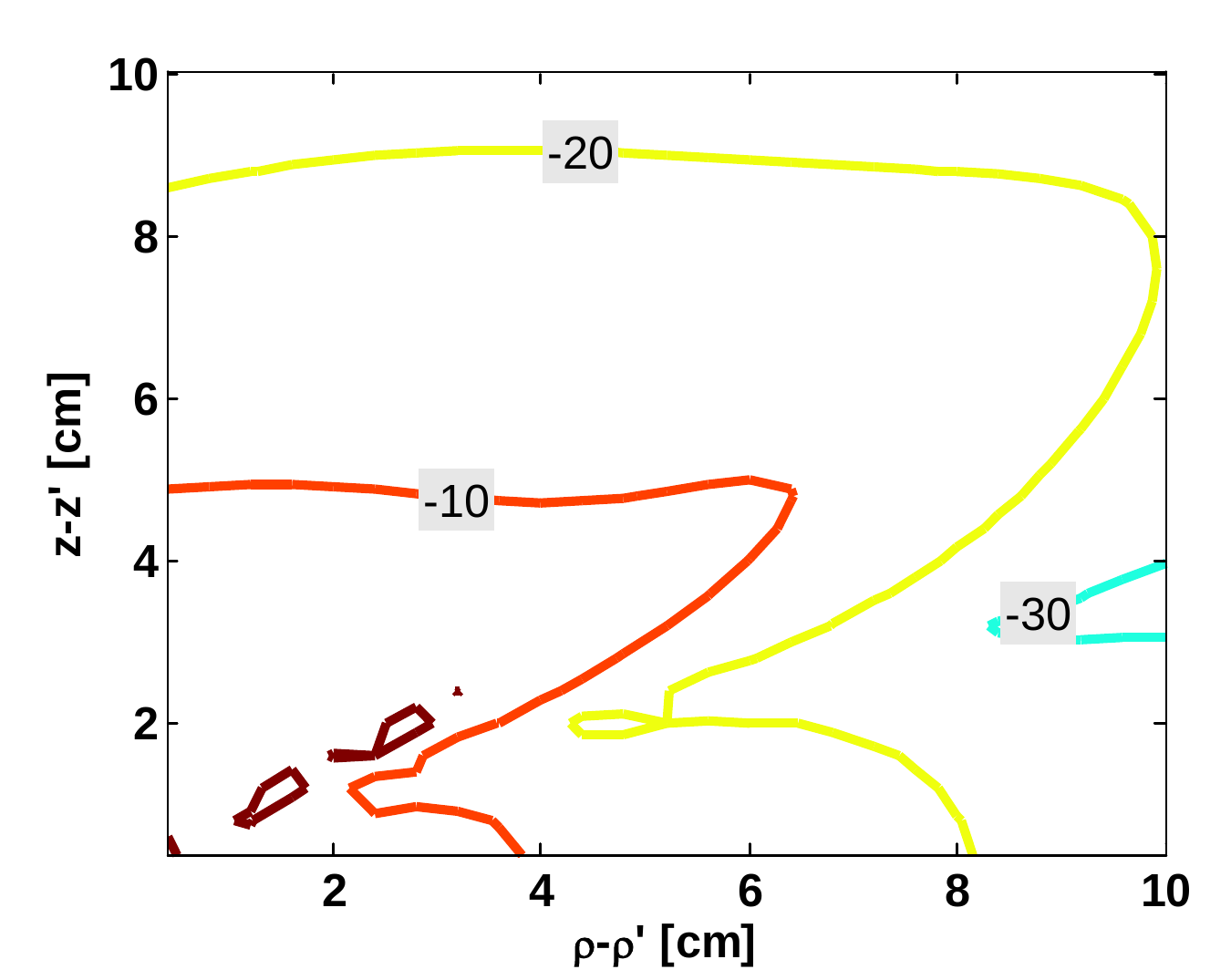}
    }
    \hfill
    \subfloat[\label{ch4.F.1k.20.1000}]{%
      \includegraphics[width=3.0in]{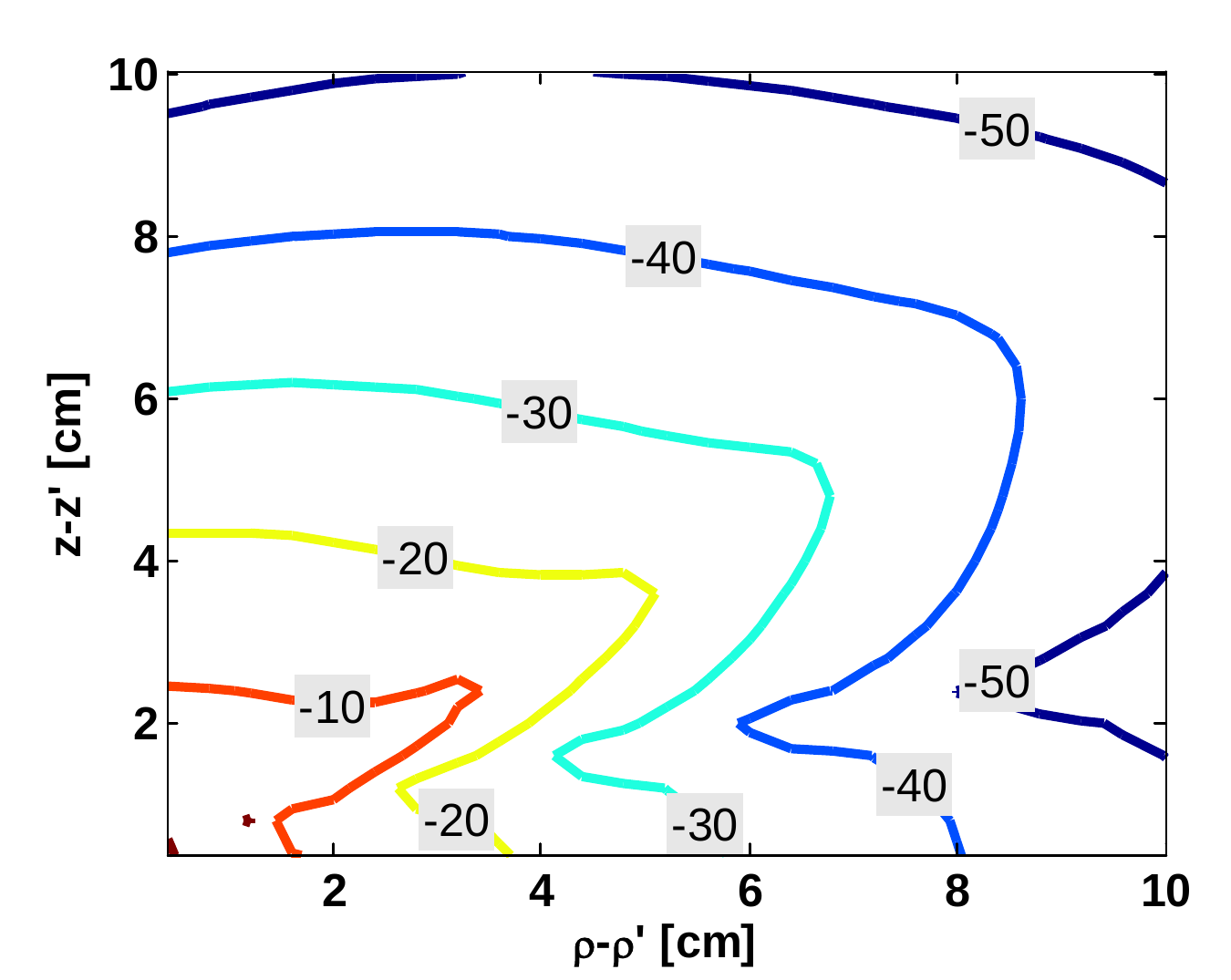}
    }\\
    \subfloat[\label{ch4.F.1k.10.2000}]{%
      \includegraphics[width=3.0in]{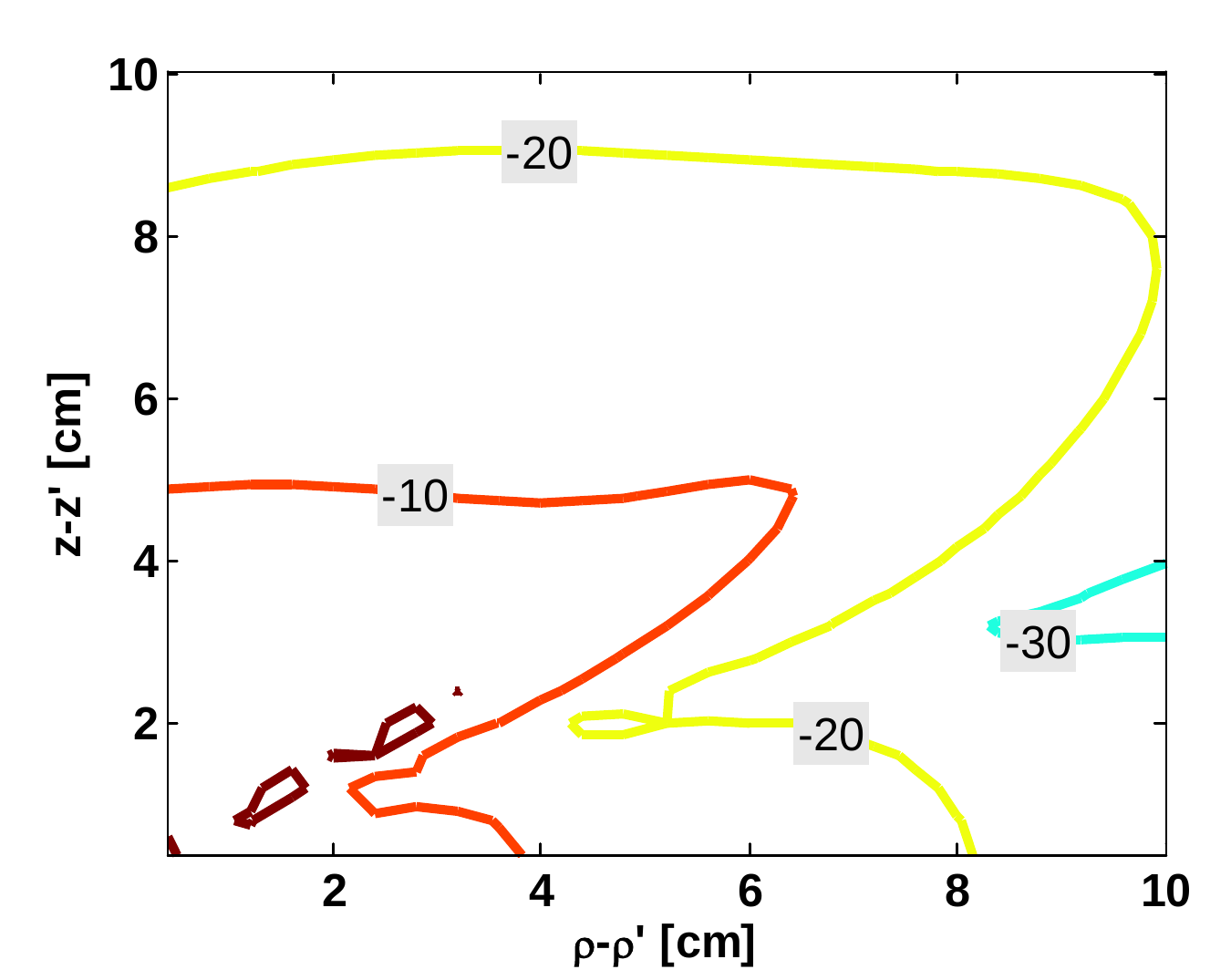}
    }
    \hfill
    \subfloat[\label{ch4.F.1k.20.2000}]{%
      \includegraphics[width=3.0in]{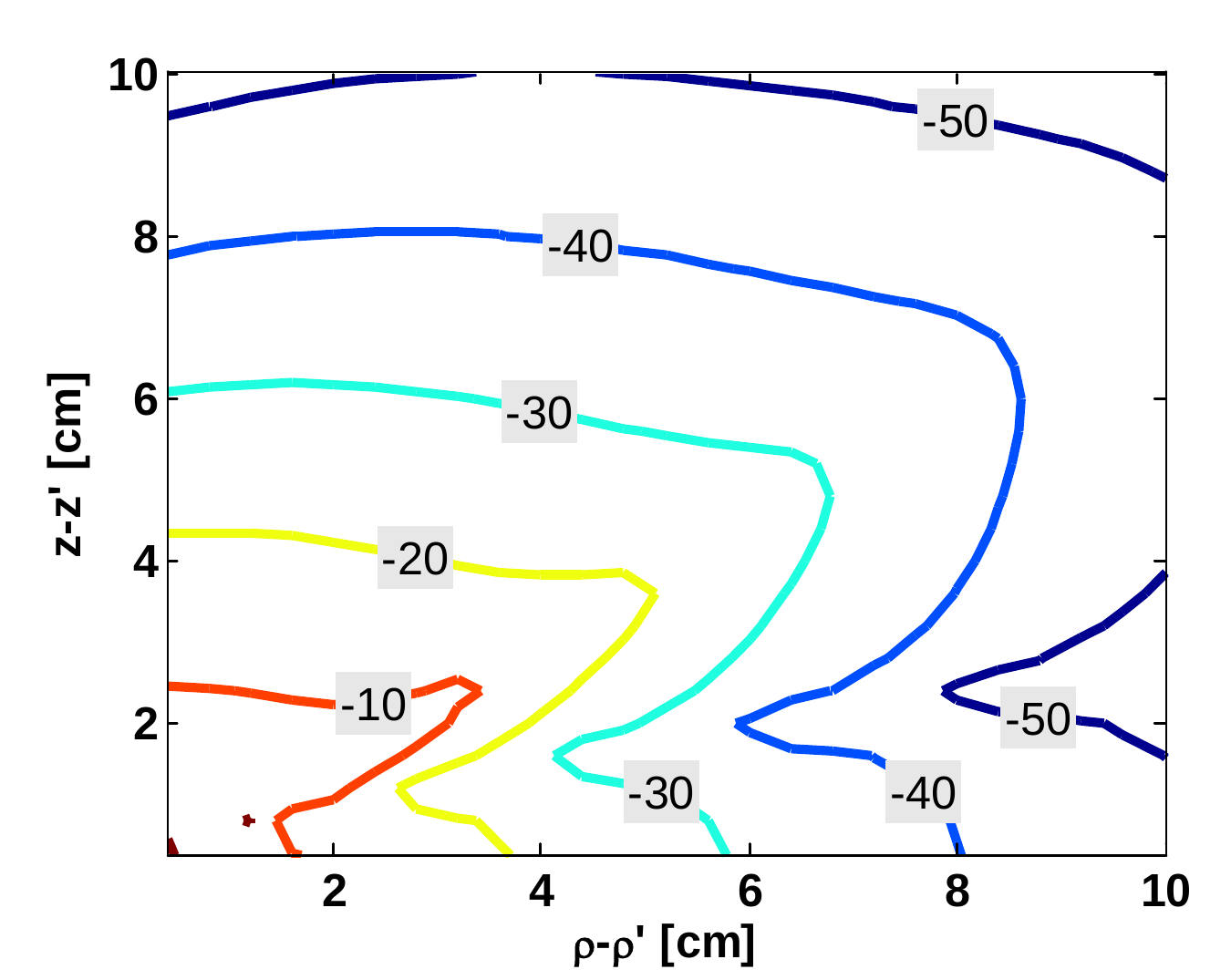}
	}    
    \caption{Relative error distribution with $\kappa=1$: (a) $n_{max}=10$, $n_{int}=1000$, (b) $n_{max}=20$, $n_{int}=1000$, (c) $n_{max}=10$, $n_{int}=2000$, and (d) $n_{max}=20$, $n_{int}=2000$.}
    \label{ch4.F.error.1k}
\end{figure}
Finally, Fig. \ref{ch4.F.1k.10.1000} -- \ref{ch4.F.1k.20.2000} show the relative error distribution for $\kappa_\epsilon=\kappa_\mu=\kappa=1$, under the assumption of $\epsilon_{p,v}=16\epsilon_0$ [F/m], $\mu_v=16\mu_0$ [H/m], and $\sigma_v=16$ [S/m], which recovers the isotropic case.
\begin{figure}[t]
	\centering
	\subfloat[\label{ch4.F.4k.10r10z}]{%
      \includegraphics[width=3.0in]{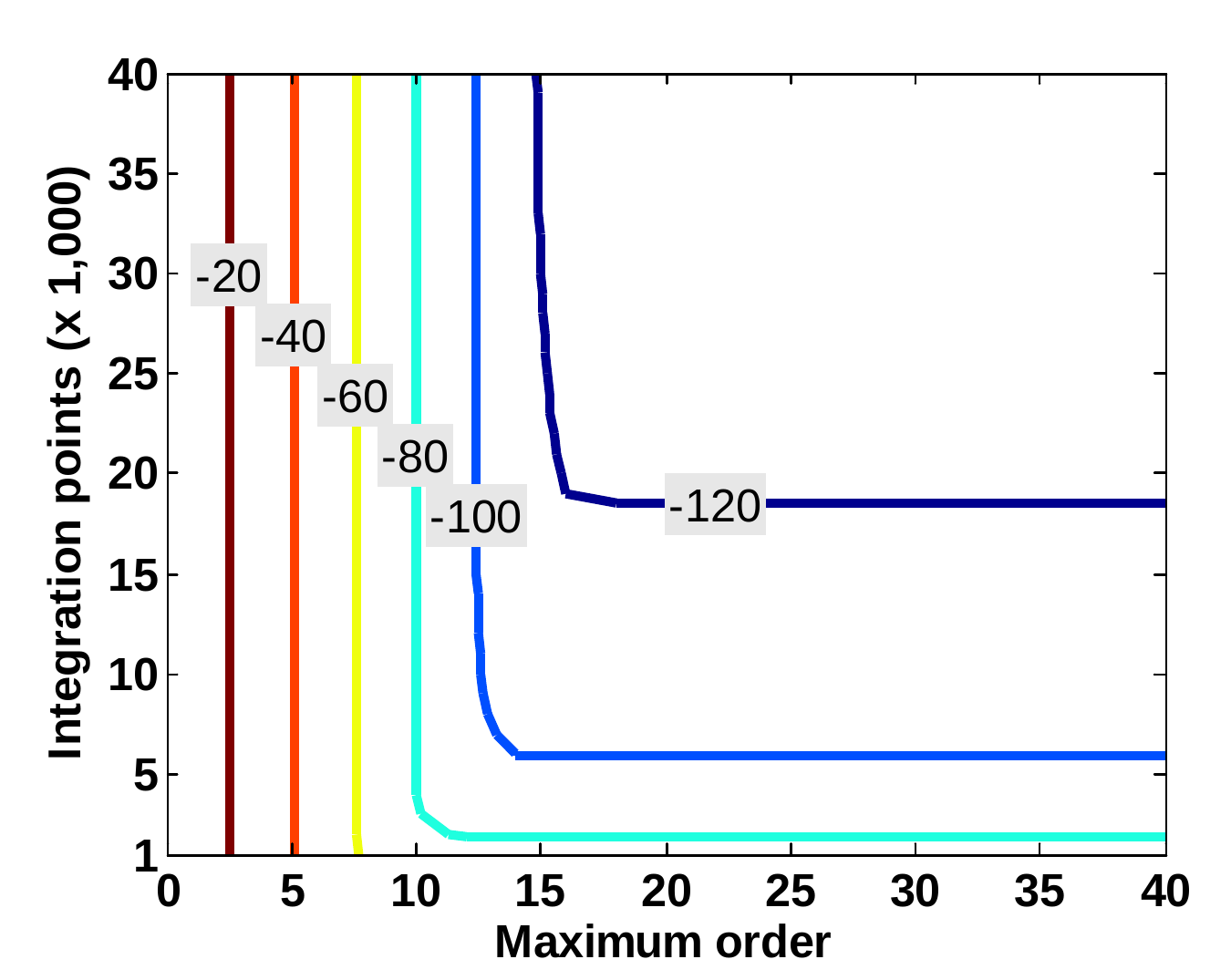}
    }
    \hfill
    \subfloat[\label{ch4.F.2k.10r10z}]{%
      \includegraphics[width=3.0in]{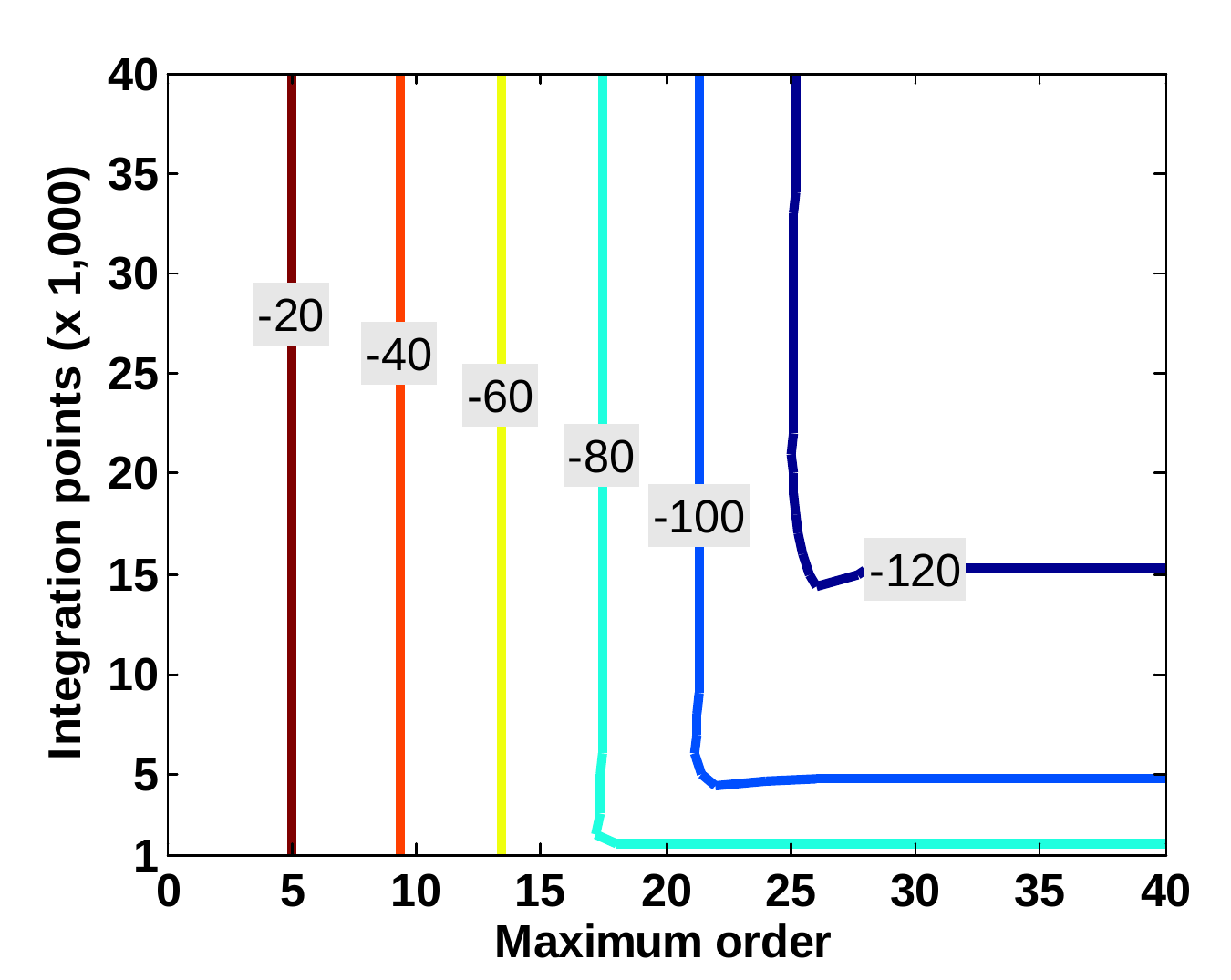}
    }\\
    \subfloat[\label{ch4.F.sqrt2k.10r10z}]{%
      \includegraphics[width=3.0in]{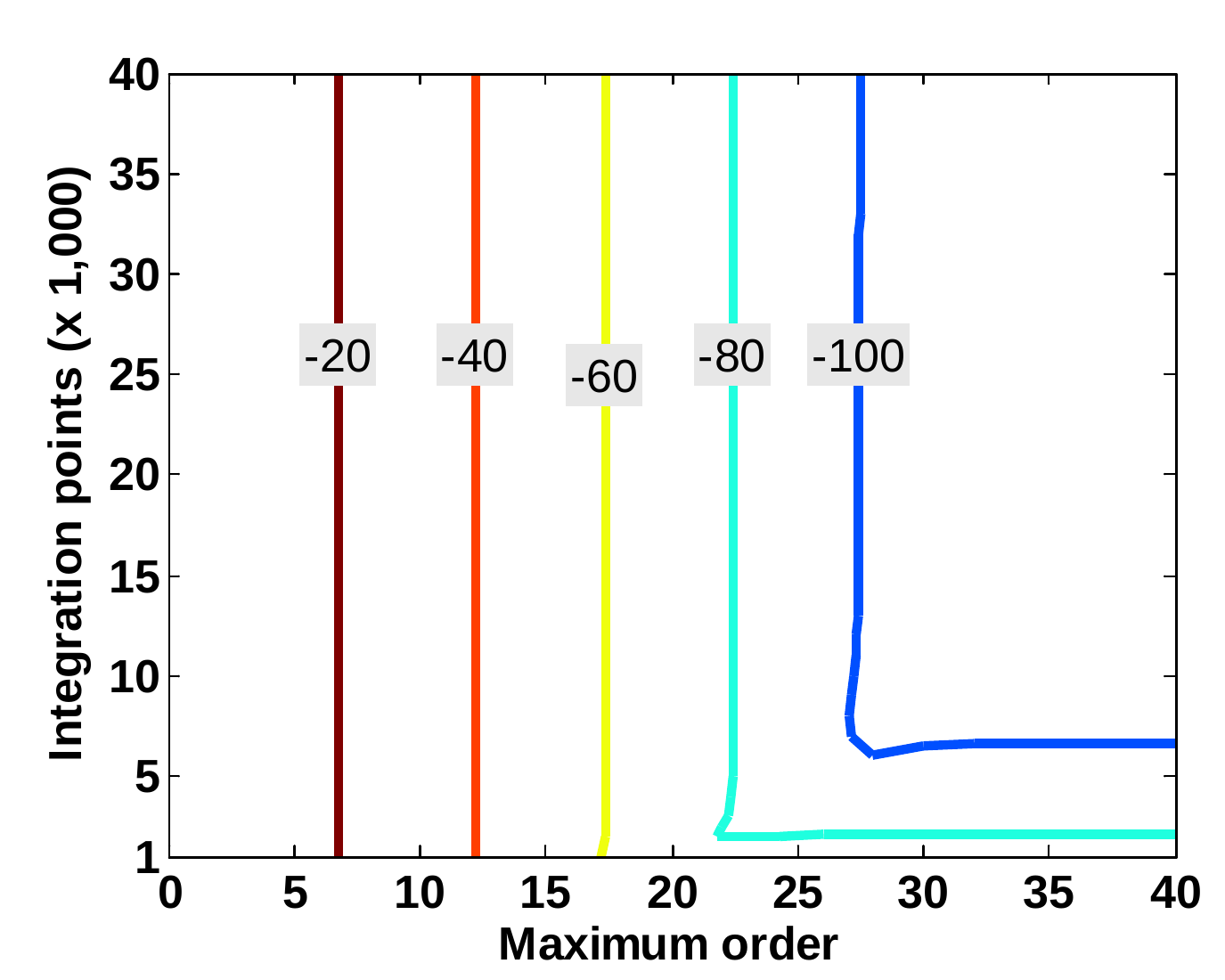}
    }
    \hfill
    \subfloat[\label{ch4.F.1k.10r10z}]{%
      \includegraphics[width=3.0in]{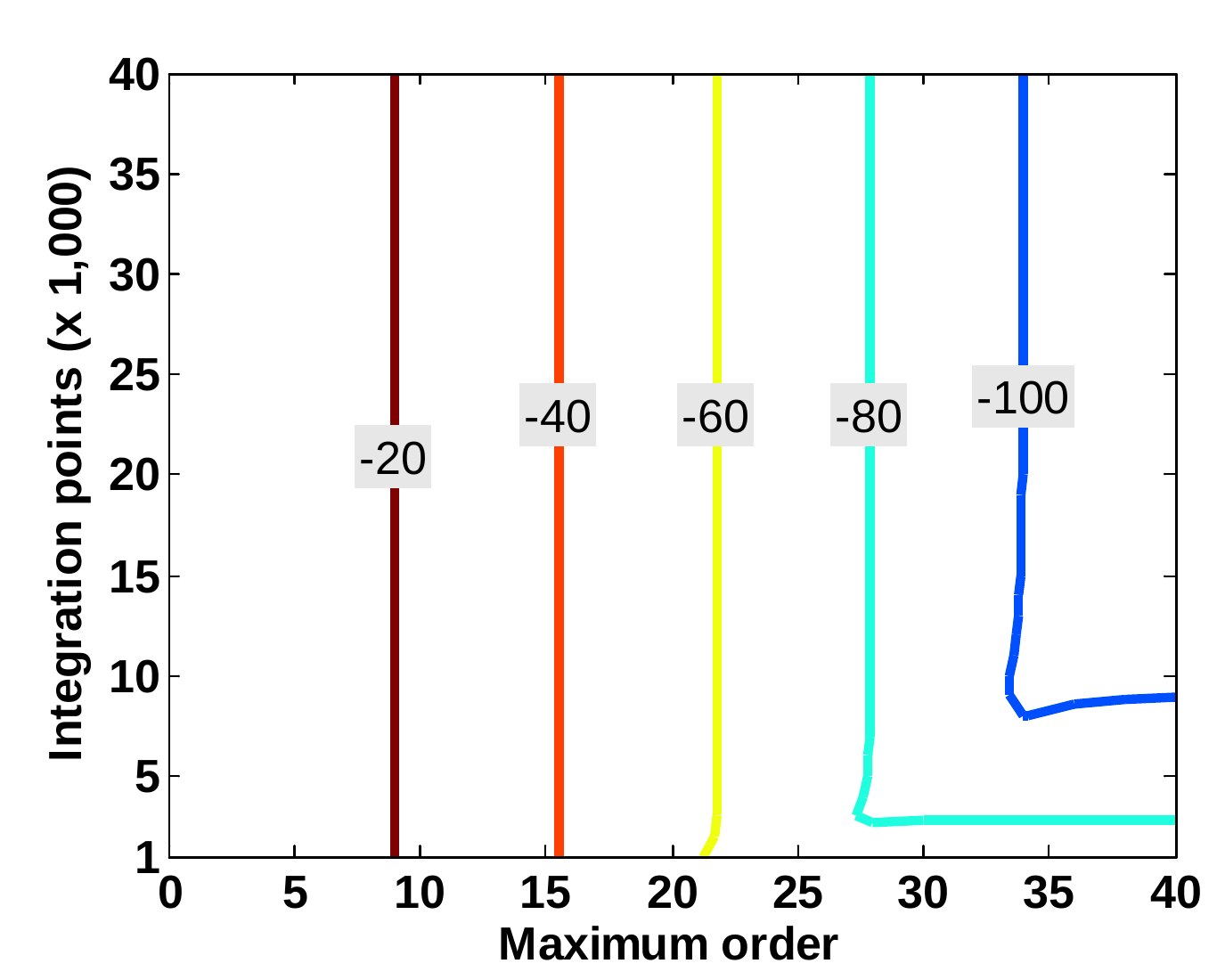}
	}    
    \caption{Relative error distribution in terms of various maximum orders $n_{max}$ and integration points $n_{int}$ with the receiver point at $\rho-\rho'=10$ cm, $\phi-\phi'=0^\circ$, and $z-z'=10$ cm: (a) $\kappa=4$, (b) $\kappa=2$, (c) $\kappa=\sqrt{2}$, and (d) $\kappa=1$.}
    \label{ch4.F.error.10r10z}
\end{figure}
In order to further scrutinize the effect of $n_{max}$ and $n_{int}$, the receiver point is next fixed at $\rho-\rho'=10$ cm, $\phi-\phi'=0^\circ$, and $z-z'=10$ cm. Figure \ref{ch4.F.4k.10r10z} -- \ref{ch4.F.1k.10r10z} show the error as $n_{max}$ and $n_{int}$ vary. 
\textcolor{\Cred}{
When $n_{max}$ is less than about 15, the relative error is not reduced despite the increase in the number of quadrature points. Therefore, $n_{max}$ should be set sufficiently large, otherwise a convergence of results with respect to the number of quadrature points would be only relative and the final results would be still inaccurate. For the types of scenarios considered here, we have observed that $n_{max} \gtrsim 30$ to provide absolute convergence.
Furthermore, it is observed that convergence is achieved faster for larger anisotropy ratios. This stems from an effective spatial ``stretching" along the $\rho$-direction effected when the anisotropy ratio increases. As can be seen in \eqref{ch1.1.E.dispersion.relation.eps} and \eqref{ch1.1.E.dispersion.relation2.mu}, $\tk_\rho^2$ and $\dk_\rho^2$ decrease with an increase in the anisotropy ratio and consequently, higher order modes (with larger $n$) exhibit a faster decay away at the receiver point.  
}

\subsection{Cylindrically layered scenarios}
\label{sec.4.2}
In this section, a number of practical cases of interest are considered to illustrate the applicability of the algorithm. In all the cases, both the relative permittivity $\epsilon_r$ and relative permeability \textcolor{\Cred}{$\mu_r$} are set to one so that $\epsilon_{p,h}=\epsilon_{p,v}=1$
and $\mu_{h}=\mu_{v}=1$, whereas the conductivity tensor (and complex permittivity tensor $\ue$, see \eqref{ch1.1.E.uniaxial.epsilon}) exhibits uniaxial anisotropy where the horizontal resistivity (reciprocal of horizontal conductivity) is set to 5 $\Omega\cdot m$ and the vertical resistivity is changed, leading to different anisotropy ratios $\kappa_\epsilon$.

\begin{figure}[!htbp]
	\centering
	\subfloat[\label{ch4.F.case1}]{%
      \includegraphics[height=3.0in]{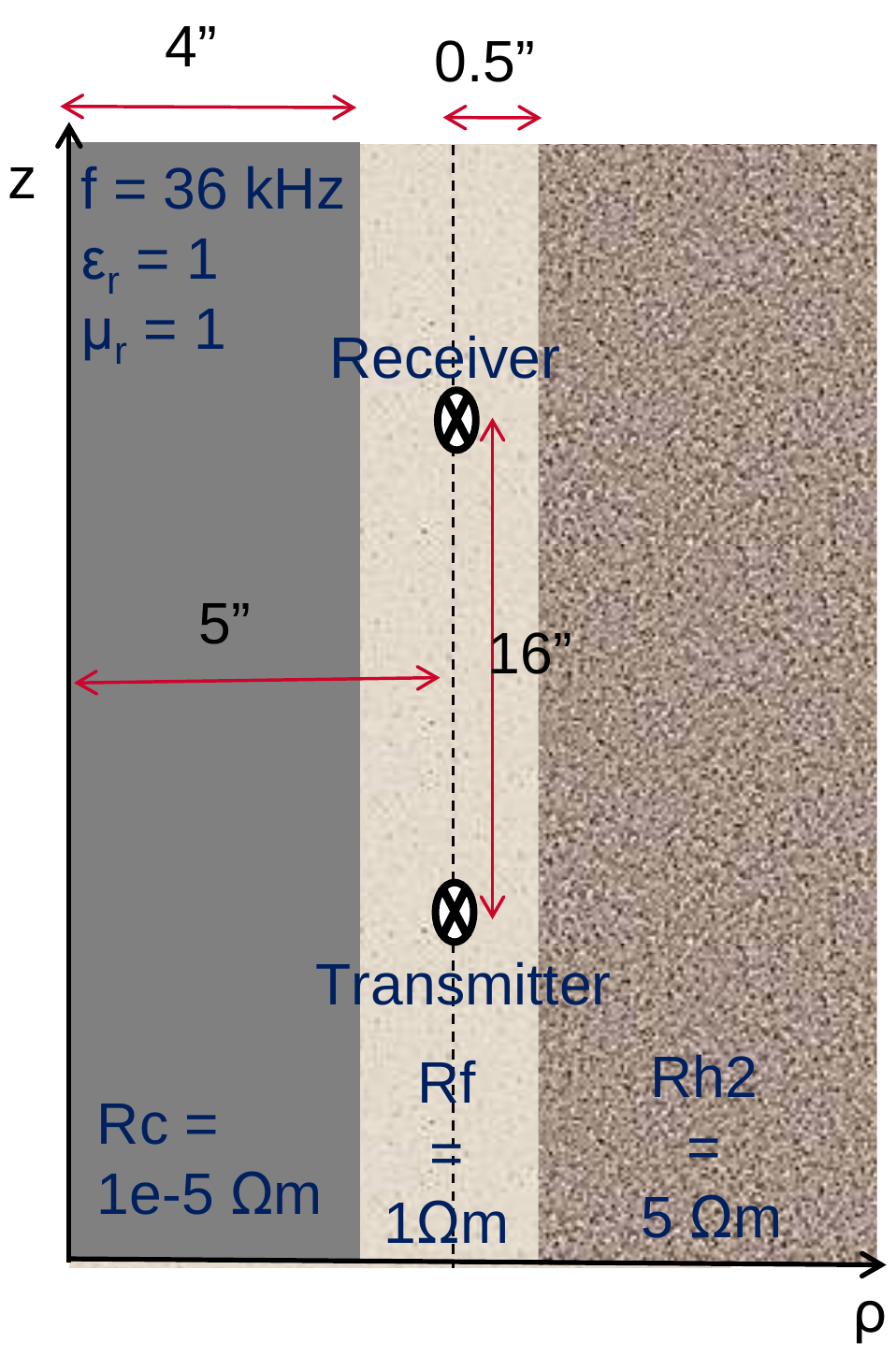}
    }
    \hspace{2 cm}
    \subfloat[\label{ch4.F.case2}]{%
      \includegraphics[height=3.0in]{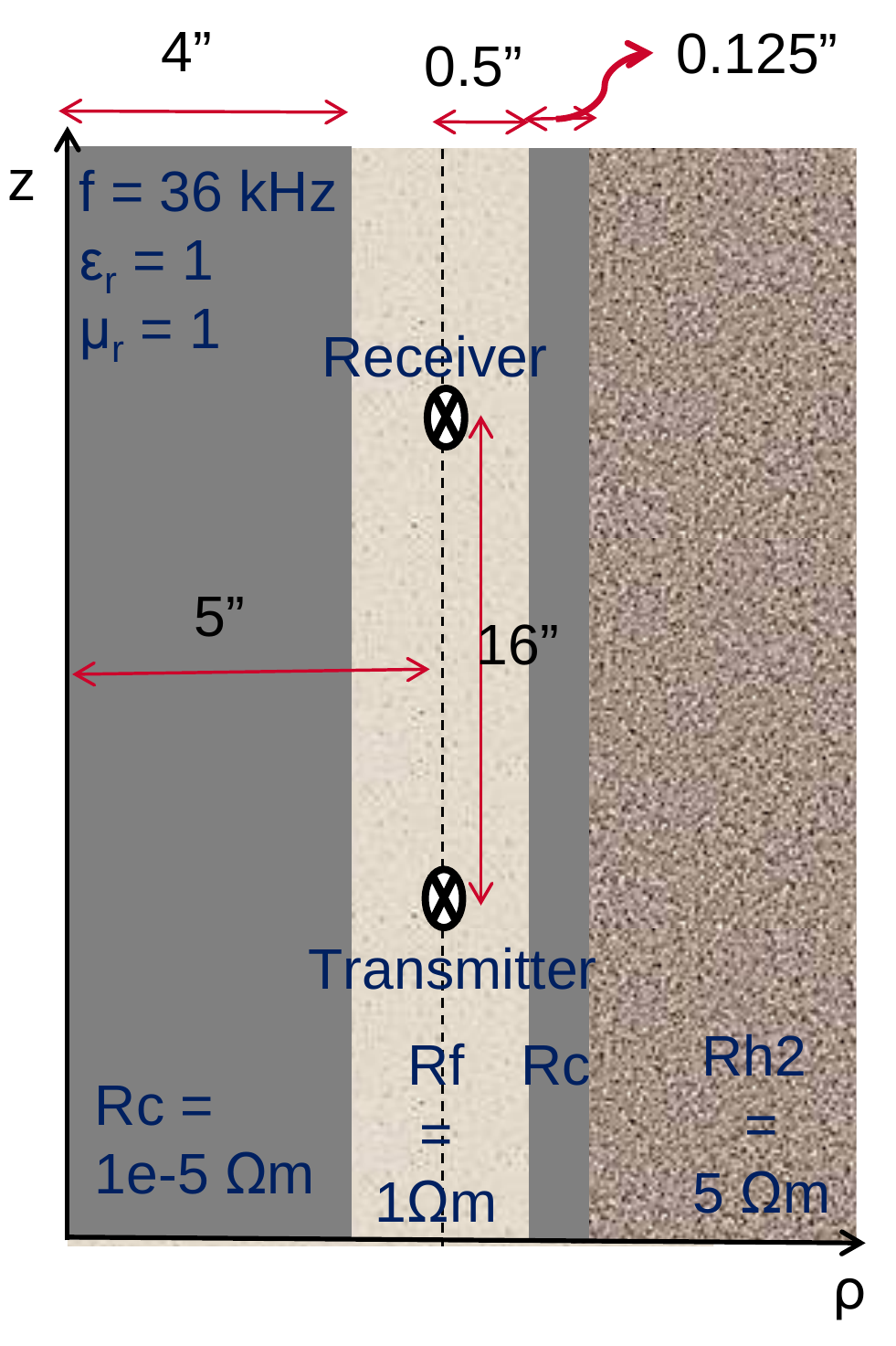}
    }
    \caption{(a) Case 1 in the $\rho z$-plane and (b) Case 2 in the $\rho z$-plane.}
    \label{ch4.F.case12}
\end{figure}

Case 1 is depicted in Fig. \ref{ch4.F.case1}. There are three layers, with the first layer representing a metallic mandrel with high conductivity, the mid-layer representing a borehole filled with an isotropic fluid, and the outermost layer representing the surrounding Earth formation with uniaxial anisotropy. Case 2 is depicted in Fig. \ref{ch4.F.case2}, where a metallic casing (third layer) is inserted between the borehole and anisotropic formation. Table \ref{ch4.T.case1} provides the comparison of corresponding results for Case 1 in terms of the square of anisotropy ratios.
The discrepancy in the magnitude of magnetic fields can be traced to FEM mesh truncation effects: the fields obtained by FEM have smaller magnitudes because the Dirichlet boundary condition at the mesh boundary moves the ground potential (originally at infinity) closer to the source location. This causes a small offset in the results. This is confirmed by Table \ref{ch4.T.case1.dif}, which shows the relative difference in the computed field magnitudes with excellent agreement. Table \ref{ch4.T.case2} and \ref{ch4.T.case2.dif} provide corresponding results for Case 2.

\begin{table}[h]
\begin{center}
\renewcommand{\arraystretch}{1.2}
\setlength{\tabcolsep}{11pt}
\caption{Comparison of magnetic fields in terms of various anisotropy ratios for Case 1.}
    \begin{tabular}{cccc}
        \hline
		Square of & Magnetic field [A/m] & Magnetic field [A/m] & Computing time \\
		anisotropy ratio  $\kappa_\epsilon^2$ & (FEM) & (Present algorithm) & (Present algorithm) \\
        \hline
         1 & 10.3116 $\angle$98.1899$^\circ$ & 10.5475 $\angle$98.1390$^\circ$ & 10 sec. \\
         2 & 10.2365 $\angle$97.5910$^\circ$ & 10.4723 $\angle$97.5486$^\circ$ & 32 sec. \\
         5 & 10.1565 $\angle$96.9883$^\circ$ & 10.3924 $\angle$96.9612$^\circ$ & 32 sec. \\
        10 & 10.1070 $\angle$96.6327$^\circ$ & 10.3428 $\angle$96.6152$^\circ$ & 31 sec. \\	
	    \hline	
    \end{tabular}
    \label{ch4.T.case1}
\vspace{1em}
\caption{Comparison of magnitude difference in magnetic fields for Case 1.}
    \begin{tabular}{ccc}
        \hline      				
		 & FEM & Present algorithm \\
        \hline
        between $\kappa_\epsilon^2=1$ and $\kappa_\epsilon^2=2$ & 0.0751  & 0.0752 \\
        between $\kappa_\epsilon^2=2$ and $\kappa_\epsilon^2=5$ & 0.0800  & 0.0799 \\
        between $\kappa_\epsilon^2=5$ and $\kappa_\epsilon^2=10$ & 0.0495  & 0.0496 \\        	
	\hline	
    \end{tabular}
    \label{ch4.T.case1.dif}
\vspace{3em}
\caption{Comparison of magnetic fields in terms of various anisotropy ratios for Case 2.}
    \begin{tabular}{cccc}
        \hline
        Square of & Magnetic field [A/m] & Magnetic field [A/m] & Computing time \\
		anisotropy ratio  $\kappa_\epsilon^2$ & (FEM) & (Present algorithm) & (Present algorithm) \\
        \hline
         1 & 46.6091 $\angle$118.4181$^\circ$ & 46.6303 $\angle$118.4324$^\circ$ & 15 sec. \\
         2 & 46.6099 $\angle$118.4234$^\circ$ & 46.6311 $\angle$118.4381$^\circ$ & 44 sec. \\
         5 & 46.6110 $\angle$118.4283$^\circ$ & 46.6321 $\angle$118.4432$^\circ$ & 44 sec. \\
	    10 & 46.6118 $\angle$118.4310$^\circ$ & 46.6329 $\angle$118.4459$^\circ$ & 44 sec. \\	
	\hline	
    \end{tabular}
    \label{ch4.T.case2}
\vspace{1em}    
\caption{Comparison of magnitude difference in magnetic fields for Case 2.}
    \begin{tabular}{ccc}
        \hline      				
		 & FEM & Present algorithm \\
        \hline
        between $\kappa_\epsilon^2=1$ and $\kappa_\epsilon^2=2$ & -0.0008  & -0.0008 \\
        between $\kappa_\epsilon^2=2$ and $\kappa_\epsilon^2=5$ & -0.0011  & -0.0010 \\
        between $\kappa_\epsilon^2=5$ and $\kappa_\epsilon^2=10$ & -0.0008  & -0.0008 \\        	
	\hline	
    \end{tabular}
    \label{ch4.T.case2.dif}
\end{center}
\end{table}

\begin{figure}[!htbp]
	\centering
	\subfloat[\label{ch4.F.case3}]{%
      \includegraphics[height=3.0in]{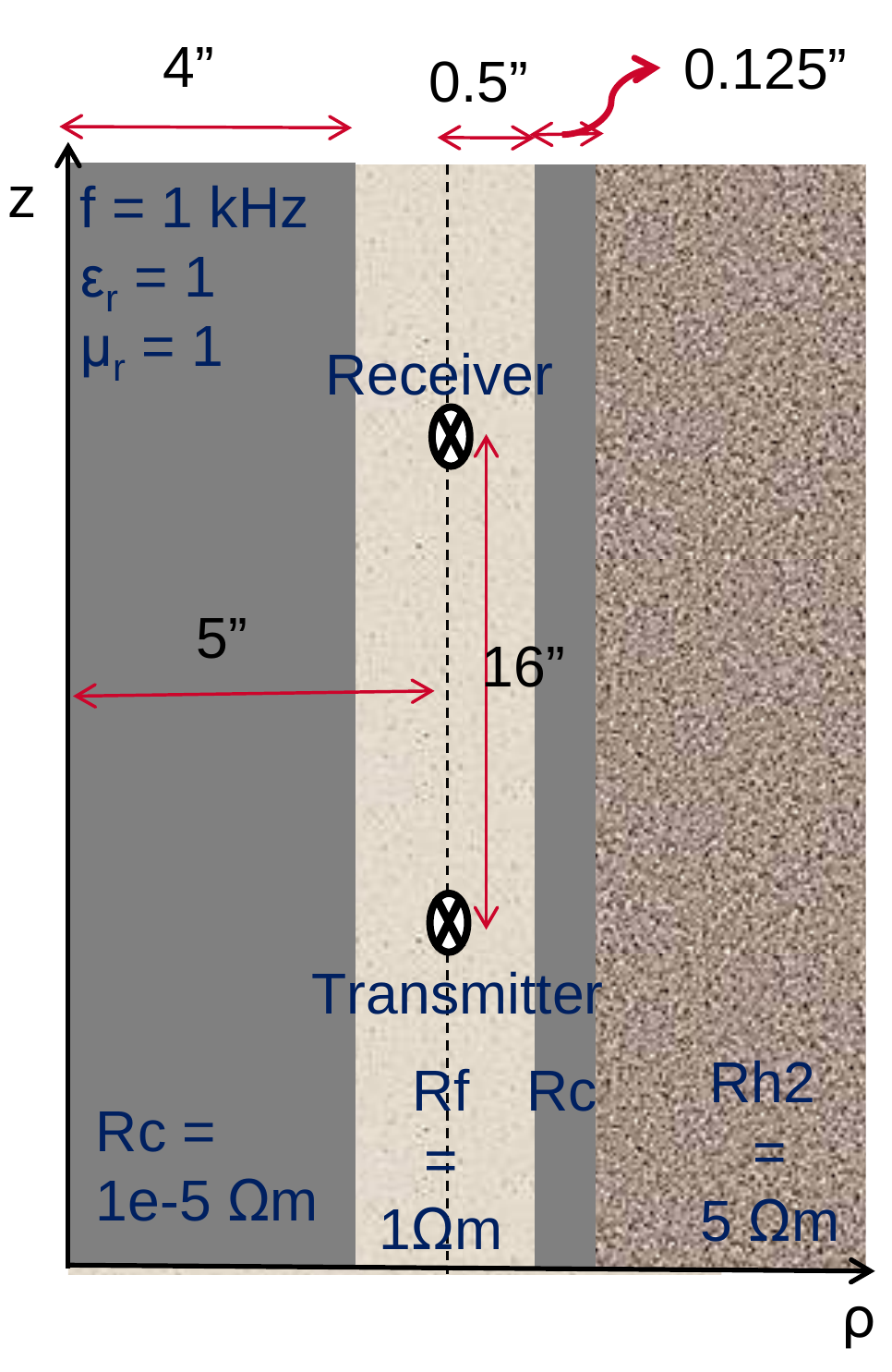}
    }
    \hspace{2 cm}
    \subfloat[\label{ch4.F.case4}]{%
      \includegraphics[height=3.0in]{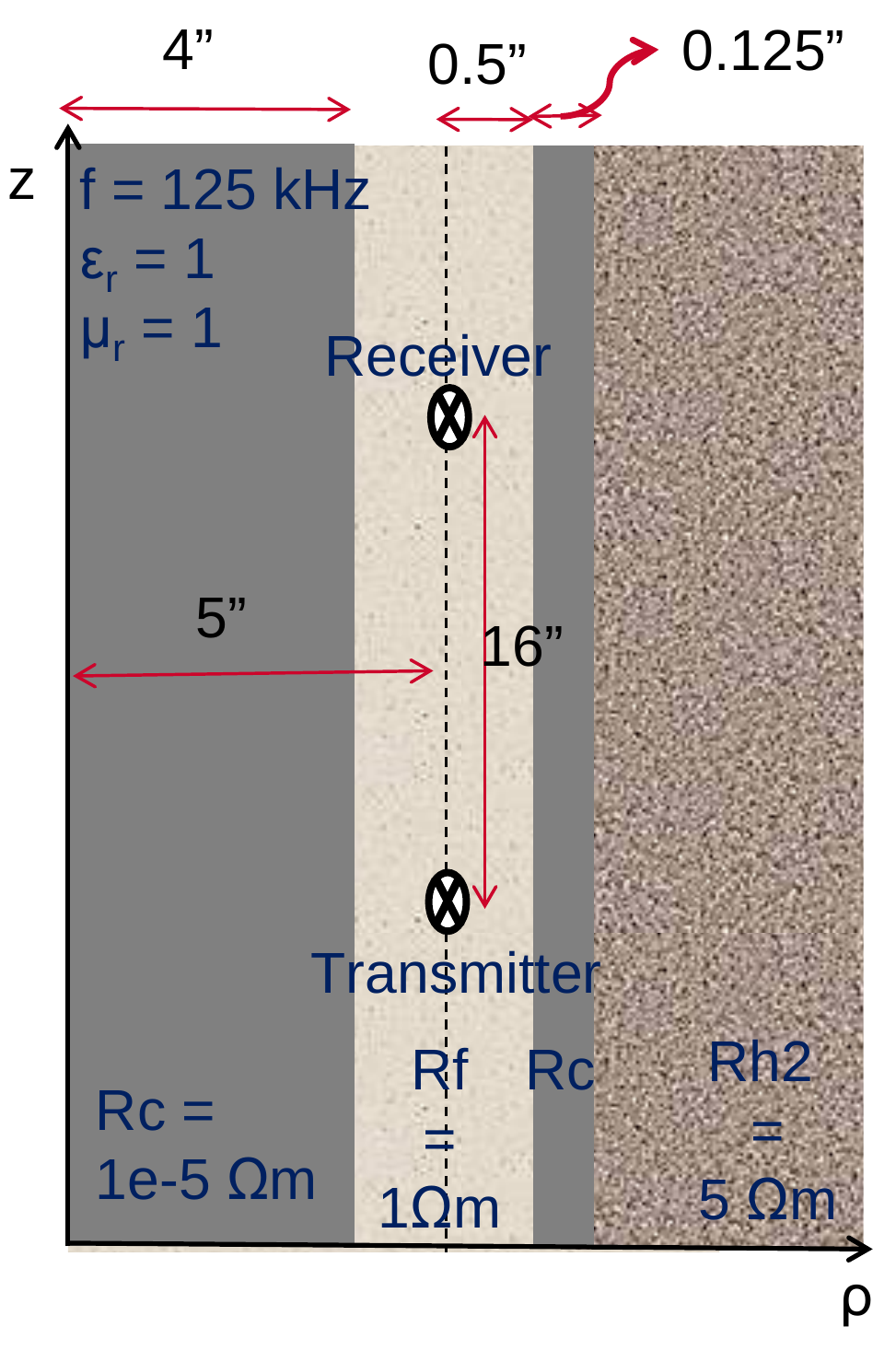}
    }
    \caption{(a) Case 3 in the $\rho z$-plane and (b) Case 4 in the $\rho z$-plane.}
    \label{ch4.F.case34}
\end{figure}

Case 3 and 4 are depicted in Figs. \ref{ch4.F.case3} and \ref{ch4.F.case4}, which are the same as Case 2 except for the operating frequencies, which for Case 3 is 1 kHz and for Case 4 is 125 kHz. Tables \ref{ch4.T.case3}, \ref{ch4.T.case3.dif}, \ref{ch4.T.case4}, and \ref{ch4.T.case4.dif} provide the corresponding results.\\

\begin{table}[h]
\begin{center}
\renewcommand{\arraystretch}{1.2}
\setlength{\tabcolsep}{11pt}
\caption{Comparison of magnetic fields in terms of various anisotropy ratios for Case 3.}
    \begin{tabular}{cccc}
        \hline
        Square of & Magnetic field [A/m] & Magnetic field [A/m] & Computing time \\
		anisotropy ratio  $\kappa_\epsilon^2$ & (FEM) & (Present algorithm) & (Present algorithm) \\
        \hline
         1 & 290.2144 $\angle$127.4332$^\circ$ & 290.2711 $\angle$127.4185$^\circ$ & 7 sec. \\
         2 & 289.9780 $\angle$127.4252$^\circ$ & 290.0564 $\angle$127.4170$^\circ$ & 22 sec. \\
         5 & 289.7678 $\angle$127.4089$^\circ$ & 289.8157 $\angle$127.4175$^\circ$ & 22 sec. \\
        10 & 289.6754 $\angle$127.3988$^\circ$ & 289.6589 $\angle$127.4189$^\circ$ & 22 sec. \\	
	\hline	
    \end{tabular}
    \label{ch4.T.case3}
\vspace{1em} 
\caption{Comparison of magnitude difference in magnetic fields for Case 3.}
    \begin{tabular}{ccc}
        \hline      				
		 & FEM & Present algorithm \\
        \hline
        between $\kappa_\epsilon^2=1$ and $\kappa_\epsilon^2=2$ & 0.2364  & 0.2147 \\
        between $\kappa_\epsilon^2=2$ and $\kappa_\epsilon^2=5$ & 0.2102  & 0.2407 \\
        between $\kappa_\epsilon^2=5$ and $\kappa_\epsilon^2=10$ & 0.0924  & 0.1568 \\        	
	\hline	
    \end{tabular}
    \label{ch4.T.case3.dif}
\vspace{3em}
\caption{Comparison of magnetic fields in terms of various anisotropy ratios for Case 4.}
    \begin{tabular}{cccc}
        \hline
        Square of & Magnetic field [A/m] & Magnetic field [A/m] & Computing time \\
		anisotropy ratio  $\kappa_\epsilon^2$ & (FEM) & (Present algorithm) & (Present algorithm) \\
        \hline
         1 & 18.7957 $\angle$110.9753$^\circ$ & 18.8074 $\angle$110.9191$^\circ$ & 7 sec. \\
         2 & 18.7959 $\angle$110.9762$^\circ$ & 18.8076 $\angle$110.9200$^\circ$ & 19 sec. \\
         5 & 18.7962 $\angle$110.9770$^\circ$ & 18.8079 $\angle$110.9209$^\circ$ & 19 sec. \\
        10 & 18.7963 $\angle$110.9774$^\circ$ & 18.8080 $\angle$110.9213$^\circ$ & 19 sec. \\	
	\hline	
    \end{tabular}
    \label{ch4.T.case4}
\vspace{1em} 
\caption{Comparison of magnitude difference in magnetic fields for Case 4.}
    \begin{tabular}{ccc}
        \hline      				
		 & FEM & Present algorithm \\
        \hline
        between $\kappa_\epsilon^2=1$ and $\kappa_\epsilon^2=2$ & -0.0002  & -0.0002 \\
        between $\kappa_\epsilon^2=2$ and $\kappa_\epsilon^2=5$ & -0.0003  & -0.0003 \\
        between $\kappa_\epsilon^2=5$ and $\kappa_\epsilon^2=10$ & -0.0001  & -0.0001 \\        	
	\hline	
    \end{tabular}
    \label{ch4.T.case4.dif}    
\end{center}
\end{table}

\begin{figure}[!htbp]
	\centering
	\subfloat[\label{ch4.F.case6}]{%
      \includegraphics[height=3.0in]{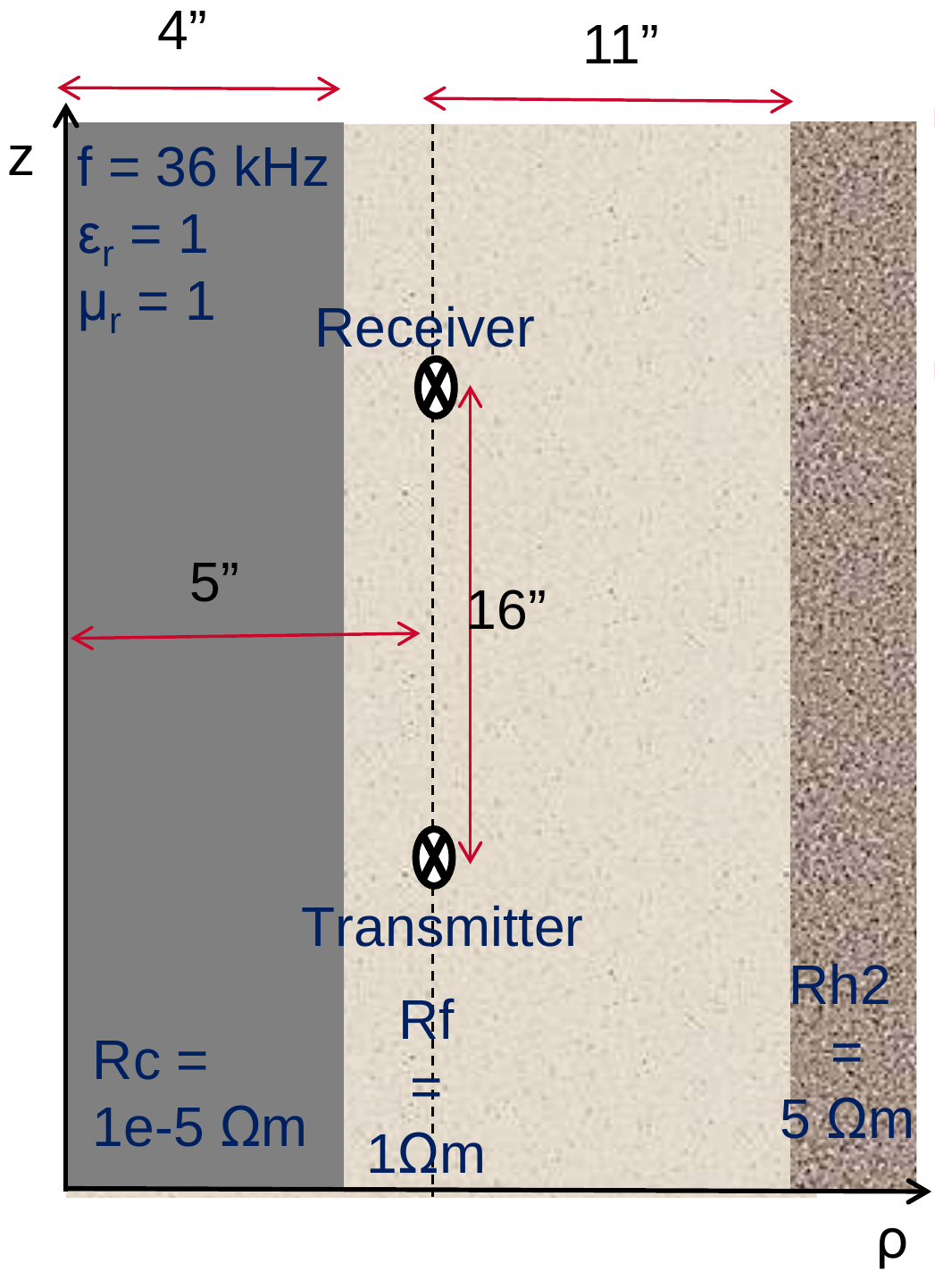}
    }
    \hspace{2 cm}
    \subfloat[\label{ch4.F.case9}]{%
      \includegraphics[height=3.0in]{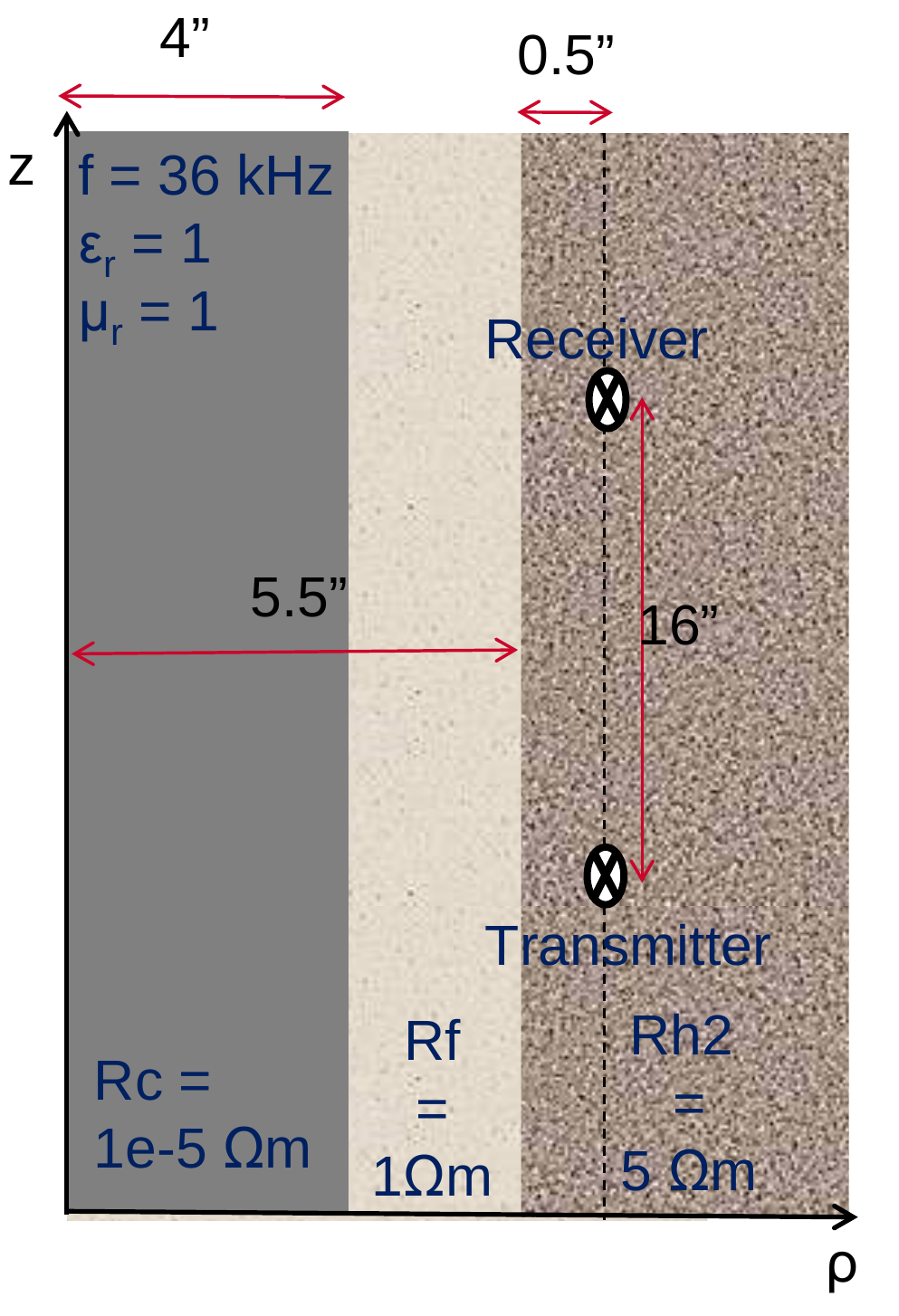}
    }
    \caption{(a) Case 5 in the $\rho z$-plane and (b) Case 6 in the $\rho z$-plane.}
    \label{ch4.F.case56}
\end{figure}

Case 5 and 6 are depicted in Figs. \ref{ch4.F.case6} and \ref{ch4.F.case9}. For Case 5, the borehole is extended to $16^{\prime\prime}$ without casing. For Case 6, both the transmitter and receiver are positioned inside the formation, which again has uniaxial anisotropy.  Tables \ref{ch4.T.case6} and \ref{ch4.T.case9} provide the comparison of corresponding results for Case 5 and Case 6 in terms of the anisotropy ratios squared. Tables \ref{ch4.T.case6.dif} and \ref{ch4.T.case9.dif} show the relative difference in the field magnitude for each case.\\

\begin{table}[h]
\begin{center}
\renewcommand{\arraystretch}{1.2}
\setlength{\tabcolsep}{11pt}
\caption{Comparison of magnetic fields in terms of various anisotropy ratios for Case 5.}
    \begin{tabular}{cccc}
        \hline
        Square of & Magnetic field [A/m] & Magnetic field [A/m] & Computing time \\
		anisotropy ratio  $\kappa_\epsilon^2$ & (FEM) & (Present algorithm) & (Present algorithm) \\
        \hline
         1 & 10.7855 $\angle$100.0572$^\circ$ & 10.7857 $\angle$100.0586$^\circ$ & 10 sec. \\
         2 & 10.6723 $\angle$99.4009$^\circ$ & 10.6721 $\angle$99.4024$^\circ$ & 29 sec. \\
         5 & 10.5553 $\angle$98.7180$^\circ$ & 10.5546 $\angle$98.7190$^\circ$ & 29 sec. \\
        10 & 10.4847 $\angle$98.3108$^\circ$ & 10.4839 $\angle$98.3113$^\circ$ & 29 sec. \\	
	\hline	
    \end{tabular}
    \label{ch4.T.case6}
\vspace{1em} 
\caption{Comparison of magnitude difference in magnetic fields for Case 5.}
    \begin{tabular}{ccc}
        \hline      				
		 & FEM & Present algorithm \\
        \hline
        between $\kappa_\epsilon^2=1$ and $\kappa_\epsilon^2=2$ & 0.1132  & 0.1136 \\
        between $\kappa_\epsilon^2=2$ and $\kappa_\epsilon^2=5$ & 0.1170  & 0.1175 \\
        between $\kappa_\epsilon^2=5$ and $\kappa_\epsilon^2=10$ & 0.0706  & 0.0707 \\        	
	\hline	
    \end{tabular}
    \label{ch4.T.case6.dif}
\vspace{3em} 
\caption{Comparison of magnetic fields in terms of various anisotropy ratios for Case 6.}
    \begin{tabular}{cccc}
        \hline
        Square of & Magnetic field [A/m] & Magnetic field [A/m] & Computing time \\
		anisotropy ratio  $\kappa_\epsilon^2$ & (FEM) & (Present algorithm) & (Present algorithm) \\
        \hline
         1 & 8.1259 $\angle$97.0379$^\circ$ & 8.1326 $\angle$97.0341$^\circ$ & 11 sec. \\
         2 & 8.0817 $\angle$96.5262$^\circ$ & 8.0814 $\angle$96.4841$^\circ$ & 32 sec. \\
         5 & 8.0276 $\angle$95.9964$^\circ$ & 8.0271 $\angle$95.9416$^\circ$ & 32 sec. \\
        10 & 7.9939 $\angle$95.6786$^\circ$ & 7.9933 $\angle$95.6240$^\circ$ & 32 sec. \\	
	\hline	
    \end{tabular}
    \label{ch4.T.case9}
\vspace{1em} 
\caption{Comparison of magnitude difference in magnetic fields for Case 6.}
    \begin{tabular}{ccc}
        \hline      				
		 & FEM & Present algorithm \\
        \hline
        between $\kappa_\epsilon^2=1$ and $\kappa_\epsilon^2=2$ & 0.0442  & 0.0512 \\
        between $\kappa_\epsilon^2=2$ and $\kappa_\epsilon^2=5$ & 0.0541  & 0.0543 \\
        between $\kappa_\epsilon^2=5$ and $\kappa_\epsilon^2=10$ & 0.0337  & 0.0338 \\        	
	\hline	
    \end{tabular}
    \label{ch4.T.case9.dif}    
\end{center}
\end{table}

\begin{figure}[t]
	\centering
	\subfloat[\label{ch4.F.case1.y.k1}]{%
      \includegraphics[width=3.0in]{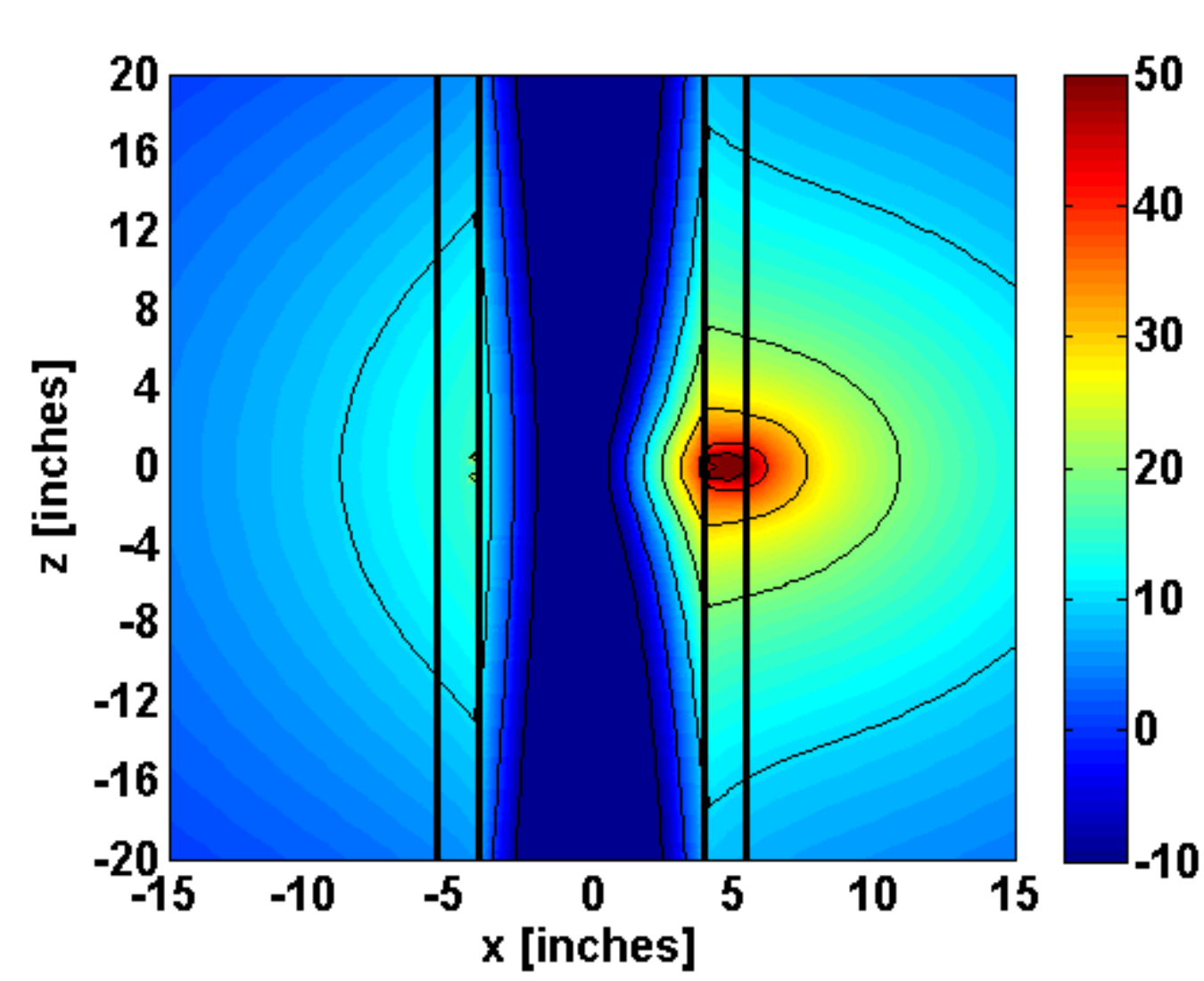}
    }
    \hfill
    \subfloat[\label{ch4.F.case2.y.k1}]{%
      \includegraphics[width=3.0in]{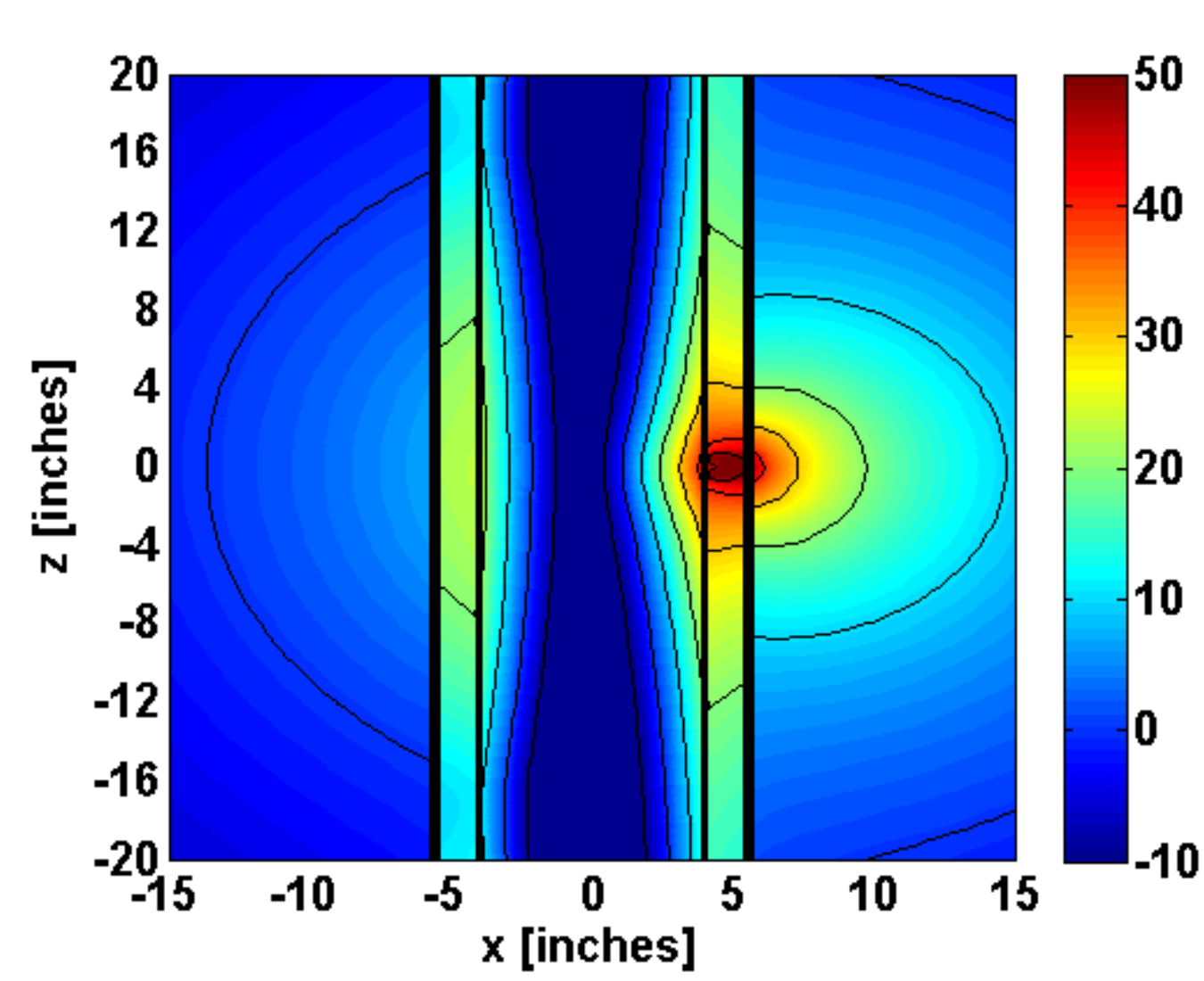}
    }\\
    \subfloat[\label{ch4.F.case1.y.k10}]{%
      \includegraphics[width=3.0in]{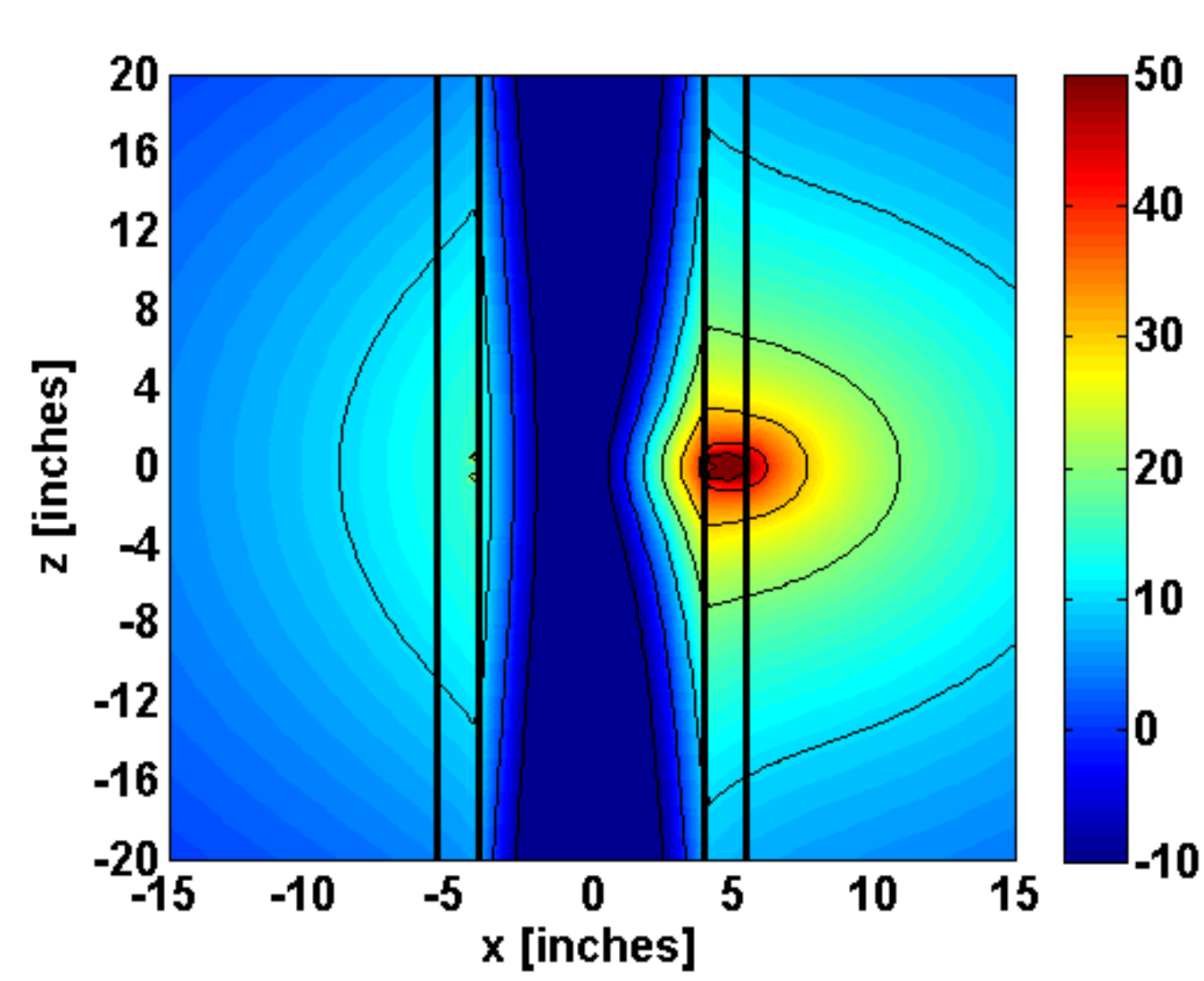}
    }
    \hfill
    \subfloat[\label{ch4.F.case2.y.k10}]{%
      \includegraphics[width=3.0in]{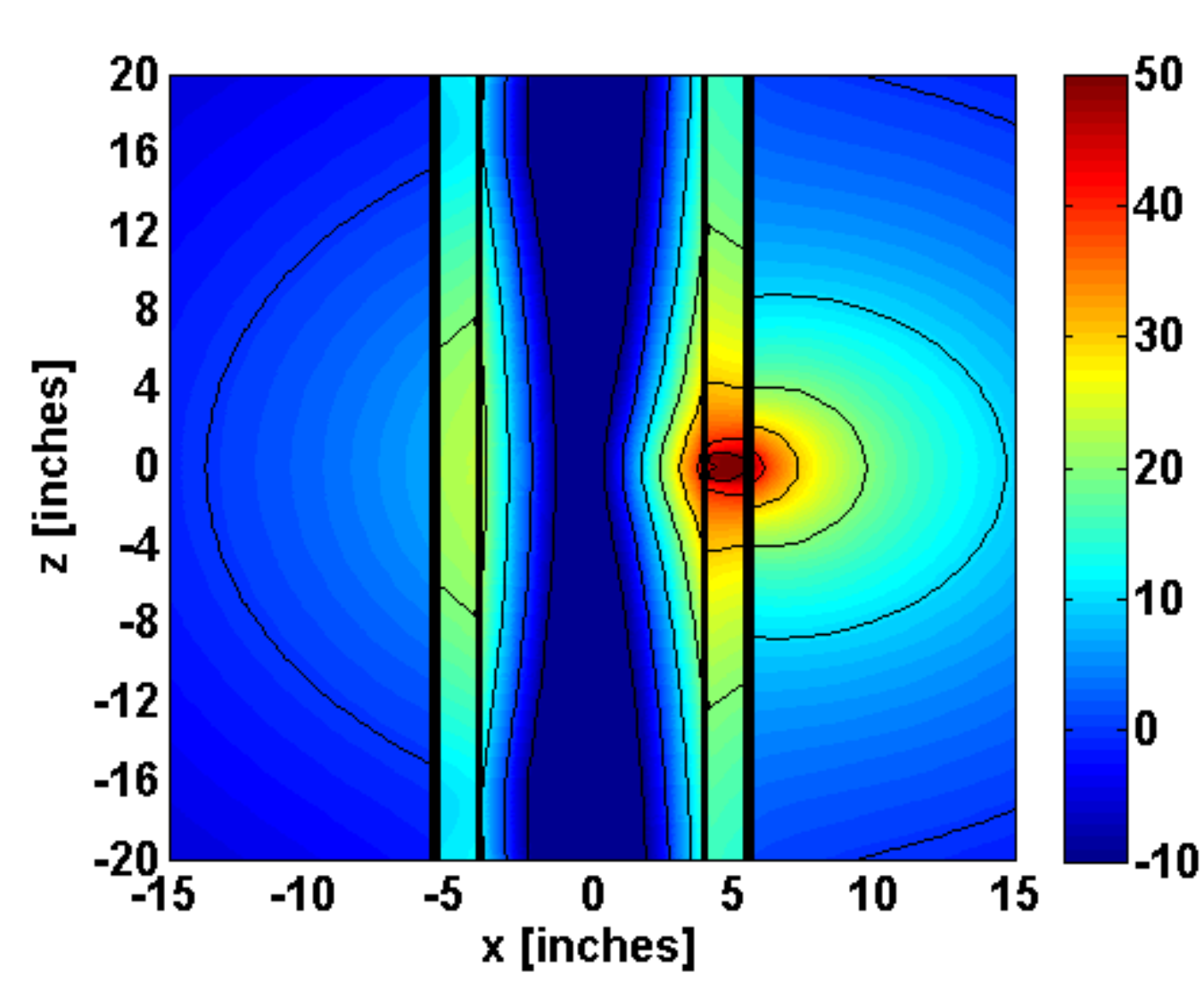}
	}    
    \caption{\textcolor{\Cblue}{Spatial distribution of the magnetic field magnitude on the $y=0^{\prime\prime}$ plane at 36 kHz: (a) Case 1 with $\kappa_\epsilon^2=1$, (b) Case 2 with $\kappa_\epsilon^2=1$, (c) Case 1 with $\kappa_\epsilon^2=10$, and (d) Case 2 with $\kappa_\epsilon^2=10$.}}
    \label{ch4.F.case12.y}
\end{figure}

\textcolor{\Cblue}{
The magnetic field magnitude in the $y=0^{\prime\prime}$ plane is shown in Figure \ref{ch4.F.case12.y}, for Cases 1 and 2. The field is plotted in a decibel scale $10\log_{10}|\mathbf{H}|$ because of the large magnitude variation. In this scale, the small differences in magnitude between $\kappa_\epsilon^2=1$ and $\kappa_\epsilon^2=10$ observed in Tables 3 and 5 are hardly distinguishable. On the other hand, these figures clearly show that Case 1 has less confinement of fields within the source layer than Case 2 due to the presence of the metallic casing in the latter case, as depicted in Figure \ref{ch4.F.case12}. 
}

\section{Conclusion}
\label{sec.5.con}
We provided a robust algorithm for the stable computation of electromagnetic fields in cylindrically stratified media with doubly uniaxial anisotropic layers.  
Range-conditioned integrands, which were originally developed for isotropic media, are extended here for uniaxial media. Associated multiplicative factors used for the stabilization are expressed \textcolor{\Cred}{as} 2$\times$2 matrices in this case. The results show that the formulation is indeed stable and have good error controllability. Illustrative scenarios were included to show applicability of the proposed algorithm to geophysical exploration problems involving borehole sensors in Earth formations with anisotropic responses.
\section*{Acknowledgement}
 We thank Halliburton Energy Services for the permission to publish this work, and Dr. Baris Guner for kindly compiling some of the comparison data.

\appendix
\section{Analytical solution in homogeneous doubly-uniaxial media}
\label{app1}
In this appendix, the closed-form analytical expressions used to obtain for electromagnetic fields in homogeneous \textcolor{\Cred}{and} doubly-uniaxial media are presented. In such media, 
Maxwell's equations with $e^{-\iu\omega t}$ time-dependence write as
\begin{flalign}
\na\times\mathbf{E}(\rr)
	&=\iu\omega\um\cdot\mathbf{H}(\rr), \label{app.1.E.Maxwell.Eq.a}\\
\na\times\mathbf{H}(\rr)
	&=-\iu\omega\ue\cdot\mathbf{E}(\rr) + \mathbf{J}(\rr). \label{app.1.E.Maxwell.Eq.b}
\end{flalign}
Permittivity values are complex-valued so as to include conductivities. For simplicity, it is assumed that the anisotropy ratio for the permeability tensor coincides with that of the complex permittivity tensor, i.e., $\kappa_\epsilon=\kappa_\mu$.
To simplify the derivation, we adopt coordinate stretching techniques. For a general exposition of the coordinate stretching, refer to \cite{Chew94:3D, Teixeira98:General, Teixeira98:Analytical}. To begin with, let us consider modified Maxwell's curl equations with stretched coordinates, i.e.,
\begin{flalign}
\tna\times\tE(\trr)&=\iu\omega\tmu\tH(\trr), \label{app.1.E.Maxwell.E}\\
\tna\times\tH(\trr)&=-\iu\omega\teps\tE(\trr)+\tJ(\trr), \label{app.1.E.Maxwell.H}
\end{flalign}
with a modified nabla operator $\tna$ defined as 
\begin{flalign}
\tna = \hat{x}\frac{\pa}{\pa \widetilde{x}}
	 + \hat{y}\frac{\pa}{\pa \widetilde{y}}
	 + \hat{z}\frac{\pa}{\pa \widetilde{z}}, \label{app.1.E.modified.nabla}
\end{flalign}
where $\widetilde{x}$, $\widetilde{y}$, and $\widetilde{z}$ are stretched coordinates defined such that
\begin{flalign}
u \rightarrow \widetilde{u}=\int_0^u s_u(u')du', \label{app.1.E.stretching.coord}
\end{flalign}
where $s_u$ is the corresponding complex stretching variable, and $u$ stands for $x$, $y$, or $z$. In the above, the fields and sources are non-Maxwellian but $\tmu$ and $\teps$ are scalars, so the medium is isotropic. Using the technique in~\cite{Teixeira98:General}, \eqref{app.1.E.Maxwell.E} and \eqref{app.1.E.Maxwell.H} are rewritten as
\begin{flalign}
\na\times\left(\ttS^{-1}\cdot\tE(\trr)\right)
	&=\iu\omega\tmu\left(\text{det}\ttS\right)^{-1}\ttS\cdot\tH(\trr), \label{app.1.E.Maxwell.E.v3}\\
\na\times\left(\ttS^{-1}\cdot\tH(\trr)\right)
	&=-\iu\omega\teps\left(\text{det}\ttS\right)^{-1}\ttS\cdot\tE(\trr)
	 +\left(\text{det}\ttS\right)^{-1}\ttS\cdot\tJ(\trr). \label{app.1.E.Maxwell.H.v3}
\end{flalign}
where a dyadic $\ttS$ is defined as
\begin{flalign}
\ttS = \hat{x}\hat{x}\left(\frac{1}{s_x}\right)
     + \hat{y}\hat{y}\left(\frac{1}{s_y}\right)
     + \hat{z}\hat{z}\left(\frac{1}{s_z}\right).
\end{flalign}
Using the relations between the stretched fields to unstretched (Maxwellian) fields,
\begin{subequations}
\begin{flalign}
\mathbf{E}(\rr)&=\ttS^{-1}\cdot\tE(\trr), \label{app.1.E.unstr.E}\\
\mathbf{H}(\rr)&=\ttS^{-1}\cdot\tH(\trr), \label{app.1.E.unstr.H}\\
\mathbf{J}(\rr)&=\left(\text{det}\ttS\right)^{-1}\ttS\cdot\tJ(\trr), \label{app.1.E.unstr.J}
\end{flalign}
\end{subequations}
\eqref{app.1.E.Maxwell.E.v3} and \eqref{app.1.E.Maxwell.H.v3} are rearranged as
\begin{flalign}
\na\times\mathbf{E}(\rr)
	&=\iu\omega\left[
					\tmu\left(\text{det}\ttS\right)^{-1}\ttS\cdot\ttS
				\right]\cdot\mathbf{H}(\rr), \label{app.1.E.Maxwell.E.v4}\\
\na\times\mathbf{H}(\rr)
	&=-\iu\omega\left[
					\teps\left(\text{det}\ttS\right)^{-1}\ttS\cdot\ttS\
				\right]\cdot\mathbf{E}(\rr) + \mathbf{J}(\rr). \label{app.1.E.Maxwell.H.v4}
\end{flalign}
These two resulting curl equations can be associated with an effective anisotropic medium and represented as
\begin{flalign}
\na\times\mathbf{E}(\rr)
	&=\iu\omega\um\cdot\mathbf{H}(\rr), \label{app.1.E.Maxwell.E.v5}\\
\na\times\mathbf{H}(\rr)
	&=-\iu\omega\ue\cdot\mathbf{E}(\rr) + \mathbf{J}(\rr), \label{app.1.E.Maxwell.H.v5}
\end{flalign}
which recover the form of the original curl equations \eqref{app.1.E.Maxwell.Eq.a} and \eqref{app.1.E.Maxwell.Eq.b}.
Therefore, electromagnetic fields in an homogeneous and uniaxial media with $\mathbf{E}(\rr)$ and $\mathbf{H}(\rr)$ can be easily obtained from $\tE(\trr)$ and $\tH(\trr)$, which are solutions in isotropic media with coordinate-stretching, by the transformations expressed in \eqref{app.1.E.unstr.E}, \eqref{app.1.E.unstr.H}, and \eqref{app.1.E.unstr.J}.
In order to determine the form of the stretching variables relevant to our problem, let us examine the effective anisotropic medium obtained above. The constitutive tensors have the form
\begin{flalign}
\um
=\left[\tmu\left(\text{det}\ttS\right)^{-1}\ttS\cdot\ttS\right]
=\tmu\ttLa, \label{app.1.E.ani.mu}\\
\ue
=\left[\teps\left(\text{det}\ttS\right)^{-1}\ttS\cdot\ttS\right]
=\teps\ttLa, \label{app.1.E.ani.epsilon}
\end{flalign}
where
\begin{flalign}
\ttLa
=s_x s_y s_z
	\begin{bmatrix}
    s_x^{-2} &        0 &        0 \\
           0 & s_y^{-2} &        0 \\
           0 &        0 & s_z^{-2} \\
    \end{bmatrix}
=
	\begin{bmatrix}
    \frac{s_y s_z}{s_x} &                   0 &                   0 \\
                      0 & \frac{s_y s_z}{s_x} &                   0 \\
                      0 &                   0 & \frac{s_x s_y}{s_z} \\
    \end{bmatrix}. \label{app.1.E.Lambda}
\end{flalign}
Using two conditions on the stretching variables for uniaxial anisotropy, and using the wavenumber expression for the modified Maxwell's equations $\tk=\omega\sqrt{\tmu\teps}$, we can set $s_x=s_y=1$, and $s_z=\kappa$. Consequently,we obtain $\tmu = \frac{\mu_h}{\kappa}$ and $\teps = \frac{\epsilon_h}{\kappa}$.
Next, let us consider the source transformation \eqref{app.1.E.unstr.J}. If the source is a point Hertzian electric dipole like
$ \mathbf{J}(\rr)=Il\hat{\alpha}'\delta(\rr-\rp)$, the coordinate stretching should be carefully treated due to the presence of the Dirac delta function. The stretched current density is expressed as
$ \tJ(\trr)=Il\hat{\widetilde{\alpha}}'\delta(\trr-\trp)$. From the Dirac delta function properties,
\begin{flalign}
\delta(\trr-\trp)
= \frac{1}{s_x s_y s_z}\delta(\rr-\rp),
\end{flalign}
and from \eqref{app.1.E.unstr.J}, 
\begin{flalign}
\hat{\alpha}'\delta(\rr-\rp)
=\left(\text{det}\ttS\right)^{-1}\ttS\cdot\hat{\widetilde{\alpha}}'\delta(\trr-\trp)
=
	\begin{bmatrix}
    s_x^{-1} &        0 &        0 \\
           0 & s_y^{-1} &        0 \\
           0 &        0 & s_z^{-1} \\
    \end{bmatrix}
    \cdot\hat{\widetilde{\alpha}}'\delta(\rr-\rp). \label{app.1.E.delta.comp}  
\end{flalign}
Since $s_x=s_y=1$ and $s_z=\kappa$, we have the source transformation $\hat{\widetilde{\alpha}}'=\ttS^{-1}\cdot\hat{\alpha}'$, and
in homogeneous isotropic media, the Cartesian field components due to the Hertzian electric dipole source can be written as 
\begin{subequations}
\renewcommand{\arraystretch}{1.4}
\begin{flalign}
    \begin{bmatrix}
    \widetilde{E}_x \\ \widetilde{E}_y \\ \widetilde{E}_z
    \end{bmatrix}
&=\frac{\iu Il}{\omega\teps} \; \frac{e^{\iu \tk\tr}}{4\pi\tr} \; \overline{\mathbf{M}}_e \cdot
    \begin{bmatrix}
    \widetilde{\alpha}_{x'} \\ \widetilde{\alpha}_{y'} \\ \widetilde{\alpha}_{z'}
    \end{bmatrix}, \label{app.1.E.Exyz}\\
    \begin{bmatrix}
    \widetilde{H}_x \\ \widetilde{H}_y \\ \widetilde{H}_z
    \end{bmatrix}
&=Il \; \frac{e^{\iu \tk\tr}}{4\pi\tr} \; \overline{\mathbf{M}}_m \cdot
    \begin{bmatrix}
    \widetilde{\alpha}_{x'} \\ \widetilde{\alpha}_{y'} \\ \widetilde{\alpha}_{z'}
    \end{bmatrix}, \label{app.1.E.Hxyz}
\end{flalign} 
\end{subequations}
where
\begin{subequations}
\renewcommand{\arraystretch}{1.4}
\begin{flalign}
\overline{\mathbf{M}}_e&=
    \begin{bmatrix}
    \tk^2+A+BX^2 & BXY        & BXZ \\
    BXY        & \tk^2+A+BY^2 & BYZ \\
    BXZ        & BYZ        & \tk^2+A+BZ^2 \\
    \end{bmatrix}, \label{app.1.E.Me.str}\\
\overline{\mathbf{M}}_m&=
    \begin{bmatrix}
      0 &  AZ & -AY \\
    -AZ &   0 &  AX \\
     AY & -AX &   0 \\
    \end{bmatrix}, \label{app.1.E.Mm.str}\\
A&=\iu \tk/\tr-1/\tr^2, \label{app.1.E.A.str}\\
B&=-\tk^2/\tr^2-3\iu \tk/\tr^3+3/\tr^4, \label{app.1.E.B.str}\\
X&=s_x(x'-x)=x'-x, \label{app.1.E.X.str}\\
Y&=s_y(y'-y)=y'-y, \label{app.1.E.Y.str}\\
Z&=s_z(z'-z)=\kappa(z'-z), \label{app.1.E.Z.str}\\
\tk&=\omega\sqrt{\mu_h \epsilon_h}/\kappa, \label{app.1.E.k.str}\\
\tr&=\left[(x'-x)^2+(y'-y)^2+\kappa^2(z'-z)^2\right]^{1/2}. \label{app.1.E.r.str}
\end{flalign}
\end{subequations}
Applying field transformations, \eqref{app.1.E.unstr.E} and \eqref{app.1.E.unstr.H}, and source transformation $\hat{\widetilde{\alpha}}'=\ttS^{-1}\cdot\hat{\alpha}'$, we obtain
\begin{subequations}
\renewcommand{\arraystretch}{1.4}
\begin{flalign}
    \begin{bmatrix}
    E_x \\ E_y \\ E_z
    \end{bmatrix}
&=\frac{\iu Il}{\omega\teps} \; \frac{e^{\iu \tk\tr}}{4\pi\tr} \; 
	\ttS^{-1} \cdot \overline{\mathbf{M}}_e \cdot \ttS^{-1} \cdot
    \begin{bmatrix}
    \alpha_{x'} \\ \alpha_{y'} \\ \alpha_{z'}
    \end{bmatrix}, \label{app.1.E.Exyz.unstr}\\
    \begin{bmatrix}
    H_x \\ H_y \\ H_z
    \end{bmatrix}
&=Il \; \frac{e^{\iu \tk\tr}}{4\pi\tr} \; 
	\ttS^{-1} \cdot \overline{\mathbf{M}}_m \cdot \ttS^{-1} \cdot
    \begin{bmatrix}
    \alpha_{x'} \\ \alpha_{y'} \\ \alpha_{z'}
    \end{bmatrix}. \label{app.1.E.Hxyz.unstr}
\end{flalign} 
\end{subequations}
Finally, applying the coordinate transformations from Cartesian to cylindrical coordinates, we obtain
\begin{subequations}
\renewcommand{\arraystretch}{1.4}
\begin{flalign}
    \begin{bmatrix}
    E_\rho \\ E_\phi \\ E_z
    \end{bmatrix}
&=\frac{\iu Il}{\omega\teps} \; \frac{e^{\iu \tk\tr}}{4\pi\tr} \; 
	\overline{\mathbf{T}}_1 \cdot \ttS^{-1} \cdot \overline{\mathbf{M}}_e \cdot \ttS^{-1}
	\cdot \overline{\mathbf{T}}_2 \cdot
    \begin{bmatrix}
    \alpha_{\rho'} \\ \alpha_{\phi'} \\ \alpha_{z'}
    \end{bmatrix}, \label{app.1.E.Erpz}\\
    \begin{bmatrix}
    H_\rho \\ H_\phi \\ H_z
    \end{bmatrix}
&=Il \; \frac{e^{\iu \tk\tr}}{4\pi\tr} \; 
	\overline{\mathbf{T}}_1 \cdot \ttS^{-1} \cdot \overline{\mathbf{M}}_m \cdot \ttS^{-1}
	\cdot \overline{\mathbf{T}}_2 \cdot
    \begin{bmatrix}
    \alpha_{\rho'} \\ \alpha_{\phi'} \\ \alpha_{z'}
    \end{bmatrix}, \label{app.1.E.Hrpz}
\end{flalign} 
\end{subequations}
where
\begin{subequations}
\renewcommand{\arraystretch}{1.4}
\begin{flalign}
\overline{\mathbf{T}}_1&=
    \begin{bmatrix}
     \cos\phi & \sin\phi & 0 \\
    -\sin\phi & \cos\phi & 0 \\
            0 &        0 & 1 \\
    \end{bmatrix}, \label{app.1.E.T1}\\
\overline{\mathbf{T}}_2&=
    \begin{bmatrix}
    \cos\phi' & -\sin\phi' & 0 \\
    \sin\phi' &  \cos\phi' & 0 \\
            0 &          0 & 1 \\
    \end{bmatrix}. \label{app.1.E.T2}
\end{flalign}
\end{subequations}




\bibliographystyle{model1-num-names}
\bibliography{JCPrefs}







\end{document}